\documentclass[a4paper,reqno,10pt,oneside]{amsart}
\usepackage{amsmath,geometry}
\usepackage{graphicx}
\usepackage{mathrsfs}
\usepackage{amssymb,amsmath, amsfonts, amsthm}
\usepackage{latexsym}
\usepackage{color} 
\usepackage[dvips]{epsfig}
\usepackage{graphicx}

\allowdisplaybreaks[1]

\numberwithin{equation}{section}

\renewcommand{\O}{\operatorname{O}}

\renewcommand{\(}{\left(}
\renewcommand{\)}{\right)}
\renewcommand{\[}{\left[}
\renewcommand{\]}{\right]}

\newtheorem{theorem}{Theorem}[section]
\newtheorem{proposition}[theorem]{Proposition}
\newtheorem{lemma}[theorem]{Lemma}
\newtheorem{remark}[theorem]{Remark}

\newcommand{\zn}{H^1(\mathbb R^N)}

\renewcommand{\d }{\delta }

\renewcommand{\O}{\Omega}

\newcommand{\U}{\mathcal{U}}
\newcommand{\p}{\mathcal{P}}

\newcommand{\beq}{\begin{equation}}
\newcommand{\eeq}{\end{equation}}
\newcommand{\beqs}{\begin{equation*}}
\newcommand{\eeqs}{\end{equation*}}
\newcommand{\beqn}{\begin{eqnarray}}
\newcommand{\eeqn}{\end{eqnarray}}
\newcommand{\beqns}{\begin{eqnarray*}}
\newcommand{\eeqns}{\end{eqnarray*}}
\newcommand{\bdoc}{\begin{document}}
\newcommand{\edoc}{\end{document}}
\newcommand{\be}{\begin{enumerate}}
\newcommand{\ee}{\end{enumerate}}
\newcommand{\bdescr}{\begin{description}}
\newcommand{\edescr}{\end{description}}
\newcommand{\ba}{\begin{array}}
\newcommand{\ea}{\end{array}}
\newcommand{\intR}{\int_{\mathbb R^N}}
\newcommand{\R}{\mathbb R^N}
\newcommand{\e}{\epsilon}
 \renewcommand{\(}{\left(}
\renewcommand{\)}{\right)}
\renewcommand{\[}{\left[}
\renewcommand{\]}{\right]}
\newenvironment{Proof}{\noindent{\bf Proof}}{\hfill$\Box$\\[2mm]}

\begin{document}
\title[Sign-changing tower of bubbles for the Brezis-Nirenberg problem]{Sign-changing tower of bubbles for the Brezis-Nirenberg problem}

\author{Alessandro Iacopetti}
\address[Alessandro Iacopetti]{Dipartimento di Matematica e Fisica, Universit\'a degli Studi di Roma Tre, L.go S. Leonardo Murialdo 1, 00146 Roma, Italy}
\email{iacopetti@mat.uniroma3.it}

\author{Giusi Vaira}
\address[Giusi Vaira] {Dipartimento di Matematica ``G. Castelnuovo", Universit\`{a} di Roma ``La Sapienza", Piazzale A. Moro 1, 00161 Roma, Italy}
\email{vaira@mat.uniroma1.it}
\thanks{Research partially supported by MIUR-PRIN project-201274FYK7\underline\ 005.}
\begin{abstract}
In this paper, we prove that the Brezis-Nirenberg problem
$$-\Delta u =|u|^{p-1}u+\e u\qquad \mbox{in}\,\, \Omega,\quad u=0\,\, \mbox{on}\,\,\ \partial\Omega,$$
where $\Omega$ is a symmetric bounded smooth domain in $\R$, $N\geq 7$ and $p=\frac{N+2}{N-2}$, has a solution with the shape of a tower of two bubbles with alternate signs, centered at the center of symmetry of the domain, for all $\e>0$ sufficiently small.
\end{abstract}

\maketitle

\section{Introduction and statement of the main result}\label{intro}

In this paper we are interested in the construction of solutions to the following problem
\begin{equation}\label{BN}
\left\{
\begin{array}{lr}
-\Delta u = |u|^{p-1}u+\e u \qquad \mbox{in}\,\, \Omega\\
u=0,\qquad\qquad\qquad\qquad\mbox{on}\,\, \partial\Omega
\end{array}
\right.
\end{equation}
where $\Omega$ is a bounded smooth domain of $\mathbb R^N$ with $N\geq 7$, $\e$ is supposed to be small and positive while $p+1=\frac{2N}{N-2}$ is the critical Sobolev exponent for the embedding of $H^1_0(\Omega)$ into $L^{p+1}(\Omega)$.\\\\
The pioneering paper on equation \eqref{BN} was written by Brezis and Nirenberg \cite{Brezis} in 1983 where the authors showed that for $N\geq 4$ and $\e\in (0, \lambda_1)$, the problem \eqref{BN} has at least one positive solution where $\lambda_1$ denotes the first eigenvalue of $-\Delta$ on $\Omega$.\\ In the case $N=3$, a similar result was proved in \cite{Brezis} but only for $\e\in (\lambda^*, \lambda_1)$ with $\lambda^*=\lambda^*(\O)>0$. Moreover by using a version of the Pohozaev Identity the authors showed that $\lambda^*(\O)=\frac 14 \lambda_1$ if $\O$ is a ball and that no positive solutions exist for $\e\in (0, \frac 14 \lambda_1)$. \\ Note that, by using again Pohozaev Identity, it is easy to check that problem \eqref{BN} has no nontrivial solutions when $\e\leq 0$ and $\Omega$ is star-shaped.\\ Since then, there has been a considerable number of papers on problem \eqref{BN}.\\ We briefly recall some of the main ones.\\ Han, in \cite{Han}, proved that the solution found by Brezis and Nirenberg blows-up at a critical point of the Robin's function as $\e$ goes to zero. Conversely, Rey in \cite{Rey} and in \cite{Rey1} proved that any $C^1-$ stable critical point of the Robin's function generates a family of positive solutions which blows-up at this point as $\e$ goes to zero.\\ After the work of Brezis and Nirenberg, Capozzi, Fortunato and Palmieri \cite{Capozzi} showed that for $N=4$, $\e>0$ and $\e\not\in \sigma(-\Delta)$ (the spectrum of $-\Delta$) problem \eqref{BN} has a nontrivial solution. The same holds if $N\geq 5$ for all $\e>0$ (see also \cite{Grunau}).\\
\\ The first multiplicity result was obtained by Cerami, Fortunato and Struwe in \cite{Cerami}, in which they proved that the number of nontrivial solutions of \eqref{BN}, for $N\geq3$, is bounded below by the number of eigenvalues of $(-\Delta, \Omega)$ belonging to $(\e, \e+ S|\Omega|^{-2/N})$, where $S$ is the best constant for the Sobolev embedding $D^{1, 2}(\R)$ into $L^{p+1}(\R)$ and $|\Omega|$ is the Lebesgue measure of $\Omega$.\\
Moreover, if $N\geq 4$, then for any $\e>0$ and for a suitable class of symmetric domain $\Omega$, problem \eqref{BN} has infinitely many solutions of arbitrarily large energy (see Fortunato and Jannelli \cite{Fortunato}). \\ If $N\geq 7$ and $\Omega$ is a ball, then for each $\e>0$, problem \eqref{BN} has infinitely many sign-changing radial solutions (see Solimini \cite{Solimini}).\\ In the papers \cite{Fortunato, Solimini}, the radial symmetry of the domain plays an essential role, therefore their methods do not work for general domains.\\
Concerning sign-changing solutions, Cerami, Solimini and Struwe  showed in \cite{Cerami2} that if $N\geq 6$ and $\e\in (0, \lambda_1)$, problem \eqref{BN} has a pair of least energy sign-changing solution. In the same paper the authors studied the multiplicity of nodal solutions proving the existence of infinitely many radial solutions when $\Omega$ is a ball centered at the origin. \\ On the other side, for $3\leq N\leq 6$ and when $\Omega$ is a ball, it can be proved that there is a $\lambda^*>0$ such that \eqref{BN} has no sign-changing radial solutions for $\e\in (0, \lambda^*)$ (see Atkinson, Brezis and Peletier \cite{Atkinson}).\\ Moreover, Devillanova and Solimini in \cite{Devillanova} showed that, if $N\geq 7$ and $\O$ is an open regular subset of $\R$, problem \eqref{BN} has infinitely many solutions for each $\e>0$.\\ For low dimensions, namely $N=4, 5, 6$ and in an open regular subset of $\R$, in \cite{Devillanova1}, Devillanova and Solimini proved the existence of at least $N+1$ pairs of solutions provided $\e$ is small enough. In \cite{Clapp}, Clapp and Weth extended this last result to all $\e>0$.\\ Neither in \cite{Devillanova, Devillanova1} nor in \cite{Clapp} there is information on the kind of  sign-changing solutions obtained.\\
Recently, in \cite{Schechter}, Schechter and Wenming Zou showed that in any bounded and smooth domain, for $N\geq 7$ and for each fixed $\e>0$, problem \eqref{BN} has infinitely many sign changing solutions. 

Concerning the profile of sign-changing solutions some results have been obtained in \cite{Ben}, \cite{Ben1} for low energy solutions, namely solutions $u_\e$ such that $\displaystyle{\int_{\Omega}|\nabla u_\e|^2\, dx \rightarrow 2 S^{\frac{N}{2}}}$, as $\e\rightarrow 0$, $S$ being the Sobolev constant for the embedding of $H^1_0(\O)$ into $L^{p+1}(\O)$. More precisely in \cite{Ben} it is proved that for $N=3$ these solutions concentrate and blow-up in two different points of $\O$, as $\e\rightarrow 0$, and have the asymptotic profile of two separate bubbles. A similar result is proved in \cite{Ben1} for $N\geq 4$ but assuming that the blow-up rate of the positive and negative part of $u_\e$ is the same.\\ Existence of nodal solutions with two nodal regions concentrating in two different points of the domain $\O$ as $\e\rightarrow 0$ has been obtained in \cite{Castro}, \cite{Micheletti} and \cite{Bartsch}. So none of these solutions look like tower of bubbles, i.e. superposition of two bubbles with opposite sign concentrating at the same point, as $\e\rightarrow 0$. Such a type of solutions is shown to exist for other semilinear problems like the almost critical Lane-Emden problem (see \cite{Ben2}, \cite{Pistoia}, \cite{Musso1}) but not, to our knowledge, for the Brezis-Nirenberg problem with the exception of the case of the ball. If $\O$ is a ball, and $N\geq 7$, in a recent paper \cite{Iacopetti} the asymptotic behaviour as $\e\rightarrow 0$ of the least energy nodal radial solution $v_\e$ is analysed and among other things, it is shown that the positive and negative part of $v_\e$ concentrate at the origin. Moreover they have the asymptotic profile of a positive and negative solution of the critical problem in $\R$ and the concentration speeds are different.\\ Hence \cite{Iacopetti} provides the first example of bubble of towers for the Brezis-Nirenberg problem.\\ Then the natural question is whether these kind of solutions exist in bounded domains other than the ball.\\ In the present paper we answer positively this question constructing a sign-changing solution of \eqref{BN} in any bounded domain symmetric with respect to $N$ orthogonal axis.

We next state our result.

\begin{theorem}\label{principale}
Let $N\geq 7$. There exists $\e_0>0$ such that for any $\e\in(0, \e_0)$ there exist positive numbers $d_{j\e}$, $j=1, 2$ and a  solution $u_\e$ of problem \eqref{BN} of the form
\begin{equation}\label{soluzioneforma}
u_\e(x)=\alpha_N\left[\left(\frac{d_{1\e} \e^{\frac{1}{N-4}}}{d_{1\e}^2 \e^{\frac{2}{N-4}}+|x|^2}\right)^{\frac{N-2}{2}}-\left(\frac{d_{2\e} \e^{\frac{3N-10}{(N-4)(N-6)}}}{d_{2\e}^2\e^{2\frac{3N-10}{(N-4)(N-6)}}+|x|^2}\right)^{\frac{N-2}{2}}\right]+\Phi_\e,
\end{equation}
where $\alpha_N:=[N(N-2)]^{\frac{N-2}{4}}$, $d_{j\e}\rightarrow \bar
d_j>0$, as $\e\rightarrow 0$, $\Phi_\e\rightarrow 0$ in
$H_0^1(\Omega)$, as $\e\rightarrow 0$. Moreover $u_\e$ is even with
respect to the variables $x_1, \ldots, x_N$.
\end{theorem}

We remark that the assumption $N\geq 7$ in our proof is crucial.
We believe that it is possible to extend our result to a general domain $\O$ with some suitable modifications.\\\\
In the case the remainder term converges to zero also in $L^\infty_{loc}(\Omega)$, then,
the asymptotic expansion and some energy estimates derived in the course of the proof allow to draw interesting consequences concerning the number and shape of the nodal domains of the solution $u_\e$.

\begin{theorem}\label{principale1}
Let $N\geq 7$ and assume that the remainder term $\Phi_\e$, appearing in Theorem \ref{principale}, is such that $\Phi_\e \rightarrow 0$ uniformly in compact subsets of $\Omega$. Then, there exists $\e_0>0$ such that for any $\e\in(0,\e_0)$, the solution $u_\e$ constructed in Theorem \ref{principale} has precisely two nodal domains $\O_\e^1$, $\O_\e^2$ such that $\O_\e^1$ contains the sphere $\mathcal S_\e^1:=\left\{x\in\R\,\,\,:\,\,\, |x|=\e^{\frac{1}{N-4}}\right\}$, $\O_\e^2$ contains the sphere $\mathcal S_\e^2:=\left\{x\in\R\,\,\:\,\,\ |x|=\e^{\frac{3N-10}{(N-4)(N-6)}}\right\}$ and $u_\e>0$ on $\O_\e^1$ and $u_\e<0$ on $\O_\e^2$.\\ Consequently, $0\in\O_\e^2$ and $\O_\e^1$ is the only nodal domain of $u_\e$ which touches $\partial\O$.
\end{theorem}
\begin{remark}
Under the assumptions of Theorem \ref{principale1} it follows that the sign-changing tower of bubble $u_\e$ constructed in Theorem \ref{principale} has two nodal domains and its nodal set does not touch $\partial\O$. By this we mean that, denoting by $$Z_\e:=\left\{x\in\O\,\,:\,\, u_\e(x)=0\right\}$$ the nodal set of $u_\e$ then $\overline Z_\e \cap \partial\O=\emptyset.$
\end{remark}
The proof of Theorem \ref{principale} is based on the Lyapunov-Schmidt reduction.\\
To describe the procedure and explain the difficulties which arise when looking for bubble towers of the Brezis-Nirenberg problem, we introduce the functions
\begin{equation}\label{Udelta}
\U_{\delta}(x)=\alpha_N \frac{\delta^{\frac{N-2}{2}}}{\left(\delta^2+|x|^2\right)^{\frac{N-2}{2}}},\qquad \delta>0
\end{equation}

with $\alpha_N:=[N(N-2)]^{\frac{N-2}{4}}$. Is is well known (see \cite{Aubin}, \cite{Caffarelli}, \cite{Talenti}) that \eqref{Udelta} are the only radial solutions of the equation
\begin{equation}\label{pb0}
-\Delta u= u^p\qquad \mbox{in}\,\, \R.
\end{equation}
We define $\varphi_\delta$ to be the unique solution to the problem
\begin{equation}\label{pbvarphi}
\left\{
\begin{array}{lr}
\Delta\varphi_\delta=0\qquad \mbox{in}\,\, \Omega\\
\varphi_\delta=\U_{\delta}\qquad\quad\mbox{on}\,\, \partial\Omega,
\end{array}
\right.
\end{equation}
and let
\begin{equation}\label{proiezione}
\p \U_\delta:=\U_\delta-\varphi_\delta
\end{equation}
be the projection of $\U_\delta$ onto $H^1_0(\Omega)$, i.e.
\begin{equation}\label{pbproiezione}
\left\{
\begin{array}{lr}
-\Delta\p\U_\delta=\U_\delta^p\qquad \mbox{in}\,\,\ \Omega\\
\p\U_\delta=0\qquad\qquad\mbox{on}\,\,\,\partial\Omega.
\end{array}
\right.
\end{equation}
Finally, let $G(x, y)$ be the Green's function associated to $-\Delta$ with Dirichlet boundary conditions and $H(x, y)$ be its regular part, namely
 $$H(x, y)=\frac{1}{|x-y|^{N-2}}-\frac{1}{\gamma_N}G(x, y), \qquad \forall\,\, x, y \in \O, \qquad \mbox{with}\,\, \gamma_N=\frac{1}{N(N-2)\omega_N},$$ where $\omega_N$ is the volume of the unit ball in $\R$.\\ The function $\tau(x):=H(x, x)$, $x\in \O$ is called {\it Robin's function}. \\\\
It is well-known that the following expansions holds (see \cite{Rey})
\begin{equation}\label{expvarphi}
\varphi_\delta(x)=\alpha_N \delta^{\frac{N-2}{2}}H(0, x)+O(\delta^{\frac{N+2}{2}})\qquad \mbox{as}\,\, \delta\rightarrow 0.
\end{equation}
Moreover, from elliptic estimates it follows that
\begin{equation}\label{stimaproiezione}
0<\varphi_\delta(x)<c\delta^{\frac{N-2}{2}}, \qquad \mbox{in}\,\, \Omega
\end{equation}
and
\begin{equation}\label{stimavarphi}
|\varphi_\delta|_{q, \O}\leq C \delta^{\frac{N-2}{2}},
\qquad q\in \left(\frac{p+1}{2}, p+1\right]
\end{equation}
see for instance \cite{Vaira} and references therein.
\\\\ We look for an approximate solution to problem \eqref{BN} which is a superposition of two standard bubbles with two different scaling parameters, namely we take $\delta_1>\delta_2$ and we look for a solution to \eqref{BN} of the form
\begin{equation}\label{sol}
u_\e(x)=\p\U_{\delta_1}-\p\U_{\delta_2}+\Phi_\e(x)
\end{equation}
where the remainder term $\Phi_\e$ is a small function which is even with respect to the variables $x_1, \ldots, x_N$.\\\\
The Lyapunov-Schmidt reduction allows us to reduce the problem of finding blowing-up solutions to \eqref{BN} to the problem of finding critical points of a functional (the reduced energy) which depends only on the concentration parameters.\\ As announced before in our case some difficulties arise which need some modification of the standard procedure to be overcome.\\ Indeed, first we remark that the solutions of problem \eqref{BN} are the critical points of the functional $J_\e: H^1_0(\O)\rightarrow \mathbb R$ defined as
\begin{equation}\label{funzionale}
J_\e(u)=\frac 12 \int_\O |\nabla u|^2\, dx-\frac{1}{p+1}\int_\O |u|^{p+1}\, dx-\frac{\e}{2}\int_\O u^2\, dx,\qquad u\in H^1_0(\O).
\end{equation}
If we apply directly the reduction method looking for a solution of the form \eqref{sol} we get that the remainder term is such that $$\|\Phi_\e\|=O\left(\e^{\frac{N-2}{N-4}+\sigma}\right)\qquad \sigma>0$$ where $\|\cdot\|$ denotes the $H^1_0(\O)$-norm, and that the reduced energy
\begin{equation*}
\mbox{{\it Reduced Energy}}\sim  J_\e (\p\U_{\delta_1}-\p\U_{\delta_2})= C+C_1\tau(0)\delta_1^{N-2}-C_2\e\delta_1^2+C_3\left(\frac{\delta_2}{\delta_1}\right)^{\frac{N-2}{2}}+H.O.T.\\
 \end{equation*}
 where $C, C_i$ are some known positive constants.\\
Since $\delta_1,\delta_2$ are proper power of $\e$ of the form $\delta_j=\e^{\gamma_j}d_j$, $d_j>0$ , after some easy computations, in order to find a critical point of the reduced energy, we get that
\begin{equation*}
\mbox{{\it Reduced Energy}} \sim C+\e^{\frac{N-2}{N-4}}\left[C_1\tau(0)d_1^{N-2}-C_2d_1^2+C_3\left(\frac{d_2}{d_1}\right)^{\frac{N-2}{2}}\right]+o(\e^{\frac{N-2}{N-4}})
\end{equation*}
with $$\gamma_1:=\frac{1}{N-4};\qquad \gamma_2:=\frac{3}{N-4}.$$
 However the function $$\Psi(d_1, d_2)=C_1\tau(0)d_1^{N-2}-C_2d_1^2+C_3\left(\frac{d_2}{d_1}\right)^{\frac{N-2}{2}}$$ has a critical point in $d_1$ but not in $d_2$ and hence in this way we cannot find a solution of our problem. \\\\ Hence we use a new idea. We split the remainder term $\Phi_\e$ in two parts: $$\Phi_\e(x)=\phi_{1,\e}(x)+\phi_{2,\e}(x)$$ such that $$\|\phi_{2,\e}\|=o(\|\phi_{1,\e}\|), \ \ \hbox{as} \ \e\to 0.$$

Usually, the remainder term $\Phi_\e$, solution of the auxiliary equation, is found with a fixed point argument. Here we have to use the Contraction Mapping Theorem twice, since we split the auxiliary equation in a system of two equations. The first one depends only on $\phi_1$ while the second one depends on both $\phi_1, \phi_2$. So we solve the first equation in $\phi_1$ and then the second one finding $\phi_2$. Then we obtain the remainder term $\Phi_\e$ which consists of two terms of different orders. Then we study the finite-dimensional problem, namely the reduced energy that consists of two functions of different orders. The lower term depends only on $d_1$ while the term of higher order depends on $d_1, d_2$. At the end we look for a critical point of this new type of reduced energy. We believe that our strategy can be used also in other contexts. \\\\
The outline of the paper is the following: in Section \ref{setting} we explain the setting of the problem.
In Section \ref{ausiliaria} we look for the remainder term $\Phi_\e$ in a suitable space. In Section \ref{biforcazione} we study the reduced energy and finally Theorem \ref{principale} and Theorem \ref{principale1} are proved in Section \ref{teorema}.\\\\
{\bf Acknowledgments:} The authors wish to thank F. Pacella for proposing the problem and for useful suggestions.

\section{Setting of the problem}\label{setting}
In what follows we let $$(u, v):=\int_\O \nabla u \cdot \nabla v\, dx,\qquad \|u\|:=\(\int_\O |\nabla u|^2\, dx\)^{\frac 12}$$ as the inner product in $H^1_0(\O)$ and its corresponding norm while we denote by
$(\cdot, \cdot)_{\zn}$ and by $\|\cdot\|_{\zn}$ the scalar product and the standard norm in $\zn$.
Moreover we denote by $$|u|_{r}:=\(\int_\O|u|^r\, dx\)^{\frac{1}{r}}$$ the $L^r(\O)$-standard norm for any $r\in [1, +\infty)$ and by $|u|_{r, A}$ the $L^r(A)$- norm when $A$ is a proper subset of $\Omega$ or when $A=\R$.\\\\

From now on we assume that $\O$ is a bounded open set with smooth boundary of $\R$, symmetric with respect to $x_1, \ldots, x_N$ and which contains the origin. Moreover we assume that $N\geq 7$.\\\\ We define then $$H_{sim}:=\left\{u\in H^1_0(\O)\,\,:\,\, u \,\, \mbox{is symmetric with respect to each variable}\,\, x_k,\ k=1,\ldots,N\right\},$$ and for $q\in [1, +\infty)$ $$L^q_{sim}:=\left\{u\in L^q(\O)\,\,:\,\,\ u \,\, \mbox{is symmetric with respect to each variable}\,\, x_k,\ k=1,\ldots,N\right\}.$$
Let $i^*:L^{\frac{2N}{N+2}}_{sim}\rightarrow H_{sim}$ be the adjoint operator of the embedding $i:H_{sim}(\O)\rightarrow L^{\frac{2N}{N-2}}_{sim}$, namely if $v\in L^{\frac{2N}{N+2}}_{sim}$ then $u=i^*(v)$ in $H_{sim}$ is the unique solution of the equation $$-\Delta u=v \qquad \mbox{in}\,\, \O\qquad v=0\qquad \mbox{on}\,\, \partial\O.$$
By the continuity of $i$ it follows that
\beq\label{stimaistar}
\|i^*(v)\|\leq C |v|_{\frac{2N}{N+2}} \qquad \forall v\in L^{\frac{2N}{N+2}}_{sim}
\eeq
for some positive constant $C$ which depends only on $N$.\\
Hence we can rewrite problem \eqref{BN} in the following way

\beq\label{pbriscritto}
\left\{
\ba{lr}
u=i^*\[f(u)+\e u\]\\
u\in H_{sim}
\ea
\right.
\eeq
where $f(s)=|s|^{p-1}s$, $p=\frac{N+2}{N-2}$.\\\\ We next describe the shape of the solution we are looking for.\\\\ Let $\delta_j=\delta_j(\epsilon)$, for $j=1, 2$ be positive parameters defined as proper powers of $\e$, multiplied by a suitable positive constant to be determined later, namely
\begin{equation}\label{deltaj}
\delta_j=\e^{\alpha_j}d_j\qquad \mbox{with}\,\,\ d_j>0
\end{equation}
and $\alpha_1:=\frac{1}{N-4}$; $\alpha_2:=\frac{3N-10}{(N-4)(N-6)}$.\\ Fixed a small $\eta>0$ we impose that the parameters $d_j$ will satisfy
\begin{equation}\label{limdj}
\eta<d_j<\frac{1}{\eta}\qquad \mbox{for}\,\,\, j=1, 2.
\end{equation}
Hence, it is immediate to see that
$$\frac{\delta_2}{\delta_1}=\e^{\frac{2(N-2)}{(N-4)(N-6)}}\frac{d_2}{d_1}\rightarrow 0\qquad \mbox{as}\,\, \e\rightarrow 0.$$
We construct solutions to problem \eqref{BN}, as predicted by Theorem \ref{principale}, which are superpositions of copies of the standard bubble defined in \eqref{Udelta} with alternating signs, properly modified (namely we consider the projection of the original bubble into $H^1_0(\O)$), centered at the origin which is the center of symmetry of $\O$ with parameters of concentrations $\delta_j$. Such an object has the shape of a tower of two bubbles.\\\\ 

Hence the solution to problem \eqref{BN} will be of the form
\begin{equation}\label{formasoluzione}
u_\e(x)=V_\e(x)+\Phi_\e(x)
\end{equation}
where
\begin{equation}\label{V}
V_\e(x):=\p\U_{\delta_1}(x)-\p\U_{\delta_2}(x).
\end{equation}
The term $\Phi_\e$ has to be thought as a remainder term of lower order, which has to be described accurately.\\
Let $Z_j$ the following functions
\begin{equation}\label{Zdelta}
Z_j(x):=\partial_{\delta_j} \U_{\delta_j}(x)=\alpha_N\frac{N-2}{2}\delta_j^{\frac{N-4}{2}}\frac{|x|^2-\delta_j^2}{\left(\delta_j^2+|x|^2\right)^{\frac N 2}}, \ \ j=1,2.
\end{equation}
We remark that the functions $Z_j$  solve the problem (see \cite{Bianchi})
\begin{equation}\label{linearizzatopb0}
-\Delta z= p|\U_\delta|^{p-1}z,\qquad \mbox{in}\,\, \R.
\end{equation}
Let $\p Z_j$ the projection of $Z_j$ onto $H^1_0(\O)$. Elliptic estimates give

\begin{equation}\label{stimaproiezionederivata}
\p Z_j(x)=Z_j(x)-\alpha_N\frac{N-2}{2}\delta_j^{\frac{N-4}{2}}H(0, x)+O(\delta_j^{\frac{N}{2}}), \ \ j=1,2,
\end{equation}
uniformly in $\Omega$. \\\\ Let us consider $$\mathcal{K}_1:={\rm span} \left\{\p Z_1\right\}\subset H_{sim}; \qquad \mathcal{K}:={\rm span} \left\{\p Z_j\,:\, j=1,2\right\}\subset H_{sim}$$ and $$\mathcal{K}_1^\bot:=\left\{\phi\in H_{sim}\,:\,\, \langle \phi, \p Z_1\rangle=0\right\}; \qquad \mathcal{K}^\bot:=\left\{\phi\in H_{sim}\,:\,\, \langle \phi, \p Z_j\rangle=0,\  j=1, 2\right\}.$$
Let $\Pi_1:H_{sim}\rightarrow \mathcal{K}_1$, $\Pi:H_{sim}\rightarrow \mathcal{K}$ and $\Pi_1^\bot: H_{sim}\rightarrow \mathcal{K}_1^\bot$, $\Pi^\bot: H_{sim}\rightarrow \mathcal{K}^\bot$ be the projections onto $\mathcal {K}_1$,  $\mathcal K$ and $\mathcal {K}_1^\bot$, $\mathcal K^\bot$, respectively.\\ In order to solve problem \eqref{BN} we will solve the couple of equations
\begin{equation}\label{aux}
\Pi^\bot\left\{V_\e+\Phi_\e-i^*\left[f(V_\e+\Phi_\e)+\e(V_\e+\Phi_\e)\right]\right\}=0
\end{equation}
\begin{equation}\label{bif}
\Pi\left\{V_\e+\Phi_\e-i^*\left[f(V_\e+\Phi_\e)+\e(V_\e+\Phi_\e)\right]\right\}=0.
\end{equation}
For any $(d_1,d_2)$ satisfying condition \eqref{limdj}, we solve first the equation \eqref{aux} in $\Phi_\e\in \mathcal{K}^\bot$ which is the lower order term in the description of the ansatz.\\\\
We start with solving the auxiliary equation \eqref{aux}.
As anticipated in the introduction, we split the remainder term as $$\Phi_\e=\phi_{1,\e}+\phi_{2,\e}$$ with $$\|\phi_{2,\e}\|=o(\|\phi_{1,\e}\|), \ \ \hbox{as} \ \e\to 0.$$
%
In order to find $\phi_{1,\e}$ and $\phi_{2,\e}$ we solve the following system of equations

\begin{equation}\label{sistemaaux}
\left\{
\begin{array}{lr}
\mathcal{R}_1+\mathcal{L}_1(\phi_1)+\mathcal{N}_1(\phi_1)=0,\\\\
\mathcal{R}_2+\mathcal{L}_2(\phi_2)+\mathcal{N}_2(\phi_1, \phi_2)=0,
\end{array}
\right.
\end{equation}
where

\begin{equation}\label{defelementi2}
\mathcal R_1:=\Pi^\bot_1\left\{\p\U_{\delta_1}-i^*\left[f(\p\U_{\delta_1})+\e\p\U_{\delta_1}\right]\right\},
\end{equation}
\beq\label{defelementi3}
\mathcal R_2:=\Pi^\bot\left\{-\p\U_{\delta_2}-i^*\left[f(V_\e)-f(\p\U_{\delta_1})-\e\p\U_{\delta_2}\right]\right\},
\eeq
\begin{equation}\label{defelementi1phi1}
\mathcal L_1(\phi_1):=\Pi^\bot_1\left\{\phi_1-i^*\left[f'(\p \U_1)\phi_1+\e \phi_1\right]\right\},
\end{equation}
\begin{equation}\label{defelementi1phi2}
\mathcal L_2(\phi_2):=\Pi^\bot\left\{\phi_2-i^*\left[f'(V_\e)\phi_2+\e \phi_2\right]\right\},
\end{equation}
\beq\label{defelementi4}
\mathcal{N}_1(\phi_1):=\Pi^{\perp}_1\{  - i^*[  f(\p\U_{\delta_1}+\phi_1) - f(\p\U_{\delta_1}) - f^\prime(\p\U_{\delta_1}) \phi_1]\},
\eeq and
\beq\label{defelementi5}
\mathcal{N}_2(\phi_1,\phi_2):= \displaystyle   \Pi^{\perp}\{   -i^*[f(V_\e+\phi_1+\phi_2)- f(V_\e)- f^\prime(V_\e) \phi_2 -f(\p\U_{\delta_1}+\phi_1)+f(\p\U_{\delta_1}) ]\}.
\eeq

We remark that it is not restrictive to consider $\mathcal R_1, \mathcal L(\phi_1), \mathcal N_1(\phi_1)\in \mathcal K_1^\bot$ since only $\delta_1$ appears and it is clear that a solution of \eqref{sistemaaux} gives a solution of \eqref{aux}.\\


%

Therefore we solve the first equation in \eqref{sistemaaux} finding a solution $\bar{\phi}_1=\bar{\phi}_1(\epsilon,d_1)$ and after that we solve the second equation in \eqref{sistemaaux} (with $\phi_1=\bar\phi_1$) finding also $\bar{\phi}_2=\bar{\phi}_2(\e,d_1,d_2)$. \\\\
Finally let us recall some useful inequality that we will use in the sequel.
Since these are known results, we omit the proof.

\begin{lemma}\label{lem00n2}
Let $\alpha$ be a positive real number.
If $\alpha \leq 1$ there holds
$$(x+y)^\alpha \leq x^\alpha + y^\alpha,$$
 for all $x, y >0$. If $\alpha \geq 1$ we have
$$(x+y)^\alpha \leq 2^{\alpha-1} (x^\alpha + y^\alpha),$$
 for all $x, y >0$.
 \end{lemma}

\begin{lemma}\label{lem0n2}
Let $q$ be a positive real number. There exists a positive constant $c$, depending only on $q$, such that for any $a, b\in \mathbb R$
\begin{equation}\label{lem1n2el}
||a+b|^q-|a|^q| \leq
\begin{cases}
c(q) \min\{|b|^q, |a|^{q-1}|b|\} &\ \hbox{if}\ 0<q<1,\\
c(q) (|a|^{q-1}|b|+|b|^q) & \ \hbox{if}\ q\geq1.
\end{cases}
\end{equation}
Moreover if $q>2$ then
\begin{equation}\label{lem1n2e2}
\left||a+b|^q-|a|^q-q |a|^{q-2}ab\right|\leq C\left(|a|^{q-2}|b|^2+|b|^q\right).
\end{equation}
\end{lemma}

\begin{lemma}\label{lem1n2}
Let $N\geq 7$. There exists a positive constant $c$, depending only on $p$, such that for any $a,b \in \mathbb R$
\begin{equation}\label{teslem1n2}
|f(a+b)-f(a)-f^\prime(a)b| \leq c |b|^{p}.
\end{equation}
\end{lemma}
\begin{lemma}\label{lem2n2}
There exists a positive constant $c$, depending only on $p$, such that for any $a,b \in \mathbb R$
\begin{equation}\label{eq1lem2n2}
|f(a-b)-f(a)+f(b)| \leq   c(p) (|a|^{p-1}|b|+|b|^p),
\end{equation}
or
\begin{equation}\label{eq2lem2n2}
|f(a-b)-f(a)+f(b)| \leq   c(p) (|b|^{p-1}|a|+|a|^p).
\end{equation}

\end{lemma}

\begin{lemma}\label{lem0n3}
Let $N\geq 7$. There exists a positive constant $c$ depending only on $p$ such that for any $a, b_1, b_2\in \mathbb R$
we get
\begin{equation}
\left|f(a+b_1)-f(a+b_2)-f'(a)(b_1-b_2)\right|\leq C\left(|b_1|^{p-1}+|b_2|^{p-1}\right)|b_1-b_2|.
\end{equation}
\end{lemma}

\section{The auxiliary equation: solution of the system \eqref{sistemaaux}}\label{ausiliaria}
We first define
\begin{equation}\label{thetaj}
\theta_1:=\frac{N-2}{N-4};\qquad\qquad \theta_2:=\frac{(N-2)^2}{(N-4)(N-6)}.
\end{equation}
We observe that $\theta_2$ is well defined since $N\geq 7$.
We also remark that having defined $\delta_j$ as in \eqref{deltaj}, $j=1,2$, the functions $\U_{\delta_j}$ depend on the parameters $d_j$, $j=1,2$.\\

In this section we solve system \eqref{sistemaaux}. More precisely, the aim is to prove the following result.

\begin{proposition}\label{auxsolving}
Let $N\geq 7$. 
For any $\eta>0$, there exist $\epsilon_0>0$ and $c>0$ such
that for all $\epsilon \in (0,\epsilon_0)$, for all $(d_1,d_2) \in
\mathbb R_+^2$ satisfying \eqref{limdj}, there exists a unique
$\bar\phi_1=\bar\phi_1(\epsilon,d_1)\in\mathcal K_1^\bot$ solution of the first equation of
\eqref{sistemaaux} such that $$\|\bar\phi_1\|\leq c
\e^{\frac{\theta_1}{2}+\sigma}$$ and there exists a unique solution
$ \bar\phi_2=\bar\phi_2(\e,d_1,d_2) \in {\mathcal K}^\bot$ of the second equation of
\eqref{sistemaaux} (with $\phi_1=\bar\phi_1$ such that $$ \|\bar\phi_2\| \leq c \
\epsilon^{\frac{\theta_2}{2}+\sigma},$$
 for some positive real number  $\sigma$ whose choice depends only on $N$.
Furthermore, $\bar\phi_1$ does not depend on $d_2$ and it is continuously differentiable with respect to $d_1$, $\bar\phi_2$ is continuously differentiable with respect to $(d_1, d_2)$.
\end{proposition}

In order to prove Proposition \ref{auxsolving} let us first consider the linear operator
$$\mathcal{L}_1: \mathcal{K}_1^\bot\rightarrow \mathcal{K}_1^{\bot}$$ defined as in
\eqref{defelementi1phi1}.\\ 
The next result provides an a-priori estimate for solutions $\phi\in\mathcal K_1^\bot$ of $\mathcal L_1(\phi)= h$, for some right-hand side $h$ with bounded $\|\cdot\|-$ norm.
\begin{lemma}\label{lem1invl1}
Let $N\geq 7$. For any $\eta>0$, there exists $\epsilon_0>0$ and $c>0$ such that for all $d_1 \in \mathbb R_+$  satisfying \eqref{limdj} for $j=1$, for all $\phi\in \mathcal{K}_1^\bot$ and for all $\epsilon \in (0,\epsilon_0)$ it holds $$\|\mathcal{L}_1(\phi)\| \geq c\|\phi\|.$$
\end{lemma}
\begin{proof}
For the proof it suffices to repeat with small changes the proof of Lemma 3.1 of \cite{Musso1}.
\end{proof}
Next result states the invertibility of the operator $\mathcal{L}_1$ and provides a uniform estimate on the inverse of the operator $\mathcal L_1$.

\begin{proposition}\label{linearprop1}
Let $N\geq 7$. For any $\eta>0$, there exists $\epsilon_0>0$ and $c>0$ such that the linear operator $\mathcal{L}_1$ is invertible and $\| \mathcal{L}_1^{-1}\|\leq c$ for all $\epsilon \in (0,\epsilon_0)$, for all $d_1 \in \mathbb R_+$ satisfying \eqref{limdj} for $j=1$.
\end{proposition}
\begin{proof}
For the proof it suffices to repeat with small changes the proof of Proposition 3.2 of \cite{Musso1}.
\end{proof}

For the linear operator $\mathcal{L}_2$ we state analogous results.

\begin{lemma}\label{lem1invl2}
Let $N\geq 7$. For any $\eta>0$, there exists $\epsilon_0>0$ and $c>0$ such that for all $(d_1,d_2) \in \mathbb R_+^2$  satisfying \eqref{limdj}, for all $\phi\in \mathcal{K}^\bot$ and for all $\epsilon \in (0,\epsilon_0)$ it holds $$\|\mathcal{L}_2(\phi)\| \geq c\|\phi\|.$$
\end{lemma}
\begin{proof}
For the proof see Lemma 3.1 of \cite{Musso1}.
\end{proof}

\begin{proposition}\label{linearprop2}
Let $N\geq 7$. For any $\eta>0$, there exists $\epsilon_0>0$ and $c>0$ such that the linear operator $\mathcal{L}_2$ is invertible and $\| \mathcal{L}_2^{-1}\|\leq c$ for all $\epsilon \in (0,\epsilon_0)$, for all $(d_1,d_2) \in \mathbb R_+^2$ satisfying \eqref{limdj}.
\end{proposition}
\begin{proof}
For the proof see Proposition 3.2 of \cite{Musso1}.
\end{proof}

The strategy is to solve the first equation of \eqref{sistemaaux} by
a fixed point argument, finding a unique $\bar\phi_1$ and then,
substituting  $\bar\phi_1$ in the second equation of
\eqref{sistemaaux}, we obtain an equation depending only on the
variable $\phi_2$. Hence, using again a fixed point argument, we
solve the second equation of \eqref{sistemaaux} uniquely.

\subsection{The solution of the first equation of \eqref{sistemaaux}}
The aim is to prove the following proposition.
\begin{proposition}\label{aux1}
Let $N\geq 7$. For any $\eta>0$, there exists $\epsilon_0>0$ and $c>0$ such that for all $\epsilon \in (0,\epsilon_0)$, for all $d_1 \in \mathbb R_+$ satisfying condition \eqref{limdj} for $j=1$, there exists a unique solution $\bar\phi_1=\bar\phi_1(\epsilon,d_1)$, $\bar\phi_1\in\mathcal{K}_1^\bot$ of the first equation in \eqref{sistemaaux} which is continuously differentiable with respect to $d_1$ and such that
\begin{equation}\label{stimaphi1}
\|\bar\phi_1\|\leq c \e^{\frac{\theta_1}{2}+\sigma},
\end{equation}
where $\theta_1$ is defined in \eqref{thetaj} and $\sigma$ is some positive real number whose choice depends only on $N$.
\end{proposition}
In order to prove Proposition \ref{aux1} we have to estimate the error term $\mathcal{R}_1$ defined in \eqref{defelementi2}. It holds the following result.

\begin{proposition}\label{errorprop1}
Let $N\geq 7$. For any $\eta>0$, there exists $\epsilon_0>0$ and $c>0$ such that  for all $\epsilon \in (0,\epsilon_0)$, for all $d_1 \in \mathbb R_+$ satisfying condition \eqref{limdj} for $j=1$, we have $$\| \mathcal{R}_1\|\leq c \ \epsilon^{\frac{\theta_1}{2}+\sigma},$$
 for some positive real number  $\sigma$ whose choice depends only on $N$.
\end{proposition}
\begin{proof}
By continuity of $\Pi_1^\bot$, by using \eqref{stimaistar} and since $\p\U_{\delta_1}$ weakly solves $-\Delta\p\U_{\delta_1}=\U_{\delta_1}^p$ in $\O$, it follows that
\begin{eqnarray*}
\|\mathcal{R}_1\|&=&\|\Pi_1^\bot\left\{\p\U_{\delta_1}-i^*\left[f(\p\U_{\delta_1})+\e\p\U_{\delta_1}\right]\right\}\|\leq C_1\|\p\U_{\delta_1}-i^*\left[f(\p\U_{\delta_1})+\e\p\U_{\delta_1}\right]\|\\
&\leq & C_2\left|f(\U_{\delta_1})-f(\p\U_{\delta_1})-\e\p\U_{\delta_1}\right|_{\frac{2N}{N+2}}\leq  C\underbrace{\left|f(\U_{\delta_1})-f(\p\U_{\delta_1})\right|_{\frac{2N}{N+2}}}_{(I)}+\underbrace{\e\left|\p\U_{\delta_1}\right|_{\frac{2N}{N+2}}}_{(II)}.
\end{eqnarray*}
Let us fix $\eta>0$. We estimate the terms $(I), (II)$.
\\\\ {\it Claim 1:}
\begin{equation}\label{eq1err1}
(I)= O(\epsilon^{\frac{N+2}{2(N-4)}}) .
\end{equation}

By using \eqref{stimaproiezione}, \eqref{stimavarphi} and by elementary inequalities we get
\begin{equation*}
\begin{array}{lll}
\displaystyle \int_{\Omega}|(\p\U_{\delta_1})^p -\U_{\delta_1}^p|^{\frac{2N}{N+2}}\ dx &\leq & \displaystyle c_1 \int_{\Omega}|\U_{\delta_1}^{p-1}\varphi_{\delta_1}|^{\frac{2N}{N+2}}\  dx + c_2 \int_{\Omega}|\varphi_{\delta_1}|^{p+1}\  dx \\[10pt]
&\leq & \displaystyle c_3 \delta_1^{\frac{N-2}{2}\frac{2N}{N+2}} \int_{\Omega}\left(\frac{\delta_1^2}{(\delta_1^2+|x|^2)^2}\right)^{\frac{2N}{N+2}}\  dx + c_2|\varphi_{\delta_1}|_{p+1, \O}^{p+1}\\[10pt]
&= & \displaystyle c_3 \delta_1^{\frac{N(N-2)}{N+2}} \int_{\Omega}\left(\frac{\delta_1^2}{(\delta_1^2+|x|^2)^2}\right)^{\frac{2N}{N+2}}\  dx + c_4 \delta_1^{N}.
\end{array}
\end{equation*}
Now for $N\geq 7$ we have $$\int_{\Omega}\left(\frac{\delta_1^2}{(\delta_1^2+|x|^2)^2}\right)^{\frac{2N}{N+2}}\  dx = O\left(\delta_1^{\frac{4N}{N+2}}\right).$$ Indeed:
\begin{equation*}
\displaystyle \int_{\Omega}\left(\frac{\delta_1^2}{(\delta_1^2+|x|^2)^2}\right)^{\frac{2N}{N+2}}\  dx = \displaystyle \delta_1^{\frac{4N}{N+2}} \int_{\Omega} \frac{1}{(\delta_1^2+|x|^2)^{\frac{4N}{N+2}}} \ dx\leq \displaystyle \delta_1^{\frac{4N}{N+2}} \int_{\Omega} \frac{1}{|x|^{\frac{8N}{N+2}}} \ dx,
\end{equation*}
and the last integral is finite since $N>6$, which implies $\frac{8N}{N+2}<N$.
Finally, since $\frac{4N^2(N-2)}{(N+2)^2} > N$, for any $N>4$, we deduce that
$$\displaystyle \int_{\Omega}|(\p\U_{\delta_1})^p -\U_{\delta_1}^p|^{\frac{2N}{N+2}}\ dx =  O\left(\delta_1^{N}\right),$$
and hence
\begin{equation}\label{numero}
|(\p\U_{\delta_1})^p -\U_{\delta_1}^p|_{\frac{2N}{N+2}}=O\left(\delta_1^{\frac{N+2}{2}}\right).
\end{equation}
Since $\delta_1=d_1 \epsilon^{\frac{1}{N-4}}$ and $d_1$ satisfies \eqref{limdj}, we get that $|(\p\U_{\delta_1})^p -\U_{\delta_1}^p|_{\frac{2N}{N+2}}=O\left(\epsilon^{\frac{N+2}{2(N-4)}}\right)$ and Claim 1 is proved.
\\\\ {\it Claim 2:}
\begin{equation}\label{eq2err1}
(II)= O(\epsilon^{\frac{N-2}{N-4}}).
\end{equation}

\begin{eqnarray}\label{utilerr1}
\nonumber
\displaystyle \int_{\Omega} \p\U_{\delta_{1}}^{\frac{2N}{N+2}} \ dx &\leq& \displaystyle \int_{\Omega} \U_{\delta_{1}}^{\frac{2N}{N+2}} \ dx=\displaystyle  \alpha_N^{\frac{2N}{N+2}}  \int_{\Omega}\frac{\delta_1^{-\frac{N(N-2)}{N+2}}}{(1+|\frac{x}{\delta_1}|^2)^{\frac{N(N-2)}{N+2}}}\  dx \\
&&=\displaystyle  \alpha_N^{\frac{2N}{N+2}} \delta_1^{\frac{4N}{N+2}} \int_{\R}\frac{1}{(1+|y|^2)^{\frac{N(N-2)}{N+2}}}\  dy + o(\delta_1^{\frac{4N}{N+2}}).
\end{eqnarray}
Thus, since $\delta_1=d_1 \epsilon^{\frac{1}{N-4}}$ and $d_1$ satisfies \eqref{limdj}, we get that
$$ \int_{\Omega} \p\U_{\delta_{1}}^{\frac{2N}{N+2}} \ dx = O\left(\epsilon^{\frac{4N}{(N+2)(N-4)}}\right),$$
and hence
$$\epsilon \left(\int_{\Omega} P\U_{\delta_{1}}^{\frac{2N}{N+2}} \ dx\right)^{\frac{N+2}{2N}} = \epsilon O\left(\epsilon^{\frac{2}{N-4}}\right)= O\left(\epsilon^{\frac{N-2}{N-4}}\right).$$ 

The proof of Claim 2 is complete.\\ \\ Hence, by \eqref{eq1err1} and \eqref{eq2err1}, we deduce that there exist a constant $c=c(\eta)>0$  and $\epsilon_0=\epsilon_0(\eta)>0$ sufficiently small such that, for all $\epsilon \in (0,\epsilon_0)$ and $d_1 \in \mathbb R_+$ satisfying \eqref{limdj} (with $j=1$)
$$\| \mathcal{R}_1\|\leq {c} \ (\epsilon^{\frac{N+2}{2(N-4)}} + \epsilon^{\frac{N-2}{N-4}})\leq c \e^{\frac{\theta_1}{2}+\sigma},$$  with $\sigma$ such that $0<\sigma< \frac{2}{N-4}$.
\end{proof}
We are ready to prove Proposition \ref{aux1}.\\
\begin{proof}[Proof of Proposition \ref{aux1}.]
Let us fix $\eta>0$ and define $\mathcal{T}_1 :{\mathcal K}_1^\bot \rightarrow {\mathcal K}_1^\bot $ as $$\mathcal{T}_1(\phi_1):=- \mathcal{L}_1^{-1}[\mathcal{N}_1(\phi_1)+\mathcal{R}_1].$$ Clearly solving the first equation of \eqref{sistemaaux} is equivalent to solving the fixed point equation $\mathcal{T}_1(\phi_1)=\phi_1$.\\ Let us define the ball
$$\displaystyle B_{1,\epsilon}:=\{\phi_1 \in \mathcal{K}^{\perp}; \ \|\phi_1\| \leq r \ \epsilon^{\frac{\theta_1}{2}+\sigma} \}\subset \mathcal{K}^\bot$$ with $r>0$ sufficiently large and $\sigma>0$.\\
We want to prove that, for $\e$ small, $\mathcal{T}_1$ is a contraction in the proper ball $\displaystyle B_{1,\epsilon}$, namely we want to prove that, for $\e$ sufficiently small
\begin{enumerate}
\item $\mathcal{T}_1(B_{1,\epsilon}) \subset B_{1,\epsilon}$;
\item $\|\mathcal{T}_1\|<1$.
\end{enumerate}

By Lemma \ref{lem1invl1} we get:
\begin{equation}\label{eq1n1}
\|\mathcal{T}_1(\phi_1)\| \leq c (\|\mathcal{N}_1(\phi_1)\| + \|\mathcal{R}_1\|)
\end{equation}
and
\begin{equation}\label{eq2n1}
\|\mathcal{T}_1(\phi_1)-\mathcal{T}_1(\psi_1)\| \leq c (\|\mathcal{N}_1(\phi_1)-\mathcal{N}_1(\psi_1)\|),
\end{equation}
for all $\phi_1, \psi_1 \in {\mathcal K}_1^\perp$. Thanks to \eqref{stimaistar} and the definition of $\mathcal{N}_1$ we deduce that
\begin{equation}\label{eq3n1}
\|\mathcal{N}_1(\phi_1)\| \leq c | f(\p\U_{\delta_1}+\phi_1) - f(\p\U_{\delta_1}) - f^\prime(\p\U_{\delta_1}) \phi_1|_{\frac{2N}{N+2}},
\end{equation}
and
\begin{equation}\label{eq4n1}
\displaystyle \|\mathcal{N}_1(\phi_1)-\mathcal{N}_1(\psi_1)\| \leq \displaystyle c | f(\p\U_{\delta_1}+\phi_1)- f(\p\U_{\delta_1}+\psi_1) - f^\prime(\p\U_{\delta_1}) (\phi_1 -\psi_1)|_{\frac{2N}{N+2}} .
\end{equation}
Now we estimate the  right-hand term in (\ref{eq1n1}).
Thanks to Lemma \ref{lem1n2} we have the following inequality:
\begin{equation}\label{eq5n1}
 | f(\p\U_{\delta_1}+\phi_1) - f(\p\U_{\delta_1}) - f^\prime(\p\U_{\delta_1}) \phi_1|\leq c |\phi_1|^p.
\end{equation}
Since $p \frac{2N}{N+2} = \frac{2N}{N-2}$ and $|\phi_1^p|_{\frac{2N}{N+2}}=|\phi_1|_{\frac{2N}{N-2}}^{p}$, from (\ref{eq5n1}) and the Sobolev inequality we deduce the following:
\begin{equation}\label{eq6n1}
 | f(\p\U_{\delta_1}+\phi_1) - f(\p\U_{\delta_1}) - f^\prime(\p\U_{\delta_1}) \phi_1|_{\frac{2N}{N+2}} \leq c_1 |\phi_1|_{\frac{2N}{N-2}}^p \leq c_2 \|\phi_1\|^p.
\end{equation}
Thanks to (\ref{eq1n1}), Proposition \ref{errorprop1}, (\ref{eq3n1}), (\ref{eq6n1}) and since $p>1$, then, there exist $c=c(\eta)>0$ and $\epsilon_0=\epsilon_0(\eta)>0$ such that  $$\|\phi_1\| \leq c \epsilon^{\frac{\theta_1}{2}+ \sigma} \Rightarrow \|\mathcal{T}_1(\phi_1) \|\leq c \epsilon^{\frac{\theta_1}{2}+ \sigma},$$
for all $\epsilon \in (0,\epsilon_0)$, for all $d_1 \in \mathbb R_+$ satisfying \eqref{limdj} (with $j=1$), for some positive real number  $\sigma$, whose choice depends only on $N$.
In other words $\mathcal{T}_1$ maps the ball $B_{1,\epsilon}$ into itself and (1) is proved.

We want to show that $\mathcal{T}_1$ is a contraction. By using Lemma \ref{lem0n3} we get that
for any $\phi_1, \psi_1\in B_{1, \e}$
\begin{equation*}
\left|f(\p\U_{\delta_1}+\phi_1)-f(\p\U_{\delta_1}+\psi_1)-f'(\p\U_{\delta_1})(\phi_1-\psi_1)\right|\leq C\left(|\phi_1|^{p-1}+|\psi_1|^{p-1}\right)|\phi_1-\psi_1|.
\end{equation*}

By direct computation $(p-1) \frac{2N}{N+2}=\frac{8N}{(N-2)(N+2)}$, so, since $ |\phi_1|^{(p-1) \frac{2N}{N+2}}, |\psi_1|^{(p-1) \frac{2N}{N+2}} \in L^{\frac{N+2}{4}}$, $|\phi_1-\psi_1|^{\frac{2N}{N+2}} \in L^{p}$ and $1= \frac{4}{N+2}+\frac{N-2}{N+2}$ by H\"older inequality we get that
\begin{eqnarray}\label{eq8n1}
\nonumber
&&\hskip-1.5cm \displaystyle \left| \left( |\phi_1|^{p-1}+|\psi_1|^{p-1}\right) (\phi_1 - \psi_1) \right|_{\frac{2N}{N+2}} \leq  \displaystyle\left[ \left(| \phi_1|_{\frac{2N}{N-2}}^{\frac{4}{N-2}} + |\psi_1|_{\frac{2N}{N-2}}^{\frac{4}{N-2}}\right)^{\frac{2N}{N+2}} \ \left( | \phi_1 - \psi_1|_{\frac{2N}{N-2}}^{{\frac{2N}{N-2}}} \right)^{\frac{N-2}{N+2}} \right]^{\frac{N+2}{2N}}\\
&=& \displaystyle \left(| \phi_1|_{\frac{2N}{N-2}}^{\frac{4}{N-2}} + |\psi_1 |_{\frac{2N}{N-2}}^{\frac{4}{N-2}}\right) \ | \phi_1 - \psi_1|_{\frac{2N}{N-2}}.
\end{eqnarray}

Hence by (\ref{eq2n1}), (\ref{eq4n1}), (\ref{eq8n1}) and Sobolev inequality we get that there exists $L\in (0,1)$ such that
$$\|\phi_1\| \leq c \epsilon^{\frac{\theta_1}{2}+ \sigma}, \|\psi_1\| \leq c \epsilon^{\frac{\theta_1}{2}+ \sigma} \Rightarrow \|\mathcal{T}_1(\phi_1)- \mathcal{T}_1(\psi_1)\| \leq L \|\phi_1-\psi_1\|.$$

Hence by the Contraction Mapping Theorem we can
uniquely solve $\mathcal{T}_1(\phi_1)=\phi_1$ in $B_{1,\epsilon}$. We denote by $\bar\phi_1 \in B_{1,\epsilon}$ this solution.  A standard argument shows that $d_1 \rightarrow \bar\phi_1(d_1)$ is a $C^1$-map (see also \cite{Musso1}). The proof is then concluded.
\end{proof}

\subsection{The proof of Proposition \ref{auxsolving}}

Before proving Proposition \ref{auxsolving} we need some preliminary results, in particular we need to improve the estimate on the solution $\bar\phi_1$ of the first equation of (\ref{sistemaaux}) found in Proposition \ref{aux1}.\\

The first preliminary result is an estimate on the error term $\mathcal{R}_2$ defined in \eqref{defelementi3}. 

\begin{proposition}\label{errorprop2}
For any $\eta>0$, there exists $\epsilon_0>0$ and $c>0$ such that  for all $\epsilon \in (0,\epsilon_0)$, for all $(d_1,d_2) \in \mathbb R_+^2$ satisfying \eqref{limdj}, we have $$\| \mathcal{R}_2\|\leq c \ \epsilon^{\frac{\theta_2}{2}+\sigma},$$ for some positive real number  $\sigma$, whose choice depends only on $N$.
\end{proposition}
\begin{proof}
By continuity of $\Pi^\bot$ and by using \eqref{stimaistar} we deduce that
\begin{eqnarray}\label{eq1err2}
\nonumber
\displaystyle \|\mathcal{R}_2\| &\leq& \displaystyle c |f(\U_{\delta_2}) + f(\p\U_{\delta_1}- \p\U_{\delta_2}) - f(\p\U_{\delta_1}) - \epsilon \p\U_{\delta_2}|_{\frac{2N}{N+2}} \\[10pt]
\nonumber
&\leq& \underbrace{\displaystyle c |f(\p\U_{\delta_1}- \p\U_{\delta_2}) - f(\p\U_{\delta_1})+ f(\p\U_{\delta_2})  |_{\frac{2N}{N+2}}}_{(I)}+\underbrace{\displaystyle c |f(\p\U_{\delta_2})-f(\U_{\delta_2})|_{\frac{2N}{N+2}}}_{(II)}  \\[10pt]
&& + \underbrace{c\epsilon |\p\U_{\delta_2}|_{\frac{2N}{N+2}}}_{(III)}.
\end{eqnarray}
Let us fix $\eta>0$. We begin estimating $(I)$. Let $\rho>0$ so that $B(0,\rho) \subset \Omega$. We decompose the domain $\Omega$ as $\Omega=A_0 \sqcup A_1 \sqcup A_2$, where
$A_0:=\Omega \setminus B(0,\rho)$,
$A_1:=B(0,\rho) \setminus B(0,\sqrt{\delta_1 \delta_2})$ and $A_2:= B(0,\sqrt{\delta_1 \delta_2})$. We evaluate $(I)$ in every set of this decomposition.

Thanks to  Lemma \ref{lem2n2}  there exists a positive constant $c$ (depending only on $p$) such that
\begin{equation}\label{eq6err2}
|f(\p\U_{\delta_1}- \p\U_{\delta_2}) - f(\p\U_{\delta_1})+ f(\p\U_{\delta_2}) |\leq c (\p\U_{\delta_1}^{p-1}\p\U_{\delta_2} +\p\U_{\delta_2}^p ).
\end{equation}
Integrating on $A_0$ and using the usual elementary inequalities (see Lemma \ref{lem00n2}) we get that
\begin{equation}\label{elemcompututil}
\begin{array}{lll}
&& \hskip-2.0cm\displaystyle \int_{A_0}  |f(\p\U_{\delta_1}- \p\U_{\delta_2}) - f(\p\U_{\delta_1})+ f(\p\U_{\delta_2})  |^{\frac{2N}{N+2}} \ dx\\
& \leq& \displaystyle C_1 \int_{A_0}(\p\U_{\delta_1}^{(p-1)(\frac{2N}{N+2})} \p\U_{\delta_2}^{\frac{2N}{N+2}} + \p\U_{\delta_2}^{p+1})  \ dx\\
&\leq &\displaystyle C_2\int_{A_0}\frac{\delta_1^{\frac{4N}{N+2}}}{(\delta_1^2+|x|^2)^{\frac{4N}{N+2}}}\frac{\delta_2^{\frac{N(N-2)}{N+2}}}{(\delta_2^2+|x|^2)^{\frac{N(N-2)}{N+2}}}\, dx+C_3\int_{A_0} \frac{\delta_2^N}{(\delta_2^2+|x|^2)^N} \, dx\\
&\leq & \displaystyle C_4 \frac{\delta_1^{\frac{4N}{N+2}}}{\rho^{\frac{8N}{N+2}}}\frac{\delta_2^{\frac{N(N-2)}{N+2}}}{\rho^{2\frac{N(N-2)}{N+2}}}+
C_5 \frac{\delta_2^N}{\rho^{2N}}
\end{array}
\end{equation}
and hence we deduce that (recall the choice of $\delta_1, \delta_2$ (see \eqref{deltaj}))

 \begin{equation}\label{eq5err2}
  |f(\p\U_{\delta_1}- \p\U_{\delta_2}) - f(\p\U_{\delta_1})+ f(\p\U_{\delta_2})  |_{\frac{2N}{N+2},A_0} \leq c  \epsilon^{\frac{3N^2-12N-4}{2(N-4)(N-6)}}\leq c \e^{\frac{\theta_2}{2}+\sigma}
\end{equation}
where $c$ depends on $\eta$ (and also on $\Omega$, $\rho$, $N$), $\sigma$ is some positive real number (to be precise we can choose $0<\sigma\leq \frac{N^2-4N-4}{(N-4)(N-6)}$).

We evaluate now $(I)$ in $A_1$. By (\ref{eq6err2}) and the usual elementary inequalities we deduce the following:
\begin{equation}\label{eq1I1A1}
\int_{A_1}  |f(P\U_{\delta_1}- \p\U_{\delta_2}) - f(\p\U_{\delta_1})+ f(\p\U_{\delta_2})  |^{\frac{2N}{N+2}} \ dx \leq c \int_{A_1}(\p\U_{\delta_1}^{(p-1)(\frac{2N}{N+2})} \p\U_{\delta_2}^{\frac{2N}{N+2}} + \p\U_{\delta_2}^{p+1} ) \ dx.
\end{equation}
Let us estimate every term:
 \begin{eqnarray}\label{eq2I1A1}
\begin{array}{lll}
 &&\hskip-1.5cm \displaystyle\int_{A_1} \p\U_{\delta_1}^{(p-1)(\frac{2N}{N+2})}\p\U_{\delta_2}^{\frac{2N}{N+2}}  dx \leq \displaystyle  \int_{A_1} \U_{\delta_1}^{(p-1)(\frac{2N}{N+2})}\U_{\delta_2}^{\frac{2N}{N+2}} dx =\alpha_N^{p+1} \displaystyle \int_{A_1}  \frac{\delta_1^{\frac{4N}{N+2}}}{(\delta_1^2+|x|^2)^{\frac{4N}{N+2}}}  \frac{\delta_2^{\frac{N(N-2)}{N+2}}}{(\delta_2^2+|x|^2)^{\frac{N(N-2)}{N+2}}} \ dx \\[12pt]
  &=& \displaystyle c_1  \int_{\sqrt{\delta_1 \delta_2}}^{\rho}  \frac{\delta_1^{\frac{4N}{N+2}}}{(\delta_1^2+r^2)^{\frac{4N}{N+2}}}  \frac{\delta_2^{\frac{N(N-2)}{N+2}}}{(\delta_2^2+r^2)^{\frac{N(N-2)}{N+2}}} r^{N-1} dr = \displaystyle c_1  \int_{\sqrt{\frac{\delta_1} {\delta_2}}}^{\frac{\rho}{\delta_2}}  \frac{\delta_1^{\frac{4N}{N+2}}}{(\delta_1^2+\delta_2^2s^2)^{\frac{4N}{N+2}}}  \frac{\delta_2^{-\frac{N(N-2)}{N+2}}}{(1+s^2)^{\frac{N(N-2)}{N+2}}} \delta_2^N s^{N-1} ds \\[12pt]
  &=& \displaystyle c_1  \int_{\sqrt{\frac{\delta_1} {\delta_2}}}^{\frac{\rho}{\delta_2}}  \frac{\delta_1^{-\frac{4N}{N+2}}}{\left[1+\left(\frac{\delta_2}{\delta_1}\right)^2s^2\right]^{\frac{4N}{N+2}}}  \frac{\delta_2^{\frac{4N}{N+2}}}{(1+s^2)^{\frac{N(N-2)}{N+2}}}  s^{N-1} ds \leq \displaystyle c_1  \left(\frac{\delta_2}{\delta_1}\right)^{\frac{4N}{N+2}}  \int_{\sqrt{\frac{\delta_1} {\delta_2}}}^{\frac{\rho}{\delta_2}}   \frac{1}{(1+s^2)^{\frac{N(N-2)}{N+2}}}  s^{N-1} ds  \\[12pt]
        &\leq& \displaystyle c_1 \left(\frac{\delta_2}{\delta_1}\right)^{\frac{4N}{N+2}}  \int_{\sqrt{\frac{\delta_1} {\delta_2}}}^{\frac{\rho}{\delta_2}}   \frac{1}{s^{\frac{N^2-5N+2}{N+2}}}   ds= \displaystyle c_2  \left(\frac{\delta_2}{\delta_1}\right)^{\frac{4N}{N+2}} \left[\left(\frac{\delta_2}{\delta_1}\right)^{\frac{N^2-6N}{2(N+2)}} -\left(\frac{\delta_2}{\rho}\right)^{\frac{N^2-6N}{(N+2)}}\right]
\\[12pt]
&\leq&\displaystyle c_3 \left(\frac{\delta_2}{\delta_1}\right)^{\frac{N}{2}}.
\end{array}
\end{eqnarray}
Moreover
\begin{equation}\label{eq3I1A1}
\int_{A_1}\p\U_{\delta_2}^{p+1}\, dx \leq \int_{A_1}\U_{\delta_2}^{p+1}\, dx \leq C_1 \int_{\sqrt{\frac{\delta_1}{\delta_2}}}^{\frac{\rho}{\delta_2}}\frac{r^{N-1}}{\left(1+r^2\right)^N}\, dr\leq C_2\left(\frac{\delta_2}{\delta_1}\right)^{\frac{N}{2}}.
\end{equation}
Thanks to the choice of $\delta_1$, $\delta_2$ we have
\begin{equation}\label{orddelta}
\left(\frac{\delta_2}{\delta_1}\right)^{\frac{N}{2}}=O(\epsilon^{\frac{N(N-2)}{(N-4)(N-6)}}).
\end{equation}
Hence, from (\ref{eq1I1A1}), (\ref{eq2I1A1}), (\ref{eq3I1A1}) and (\ref{orddelta}) we deduce that
 \begin{equation}\label{eq4I1A1}
  |f(\p\U_{\delta_1}- \p\U_{\delta_2}) - f(\p\U_{\delta_1})+ f(\p\U_{\delta_2})  |_{\frac{2N}{N+2},A_1} \leq c  \epsilon^{\frac{(N-2)(N+2)}{2(N-4)(N-6)}}\leq c \e^{\frac{\theta_2}{2}+\sigma},
\end{equation}
where $c$ depends on $\eta$, $\sigma$ is some positive real number (to be precise we can choose $0<\sigma\leq \frac{2(N-2)}{(N-4)(N-6)}$).

Now we evaluate $(I)$ in $A_2$. To do this we apply (\ref{eq2lem2n2}) of Lemma \ref{lem2n2}, so  there exists a constant $c>0$ such that
\begin{equation}\label{eq8err2}
|f(\p\U_{\delta_1}- \p\U_{\delta_2}) - f(\p\U_{\delta_1})+ f(\p\U_{\delta_2}) |\leq c (\p\U_{\delta_2}^{p-1}\p\U_{\delta_1} +\p\U_{\delta_1}^p ).
\end{equation}
Thanks to (\ref{eq8err2}) and the usual elementary inequalities we deduce the following:
\begin{equation}\label{eq9err2}
\int_{A_2}  |f(\p\U_{\delta_1}- \p\U_{\delta_2}) - f(\p\U_{\delta_1})+ f(\p\U_{\delta_2})  |^{\frac{2N}{N+2}} \ dx \leq c \int_{A_2}(\p\U_{\delta_2}^{(p-1)(\frac{2N}{N+2})} \p\U_{\delta_1}^{\frac{2N}{N+2}} + \p\U_{\delta_1}^{p+1} ) \ dx.
\end{equation}
We estimate the first term
 \begin{equation}\label{eq11err2}
\begin{array}{lll}
 &&\hskip-1.5cm\displaystyle\int_{A_2} \p\U_{\delta_2}^{(p-1)(\frac{2N}{N+2})}\p\U_{\delta_1}^{\frac{2N}{N+2}}  dx \leq \displaystyle  \int_{A_2} \U_{\delta_2}^{(p-1)(\frac{2N}{N+2})}\U_{\delta_1}^{\frac{2N}{N+2}} dx =\alpha_N^{p+1} \displaystyle \int_{A_2}  \frac{\delta_2^{\frac{4N}{N+2}}}{(\delta_2^2+|x|^2)^{\frac{4N}{N+2}}}  \frac{\delta_1^{\frac{N(N-2)}{N+2}}}{(\delta_1^2+|x|^2)^{\frac{N(N-2)}{N+2}}} \ dx \\[12pt]
  &=& \displaystyle c_1  \int_0^{\sqrt{\frac{\delta_2} {\delta_1}}} \frac{\delta_2^{\frac{4N}{N+2}}}{(\delta_2^2+\delta_1^2s^2)^{\frac{4N}{N+2}}}  \frac{\delta_1^{-\frac{N(N-2)}{N+2}}}{(1+s^2)^{\frac{N(N-2)}{N+2}}} \delta_1^N s^{N-1} ds = \displaystyle c_1 \int_0^{\sqrt{\frac{\delta_2} {\delta_1}}} \frac{\delta_2^{\frac{4N}{N+2}}}{\left[\left(\frac{\delta_2}{\delta_1}\right)^2+s^2\right]^{\frac{4N}{N+2}}}  \frac{\delta_1^{-\frac{4N}{N+2}}}{(1+s^2)^{\frac{N(N-2)}{N+2}}}  s^{N-1} ds \\[30pt]
    &\leq & \displaystyle c_1  \left(\frac{\delta_2}{\delta_1}\right)^{\frac{4N}{N+2}}  \int_0^{\sqrt{\frac{\delta_2} {\delta_1}}}    \frac{1}{s^{\frac{8N}{N+2}}(1+s^2)^{\frac{N(N-2)}{N+2}}}  s^{N-1} ds \leq  \displaystyle c_1  \left(\frac{\delta_2}{\delta_1}\right)^{\frac{4N}{N+2}}  \int_0^{\sqrt{\frac{\delta_2} {\delta_1}}}    \frac{s^{\frac{N^2-7N-2}{N+2}}}{(1+s^2)^{\frac{N(N-2)}{N+2}}}   ds \\[12pt]
        &\leq& \displaystyle c_1  \left(\frac{\delta_2}{\delta_1}\right)^{\frac{4N}{N+2}}   \int_0^{\sqrt{\frac{\delta_2} {\delta_1}}}    {s^{\frac{N^2-7N-2}{N+2}}} ds = \displaystyle c_2  \left(\frac{\delta_2}{\delta_1}\right)^{\frac{4N}{N+2}} \left(\frac{\delta_2}{\delta_1}\right)^{\frac{N^2-6N}{2(N+2)}}\\[8pt]
        &=&\displaystyle c_2 \left(\frac{\delta_2}{\delta_1}\right)^{\frac{N}{2}}.
\end{array}
\end{equation}
By making similar computations as before we get that
\begin{equation}\label{eqerrore2post}
\int_{A_2}\p\U_{\delta_1}^{p+1}\, dx \leq c_3 \left(\frac{\delta_2}{\delta_1}\right)^{\frac{N}{2}}.
\end{equation}
So from (\ref{eq9err2}) and (\ref{eq11err2}) we deduce that
\begin{equation}\label{eq12err2}
  |f(\p\U_{\delta_1}- \p\U_{\delta_2}) - f(\p\U_{\delta_1})+ f(\p\U_{\delta_2})  |_{\frac{2N}{N+2},A_2}  \leq c \epsilon^{\frac{(N+2)(N-2)}{2(N-4)(N-6)}}\leq c \e^{\frac{\theta_2}{2}+\sigma},
\end{equation}
where $c$ depends on $\eta$, $\sigma$ is some positive real number (to be precise we can choose $0<\sigma\leq \frac{2(N-2)}{(N-4)(N-6)}$).
Hence from (\ref{eq5err2}), (\ref{eq4I1A1}) and (\ref{eq12err2}) we deduce that
\begin{equation}\label{eqI1F}
(I) \leq c \epsilon^{\frac{\theta_2}{2} + \sigma},
\end{equation}
for some positive constant $c$, for some positive real number $\sigma$ depending only on $N$.\\
Now by making similar computations as for (I) of Proposition \ref{errorprop1} (see \eqref{numero}) we get that
$$(II)=O\left(\delta_2^{\frac{N+2}{2}}\right),$$
and hence we deduce that
$$(II) \leq c\e^{\frac{(3N-10)(N+2)}{2(N-4)(N-6)}} \leq c \e^{\frac{\theta_2}{2}+\sigma},$$
where $c$, $0<\sigma\leq \frac{N^2-12}{(N-4)(N-6)}$.

It remains to estimate $(III)$.\\

From (\ref{utilerr1}), exchanging $\delta_1$ with $\delta_2$ we get:
\begin{equation*}
\begin{array}{lll}
\displaystyle \int_{\Omega} \p\U_{\delta_{2}}^{\frac{2N}{N+2}} \ dx &\leq& \displaystyle  \alpha_N^{\frac{2N}{N+2}} \delta_2^{\frac{4N}{N+2}} \int_{\R}\frac{1}{(1+|y|^2)^{\frac{N(N-2)}{N+2}}}\  dy.
\end{array}
\end{equation*}
Hence we deduce that $(III) \leq c \ \epsilon \delta_2^2$, and thanks to the choice $\delta_2$, by an elementary computation, we get that:
$$(III)\leq c \ \epsilon^{\frac{(N-2)^2}{(N-4)(N-6)}}\leq c\e^{\frac{\theta_2}{2}+\sigma},$$ where $c$, $0<\sigma\leq \frac{(N-2)^2}{2(N-4)(N-6)}$.
Finally, putting together all these estimates we deduce that there exist a positive constant $c=c(\eta)>0$ and $\epsilon_0=\epsilon_0(\eta)>0$ such that for all $\epsilon \in (0,\epsilon_0)$, for all $(d_1,d_2) \in \mathbb R_+^2$ satisfying \eqref{limdj}
$$ \|\mathcal{R}_2\| \leq c \epsilon^{\frac{\theta_2}{2} + \sigma},$$
for some positive real number $\sigma$ (whose choice depends only on $N$).
The proof is complete.
\end{proof}

Now we prove a technical result on the behavior of the $L^\infty$-norm of $\bar\phi_1$, which will be useful in the sequel.

\begin{lemma}\label{techlembehvphi1}
Let $\eta$ be a small positive real number and let $\bar\phi_1\in\mathcal K_1^\bot$ be the solution of the first
equation in \eqref{sistemaaux}, found in Proposition \ref{aux1}. Then, as $\e\rightarrow 0^+$, we have
$$|\bar\phi_1|_\infty =o(\e^{-\frac{N-2}{2(N-4)}}), $$
uniformly with respect to $d_1$ satisfying \eqref{limdj} for $j=1$.
\end{lemma}

\begin{proof}
Let us fix a small $\eta>0$ and remember that $\delta_1=\epsilon^{\frac{1}{N-4}}d_1$ (see \eqref{deltaj}), with $d_1$ satisfying \eqref{limdj} for $j=1$. We observe that by definition, since  $\bar\phi_1\in\mathcal
K_1^\bot$ solves the first equation of (\ref{sistemaaux}), then, for
all $\epsilon>0$ sufficiently small, there exists a constant
$c_\epsilon$ (which depends also on $d_1$) such that $\bar\phi_1$ weakly solves
\begin{equation}\label{eqdebolec}
-\Delta \bar\phi_1 = \epsilon \bar\phi_1 +\e \p\U_{\delta_1}+
f(\p\U_{\delta_1}+ \bar\phi_1) -f(\U_{\delta_1})- c_\epsilon \Delta\p Z_1.
\end{equation}
Testing \eqref{eqdebolec} with $\p Z_1$, taking into account that $\bar\phi_1\in\mathcal K_1^\bot$ and the definition of $\p Z_1$, we have that
\begin{equation}\label{eqprimoresto}
\begin{array}{lll}
\displaystyle c_\e  \int_\O p U_{\delta_1}^{p-1} \p Z_1 Z_1 \, dx &=&\displaystyle -\e \int_\O \bar\phi_1 \p Z_1\, dx
-\e\int_\O \p\U_{\delta_1}\p Z_1\, dx
-\int_\O\left[f(\p\U_{\delta_1})-f(\U_{\delta_1})\right]\p Z_1\, dx\\
&&\displaystyle -\int_\O\left[f(\p\U_{\delta_1}+\bar\phi_1)-f(\p\U_{\delta_1})\right]\p Z_1 \,
dx.
\end{array}
\end{equation}

By definition, if we set $\psi:=Z_1-PZ_1$, then $\psi$ is an harmonic function and $\psi=Z_1$ on $\partial\Omega$, therefore, by elementary elliptic estimates, for all sufficiently small $\e>0$, for any $d_1\in]\eta,\frac{1}{\eta}[$ we have that $|\psi|_{\infty,\Omega}\leq C \delta_1^{\frac{N-4}{2}}$, for some positive  constant $C=C(N,\Omega)$ depending only on $N$ and $\Omega$, and hence
$$ \int_\O p U_{\delta_1}^{p-1} \p Z_1 Z_1 \, dx = \int_\O p U_{\delta_1}^{p-1} Z_1^2 \, dx - \int_\O p U_{\delta_1}^{p-1} \psi Z_1 \, dx .$$
Now
\begin{eqnarray*}
\int_\O p U_{\delta_1}^{p-1} Z_1^2\, dx &=&
c_N \delta_1^{N-4} \delta_1^2 \int_\O\frac{1}{(\delta_1^2+|x|^2)^2}\frac{(|x|^2-\delta_1^2)^2}{(\delta_1^2+|x|^2)^N}\,
dx\\
&=& c_N \delta_1^{-2}\int_{\mathbb R^N}\frac{(|y|^2-1)^2}{(1+|y|^2)^{N+2}}\, dy +
O(\delta_1^{N-2})\\
&=& A_N \delta_1^{-2}+ o(1),\ \ \ \hbox{as} \ \epsilon\to 0.
\end{eqnarray*}

By using 
the property $|\psi|_{\infty,\Omega}\leq C \delta_1^{\frac{N-4}{2}}$, by the same computations, we see that
$$\int_\O p U_{\delta_1}^{p-1} \psi Z_1 \, dx=O(\delta_1^{N-4}), \ \ \hbox{as} \ \e\to 0.$$
Therefore, we get that
\begin{equation}\label{uguaglianzaimpoz1}
 \int_\O p U_{\delta_1}^{p-1} \p Z_1 Z_1 \, dx =A_N \delta_1^{-2}+ o(1),\ \ \ \hbox{as} \ \epsilon\to 0.
\end{equation}
Moreover, reasoning as before, we have
\begin{eqnarray*}
\int_\O Z_1^2\, dx &=&
c_N\delta_1^{N-4}\int_\O\frac{(|x|^2-\delta_1^2)^2}{(\delta_1^2+|x|^2)^N}\,
dx\\
&=& c_N \int_{\mathbb R^N}\frac{(|y|^2-1)^2}{(1+|y|^2)^N}\, dy +
O(\delta_1^{N-2})\\
&=& B_N + o(1), \ \ \hbox{as} \ \e\to 0,
\end{eqnarray*}

and, by an analogous computation

$$|Z_1|_{\frac{2N}{N-2}}\leq c_{N}\left[\int_\O \d_1^{\frac{N(N-4)}{N-2}}\frac{||x|^2-\d_1^2|^{\frac{2N}{N-2}}}{(\d_1^2+|x|^2)^{\frac{N^2}{N-2}}}\, dx\right]^{\frac{N-2}{2N}}\leq C_{N} \delta_1^{-1},$$

and hence, since $\p Z_1=Z_1 - \psi$, by elementary estimates, we get that for all sufficiently small $\epsilon>0$
\begin{equation}\label{stimanormeprojder}
|\p Z_1 |_2^2\leq 2B_N,\ \ \ |\p Z_1 |_{\frac{2N}{N-2}}\leq 2C_N \delta_1^{-1} .
\end{equation}
Thanks to \eqref{stimanormeprojder}, applying H\"older inequality, Poincar\`e inequality, taking into account of \eqref{numero}, the asymptotic expansion of $|\p\U_{\delta_1}|_2$ (see  Lemma \ref{lemma6} and its proof), the choice of $\delta_1$ (see \eqref{deltaj}) and since $\bar\phi_1 \in B_{1,\epsilon}$, we have the following inequalities

$$\e\int_\O |\bar\phi_1||\p Z_1|\, dx \leq \e
|\bar\phi_1|_2|\p Z_1|_2\leq c_1\e \|\bar\phi_1\||\p Z_1|_2\leq c_2 \e^{\frac{\theta_1}{2}+1+\sigma}$$
$$\e\int_\O
\p\U_{\delta_1}|\p Z_1|\, dx \leq \e |\p\U_{\delta_1}|_2|\p Z_1|_2 \leq c
\e\delta_1,
$$
$$\int_\O |f(\p\U_{\delta_1})-f(\U_{\delta_1})||\p Z_1|\, dx \leq
|f(\p\U_{\delta_1})-f(\U_{\delta_1})|_{\frac{2N}{N+2}}|\p Z_1|_{\frac{2N}{N-2}}\leq
c\delta_1^{\frac{N+2}{2}}\delta_1^{-1}=c \d_1^{\frac N 2}.$$

Moreover, taking into account of Lemma \ref{lem1n2} and Sobolev inequality, we get that
\begin{eqnarray*}
&&\int_\O |f(\p\U_{\delta_1}+\bar\phi_1)-f(\p\U_{\delta_1})||\p Z_1|\, dx\\
&\leq &
|f(\p\U_{\delta_1}+\bar\phi_1)-f(\p\U_{\delta_1})|_{\frac{2N}{N+2}}|\p Z_1|_{\frac{2N}{N-2}}\\
&\leq& |f(\p\U_{\delta_1}+\bar\phi_1)-f(\p\U_{\delta_1})-f'(\p\U_{\delta_1})\bar\phi_1|_{\frac{2N}{N+2}}|\p Z_1|_{\frac{2N}{N-2}}+|f'(\p\U_{\delta_1})\bar\phi_1|_{\frac{2N}{N+2}}|\p Z_1|_{\frac{2N}{N-2}}\\
&\leq &
c\big||\bar\phi_1|^p\big|_{\frac{2N}{N+2}}|\p Z_1|_{\frac{2N}{N-2}}+|f'(\p\U_{\delta_1})\bar\phi_1|_{\frac{2N}{N+2}}|\p Z_1|_{\frac{2N}{N-2}}\\
&\leq &
c_1\left(|\bar\phi_1|_{\frac{2N}{N-2}}^{\frac{N+2}{N-2}}|\p Z_1|_{\frac{2N}{N-2}}+ |\p\U_{\delta_1}|_{\frac{2N}{N-2}}^{\frac{4}{N-2}}|\bar\phi_1|_{\frac{2N}{N-2}}|\p Z_1|_{\frac{2N}{N-2}}\right)
\\
&\leq &
c_2\left(\|\bar\phi_1\|^{\frac{N+2}{N-2}}|\p Z_1|_{\frac{2N}{N-2}}+ |\p\U_{\delta_1}|_{\frac{2N}{N-2}}^{\frac{4}{N-2}}\|\bar\phi_1\||\p Z_1|_{\frac{2N}{N-2}}\right)
\\
&\leq & c_3\e^{\frac{\theta_1}{2}+\sigma}\delta_1^{-1}= c_4 \e^{\frac{N-2}{2(N-4)}+\sigma-\frac{1}{N-4}} \leq c_4\e^{\frac 12}.
\end{eqnarray*}

Thus, from \eqref{eqprimoresto}, \eqref{uguaglianzaimpoz1} and the previous estimates, we get that for all sufficiently small $\epsilon>0$
\begin{equation}\label{stimaceps}
\begin{array}{lll}
\displaystyle |c_\e |&\leq &\displaystyle \frac{1}{A\delta_1^{-2}+o(1)}\left[\left|\e \int_\O \bar\phi_1
\p Z_1\,
dx\right|+\left|\e\int_\O \p\U_{\delta_1}\p Z_1\, dx\right|\right.\\[16pt]
&&\displaystyle\left.+\left|\int_\O\left[f(\p\U_{\delta_1})-f(\U_{\delta_1})\right]\p Z_1\,
dx\right|+\left|\int_\O\left[f(\p\U_{\delta_1}+\bar\phi_1)-f(\p\U_{\delta_1})\right]\,
dx\right|\right]\\[12pt]
&&\displaystyle \leq c\e^{\frac{2}{N-4}+\frac{1}{2}},
\end{array}
\end{equation}
uniformly with respect to $d_1$ satisfying $\eta<d_1<\frac{1}{\eta}$.

We observe that $\bar\phi_1$ is a classical solution of \eqref{eqdebolec}. This comes from the fact that $\bar\phi_1 \in H_0^1(\Omega)$ weakly solves \eqref{eqdebolec}, taking into account the smoothness of $\p\U_{\delta_1}$, $\U_{\delta_1}$, $\p Z_1$, from standard elliptic regularity theory and the application of a well-known lemma by Brezis and Kato.\\

We consider the quantity $\sup_{d_1 \in ]\eta,\frac{1}{\eta}[}\left(\frac{|\bar\phi_1|_\infty}{|\U_{\delta_1}|_\infty}\right)$, which is defined for all $\epsilon \in (0,\e_0)$, where $\epsilon_0>0$ is given by Proposition \ref{aux1}.
We want to prove that
\begin{equation}\label{sefunge}
 \displaystyle \lim_{\e\to0^+} \sup_{d_1 \in ]\eta,\frac{1}{\eta}[}\left(\frac{|\bar\phi_1|_\infty}{|\U_{\delta_1}|_\infty}\right)=
0.
\end{equation}
It is clear that \eqref{sefunge} implies the thesis. In fact, we recall that, thanks to the definition \eqref{Udelta} and the choice of $\delta_1$ (see \eqref{deltaj}), 
 for any $d_1 \in]\eta,\frac{1}{\eta}[$ we have $${\alpha_N} \eta^{\frac{N-2}{2}} \epsilon^{-\frac{N-2}{2(N-4)}}<  |\U_{\delta_1}|_{\infty}< {\alpha_N} \eta^{-\frac{N-2}{2}} \epsilon^{-\frac{N-2}{2(N-4)}}.$$ Hence, by this estimate and \eqref{sefunge}, we get that
$$0\leq \sup_{d_1 \in ]\eta,\frac{1}{\eta}[}\frac{|\bar\phi_1|_\infty}{\epsilon^{-\frac{N-2}{2(N-4)}}}=\sup_{d_1 \in ]\eta,\frac{1}{\eta}[}\left(\frac{|\bar\phi_1|_\infty}{|\U_{\delta_1}|_\infty} \cdot\frac{|\U_{\delta_1}|_\infty}{\epsilon^{-\frac{N-2}{2(N-4)}}}\right) \leq \sup_{d_1 \in ]\eta,\frac{1}{\eta}[}\left(\frac{|\bar\phi_1|_\infty}{|\U_{\delta_1}|_\infty}\right)  {\alpha_N} \eta^{-\frac{N-2}{2}}\rightarrow 0, $$
 as $\epsilon\rightarrow 0^+$, and we are done.\\

In order to prove \eqref{sefunge} we argue by contradiction. Assume that \eqref{sefunge} is false. Then, there exists a positive number $ \tau \in \mathbb{R}^+$, a sequence $(\e_k)_k \subset \mathbb{R}^+$, $\e_k\to0$ as $k\rightarrow +\infty$, such that
\begin{equation}\label{sefunge2}
 \sup_{d_1 \in ]\eta,\frac{1}{\eta}[}\left(\frac{|\bar\phi_{1,k}|_\infty}{|\U_{\delta_{1,k}}|_\infty}\right)> \tau,
\end{equation}
for any $k\in \mathbb{N}$, where, $\bar\phi_{1,k}:=\bar\phi_{1}(\epsilon_k,d_1) \in B_{1,\epsilon_k}$ and $\delta_{1,k}:=\epsilon_k^{\frac{1}{N-4}}d_1$. We observe that \eqref{sefunge2} contemplates the possibility that $ \sup_{d_1 \in ]\eta,\frac{1}{\eta}[}\left(\frac{|\bar\phi_{1,k}|_\infty}{|\U_{\delta_{1,k}}|_\infty}\right)=+\infty$. 
From (\ref{sefunge2}), for any $k \in \mathbb{N}$, thanks to the definition of $\sup$, we get that there exists $d_{1,k} \in ]\eta,\frac{1}{\eta}[$ such that
$$\left(\frac{|\bar\phi_{1,k}|_\infty}{|\U_{\delta_{1,k}}|_\infty}\right) (d_{1,k})> \frac{\tau}{2}.$$
Hence, if we consider the sequence $\left(\frac{|\bar\phi_{1,k}|_\infty}{|\U_{\delta_{1,k}}|_\infty} (d_{1,k})\right)_k$, then, up to a subsequence, as $k\rightarrow +\infty$, there are only two possibilities:
\begin{description}
\item[(a)] $ \frac{|\bar\phi_{1,k}|_\infty}{|\U_{\delta_{1,k}}|_\infty} (d_{1,k})\to +\infty$;
\item[(b)] $ \frac{|\bar\phi_{1,k}|_\infty}{|\U_{\delta_{1,k}}|_\infty} (d_{1,k}) \to l $, for some $l \geq  \frac{\tau}{2}>0$.
\end{description}
We will show that (a) and (b) cannot happen.\\\\\
Assume (a). 
We point out that, since $\eta>0$ is fixed, then, $d_{1,k} \in ]\eta,\frac{1}{\eta}[$ for all $k$, in particular this sequence stays definitely away from $0$ and from $+\infty$. Hence, in order to simplify the notation of this proof, we omit the dependence from $d_{1,k}$ in $\bar\phi_{1,k}(d_{1,k})$ and in $\delta_{1,k}(d_{1,k})=\e_k^{\frac{1}{N-4}}d_{1,k}$ and thus we simply write $\bar\phi_{1,k}$, $\delta_{1,k}$. In particular, we observe that, for any fixed $k$, $\bar\phi_{1,k}$ is a function depending only on the space variable $x \in \Omega$.

Then, for any $k\in\mathbb{N}$, let $a_k\in \Omega$ such that
$|\bar\phi_{1,k}(a_k)|=|\bar\phi_{1,k}|_\infty$ and set
$M_k:=|\bar\phi_{1,k}|_\infty$. Thanks to the assumption (a), since $|\U_{\delta_{1,k}}|_{\infty}={\alpha_N}\delta_{1,k}^{-\frac{N-2}{2}}={\alpha_N} \e_k^{-\frac{N-2}{2(N-4)}}d_{1,k}^{-\frac{N-2}{2}}$, we get that $M_k \rightarrow +\infty$, as $k \to +\infty$. We consider the rescaled function
$$\widetilde\phi_{1,k}(y):=
\frac{1}{M_k}\bar\phi_{1,k}\left(a_k +
\frac{y}{M_k^{\beta}}\right),\qquad \beta=\frac{2}{N-2}$$
defined for $y \in \widetilde
\Omega_k:=M_k^{\frac{2}{N-2}}(\Omega-a_k)$.
Moreover let us set $$\widetilde{\p\U}_{1,k}(y):=\frac{1}{M_k}
\p\U_{\delta_{1,k}} \left(a_k  +
\frac{y}{M_k^{\beta}}\right); \quad \widetilde{\U}_{1,k}(y):=\frac{1}{M_k}
\U_{\delta_{1,k}} \left(a_k  +
\frac{y}{M_k^{\beta}}\right);$$
$$ \ \ \widehat{\p Z}_{1,k}(y):=\frac{1}{M_k^{2\beta+1}}
\p Z_{1,k} \left(a_k  +
\frac{y}{M_k^{\beta}}\right).$$ Since we are assuming (a) it is clear that
$|\widetilde{\p\U}_{1,k}|_{\infty,\widetilde\Omega_k},  |\widetilde{\U}_{1,k}|_{\infty,\widetilde\Omega_k}  \rightarrow
0$, as $k\to +\infty$. Moreover, thanks to the definition of $Z_1$, and since $\p Z_1=Z_1-\psi$, with $|\psi|_{\infty,\Omega}\leq C \delta_1^{\frac{N-4}{2}}$, we have that $|PZ_{1,k}|_\infty \simeq |Z_{1,k}|_\infty \simeq \delta_{1,k}^{-\frac{N}{2}}$, and hence, thanks to (a), we have $\frac{1}{M_k^{2\beta+1}}=o( \delta_{1,k}^{\frac{N+2}{2}})$, which implies that $ |\widehat{\p Z}_{1,k}|_{\infty,\widetilde\Omega_k}  \rightarrow
0$, as $k\to +\infty$. In particular, thanks to \eqref{stimaceps}, the same conclusion holds for $c_{\epsilon_k}(d_{1,k}) \widehat{\p Z}_{1,k}$. 
Taking into account that $2\beta+1=p$, by elementary computations, we see
that $\widetilde\phi_{1,k}$ solves

\begin{equation}\label{PSBN}
\left\{
\begin{array}{lr}
-\Delta  \widetilde\phi_{1,k} =
\frac{\epsilon_k}{M_k^{2\beta}}  \widetilde\phi_{1,k}
+\e_k \frac{ \widetilde{\p\U_{\delta_{1,k}}}}{M_k^{2\beta}}+
f(\widetilde{\p\U}_{\delta_{1,k}}+\widetilde\phi_{1,k}
)-f(\widetilde{\U}_{\delta_{1,k}}) + c_{\epsilon_k}(d_{1,k}) \widehat{\p Z}_{1,k} \qquad \ \mbox{in}\,\, \widetilde{\O}_k,\\\\
\widetilde\phi_{1,k}=0\qquad\qquad\qquad\qquad\qquad\qquad\qquad\qquad\qquad\qquad\qquad\qquad\qquad\mbox{on}\,\,
\partial\widetilde\Omega_k.
\end{array}
\right.
\end{equation}

Let us denote by $\Pi$ the limit domain of $\widetilde\Omega_k$.
Since $M_k \rightarrow +\infty$, as $k \to +\infty$, we have that $\Pi$ is the whole $\R$ or an half-space. Moreover, since the family $(\widetilde\phi_{1,k})_k$ is
uniformly bounded and solves \eqref{PSBN}, then, by the same proof of Lemma 2.2 of \cite{Ben1}, we get that $0
\in \Pi$ (in particular $0 \notin \partial \Pi$), and, by standard
elliptic theory, it follows that, up to a subsequence,  as $k\rightarrow + \infty$, we have that
$\widetilde\phi_{1,k}$ converges in $C_{loc}^2(\Pi)$ to a
function $w$ which
satisfies

\begin{equation}\label{scplc}
-\Delta w =f(w) \ \hbox{in} \ \Pi, \ \ \ w(0)=1\ (\hbox{or} \ w(0)=-1),\ \ \ |w|\ \leq 1 \ \hbox{in}\ \Pi, \ \ \ w=0 \
\hbox{on} \ \partial \Pi.
\end{equation}

We observe that, thanks to the definition of the chosen rescaling, by elementary computations (see Lemma 2 of \cite{Iacopetti}), it holds $\|\widetilde
\phi_{1,k}\|_{\widetilde\Omega_\epsilon}^2 =
\|\bar\phi_{1,k}\|_\Omega^2$. Now, since $\|\bar\phi_{1,k}\|\leq c \e_k^{\frac{\theta_1}{2}+\sigma}$, where $c$ depends only on $\eta$ and $\sigma$ is some positive number (see Proposition \ref{aux1}), we have $\|\widetilde \phi_{1,k}\|_{\widetilde\Omega_k}^2 = \|\bar\phi_{1,k}\|_\Omega^2\rightarrow 0$, as
$k\rightarrow +\infty$. Hence, since $\widetilde\phi_{1,k}\rightarrow w$ in $C_{loc}^2(\Pi)$, by Fatou's lemma, it follows that

\begin{equation}\label{stimaenwlim}
\|w\|_\Pi^2 \leq \liminf_{k\rightarrow +\infty} \|\widetilde \phi_{1,k}\|_{\widetilde\Omega_k}^2=0.
\end{equation}

Therefore, since $\|w\|_\Pi^2=0$ and $w$ is smooth, it follows that $w$ is constant, and from $w(0)=1$ (or $w(0)=-1$) we get that $w\equiv1$ (or $w\equiv -1$) in $\Pi$. But, since $w$ is constant and solves $-\Delta w =f(w)$ in $\Pi$, then necessarily $f(w)\equiv0$ in $\Pi$, and hence $w$ must be the null function, but this contradicts $w\equiv1$ (or $w\equiv -1$).

Alternatively, if $\Pi$ is an half-space, by using the boundary condition $w=0$ on $\partial \Pi$, we contradicts $w\equiv1$ (or $w\equiv -1$). Hence, the only possibility is $\Pi=\R$. In this case, since $w$ solves \eqref{scplc} and $\|w\|_\Pi^2 \leq 2S^{N/2}$, it is well known that $w$ cannot be sign-changing and hence, assuming without loss of generality that $w(0)=1$, $w$ must be a positive function of the form $\U_{\delta_N}$ (see \eqref{Udelta}), for some $\delta_N$ such that $\U_{\delta_N}(0)=1$, and this contradicts $w\equiv 1$.\\ Hence (a) cannot happen.\\\\
Assume (b). Using the same convention on the notation as in previous case, we deduce 
that there exist two positive uniform constants $c_1$, $c_2$ such that
\begin{equation}\label{uniformmaggcost}
c_1\delta_{1,k}^{-\frac{N-2}{2}}\leq |\bar\phi_{1,k}|_\infty\leq c_2 \delta_{1,k}^{-\frac{N-2}{2}},
\end{equation}
for all sufficiently large $k$. In particular, it still holds that $M_k \rightarrow +\infty$, as $k \to +\infty$.
We consider the same rescaled functions $\widetilde{\phi}_{1, k}$ as in (a) and, as before, we denote by $\Pi$ the limit domain of
$\widetilde{\O}_k$.\\
Now, up to a subsequence, since $\widetilde{\p\U}_{1,k}$ and
$\widetilde{\U}_{1,k}$ are uniformly bounded we see that they converge in $C^2_{loc}(\Pi)$ to a bounded
function which we denote, respectively, by $\overline{\p\U}$ and
$\overline\U$ (one of them or both could be eventually the null function).  In fact $\widetilde{\U}_{1,k}$ is uniformly bounded and solves $-\Delta \widetilde{\U}_{1,k}=\widetilde{\U}_{1,k}^p$ on  $ \widetilde\Omega_k$, and so by standard elliptic theory we get that $\widetilde{\U}_{1,k}$ converges in $C_{loc}^2(\Pi)$ to some non-negative bounded function $\overline\U$ which solves $-\Delta\overline\U=\overline\U^p$ in $\Pi$. Now, taking into account that $\widetilde{\U}_{1,k}\rightarrow \overline \U$ in $C_{loc}^2(\Pi)$, the same argument applies to $\widetilde{\p\U_{1,k}}$, which solves
\begin{equation*}
\left\{
\begin{array}{lr}
-\Delta\widetilde{\p\U}_{1,k}=\widetilde{\U}_{1,k}^p\qquad \mbox{in}\,\,\ \widetilde\Omega_k,\\
\widetilde{\p\U}_{1,k}=0\qquad\qquad\mbox{on}\,\,\,\partial\widetilde\Omega_k,
\end{array}
\right.
\end{equation*}
and hence $\widetilde{\p\U}_{1,k}$ converges in $C_{loc}^2(\Pi)$ to some non-negative bounded function $\overline{\p\U}$ satisfying $-\Delta \overline{\p\U}=\overline \U^p$ in $\Pi$, $\overline{\p\U}=0$ on $\partial\Pi$.

We point out that as in (a), but using \eqref{uniformmaggcost}, we still have $c_{\e_k}(d_{1,k}) |\widehat{\p Z}_{1,k}|_{\infty,\widetilde\Omega_k}  \rightarrow
0$, as $k \to +\infty$. Moreover, by the proof of Lemma 2.2 of \cite{Ben1}, it also holds that $0\in \Pi$.
\\\\ Hence, by standard elliptic theory, we have that $\widetilde{\phi}_{1, k}$
converges in $C^2_{loc}(\Pi)$ to a function $w$ which solves
\begin{equation}\label{PSBN1}
\left\{
\begin{array}{lr}
-\Delta  w = f(\overline{\p\U}+w
)-f(\overline{\U})\qquad \mbox{in}\,\, \Pi,\\\\
w=0\ \ \qquad\qquad\qquad\qquad\qquad\mbox{on}\,\,
\partial\Pi,\\\\
w(0)=1 \ (\hbox{or} \ w(0)=-1 ).
\end{array}
\right.
\end{equation}

As in \eqref{stimaenwlim} we have $\|w\|_\Pi^2=0$
 and hence, since $w$ is smooth, the only possibility is $w\equiv 1$ (or $w\equiv -1)$ because of the condition $w(0)=1$ (or $w(0)=-1$). 
 Moreover, thanks to the definition of the chosen rescaling, it also holds $|\widetilde{\phi}_{1, k}|_{\frac{2N}{N-2},\widetilde \Omega_k}=|\bar{\phi}_{1,k}|_{\frac{2N}{N-2}, \Omega}$ (for the proof see Lemma 2 of \cite{Iacopetti}). Therefore, since  $|\bar{\phi}_{1,k}|_{\frac{2N}{N-2},\Omega} \rightarrow 0$ (because $\|\bar\phi_{1,k}\|\leq c \e_k^{\frac{\theta_1}{2}+\sigma}$, where $c>0$ depends only on $\eta$) and $\widetilde{\phi}_{1, k}\rightarrow w$ in $C_{loc}^2(\Pi)$, as $k \to +\infty$, then, by Fatou's Lemma, it follows that $|w|_{\frac{2N}{N-2},\Pi}=0$, and thus it cannot happen that $w\equiv 1$ (or $w\equiv -1$).\\
 Hence (a) and (b) cannot happen, and the proof is then concluded.
\end{proof}

We are now in position to prove Proposition \ref{auxsolving}.

\begin{proof}[Proof of Proposition \ref{auxsolving}]
Let us fix $\eta>0$ and let $\bar\phi_1 \in \mathcal{K}_1^\perp \cap B_{1,\epsilon}$ be the unique solution of the first equation of \eqref{sistemaaux} found in Proposition \ref{aux1}. We define the operator $\mathcal{T}_2 :\mathcal{K}^{\perp} \rightarrow \mathcal{K}^{\perp} $ as $$\mathcal{T}_2(\phi_2):=- \mathcal{L}_2^{-1}[\mathcal{N}_2(\bar\phi_1,\phi_2)+\mathcal{R}_2].$$\\
In order to find a solution of the second equation of \eqref{sistemaaux} we solve the fixed point problem $\mathcal{T}_2(\phi_2)=\phi_2$.  Let us define the proper ball

$$\displaystyle B_{2,\epsilon}:=\{\phi_2 \in \mathcal{K}^{\perp}; \ \|\phi_2\| \leq r\ \epsilon^{\frac{\theta_2}{2}+\sigma} \}$$
for $r>0$ sufficiently large and $\sigma>0$ to be chosen later.\\
From Lemma \ref{lem1invl2}, there exists $\epsilon_0=\epsilon_0(\eta)>0$ and $c=c(\eta)>0$ such that:
\begin{equation}\label{eq1n2}
\|\mathcal{T}_2(\phi_2)\| \leq c (\|\mathcal{N}_2(\bar\phi_1,\phi_2)\| + \|\mathcal{R}_2\|),
\end{equation}
and
\begin{equation}\label{eq2n2}
\|\mathcal{T}_2(\phi_2)-\mathcal{T}_2(\psi_2)\| \leq c (\|\mathcal{N}_2(\bar\phi_1, \phi_2)-\mathcal{N}_2(\bar\phi_1, \psi_2)\|),
\end{equation}
for all $\phi_2, \psi_2 \in \mathcal{K}^\perp$, for all $(d_1,d_2) \in \mathbb R_+^2$  satisfying \eqref{limdj} and for all $\epsilon \in (0,\epsilon_0)$.

We begin with estimating the right hand side of (\ref{eq1n2}).\\ Thanks to Proposition \ref{errorprop2}  we have that
$$\| \mathcal{R}_2\| \leq c \epsilon^{\frac{\theta_2}{2}+\sigma},$$
for all $\epsilon \in (0,\epsilon_0)$, for all $(d_1,d_2) \in \mathbb R_+^2$ satisfying \eqref{limdj}. Thus it remains only to estimate $\|\mathcal{N}_2(\bar\phi_1,\phi_2)\|$. Thanks to \eqref{stimaistar} and the definition of $\mathcal{N}_2$ we deduce:
\begin{equation}\label{eq3n2}
\|\mathcal{N}_2(\bar\phi_1,\phi_2)\| \leq c | f(V_\e+\bar\phi_1+\phi_2) - f(V_\e) - f^\prime(V_\e) \phi_2-f(\p\U_{\delta_1}+\bar\phi_1)+f(\p\U_{\delta_1}) |_{\frac{2N}{N+2}}.
\end{equation}
We estimate the right-hand side of (\ref{eq3n2}): 
\begin{equation*}
\begin{array}{lll}
&& | f(V_\e+\bar\phi_1+\phi_2) - f(V_\e) - f^\prime(V_\e) \phi_2 -f(\p\U_{\delta_1}+\bar\phi_1)+f(\p\U_{\delta_1}) |_{\frac{2N}{N+2}} \\[10pt]
&\leq&   | f(V_\e+\bar\phi_1+\phi_2) - f(V_\e+\bar\phi_1) - f^\prime(V_\e+\bar\phi_1)  \phi_2 |_{\frac{2N}{N+2}}+|  (f^\prime(V_\e+\bar\phi_1) - f^\prime(V_\e) )\phi_2  |_{\frac{2N}{N+2}}\\[10pt]
&&+ |f(V_\e+\bar\phi_1)-f(V_\e)-f(\p\U_{\delta_1}+\bar\phi_1)+ f(\p\U_{\delta_1})|_{\frac{2N}{N+2}}
\end{array}
\end{equation*}
In order to estimate the last three terms, by  Lemma \ref{lem0n2} and Lemma \ref{lem1n2} we deduce that:
\begin{equation} \label{eq4bisn2}
|f(V_\e+\bar\phi_1+\phi_2) - f(V_\e+\bar\phi_1) -f'(V_\e+\bar\phi_1) \phi_2| \leq  c |\phi_2|^p
\end{equation}
and
\begin{equation}\label{eq4bisn21}
|(f'(V_\e+\bar\phi_1)  - f'(V_\e))\phi_2|\leq  c |\bar\phi_1|^{p-1}|\phi_2|.
\end{equation}
Since $\frac{2N}{N+2} \cdot p = p+1$ we get that
\begin{equation*}
\begin{array}{lll}
\displaystyle \int_{\O}|f(V_\e+\bar\phi_1+\phi_2) - f(V_\e+\bar\phi_1) -f^\prime(V_\e+\bar\phi_1) \phi_2|^{\frac{2N}{N+2}} \ dx &\leq & \displaystyle c  \int_{\O} |\phi_2|^{p+1} \ dx,
 \\[12pt]
\end{array}
\end{equation*}
and applying Sobolev inequality we deduce that
\begin{equation}\label{eq4n2}
| f(V_\e+\bar\phi_1+\phi_2) - f(V_\e+\bar\phi_1) -f^\prime(V_\e+\bar\phi_1) \phi_2|_{\frac{2N}{N+2}} \leq c \|\phi_2\|^p.
\end{equation}
By (\ref{eq4bisn21}) we get that
\begin{equation*}
\begin{array}{lll}
\displaystyle \int_{\O}| (f^\prime(V_\e+\bar\phi_1)  - f^\prime(V_\e))\phi_2|^{\frac{2N}{N+2}} \ dx &\leq & \displaystyle c  \int_{\O} |\bar\phi_1|^{(p-1)\frac{2N}{N+2}}|\phi_2|^{\frac{2N}{N+2}} \ dx.
 \\[12pt]
\end{array}
\end{equation*}
We observe that $\phi_1^{(p-1)\frac{2N}{N+2}} \in L^{\frac{N+2}{4}}$, $\phi_2^{\frac{2N}{N+2}} \in L^p $ and $p$, $\frac{N+2}{4}$ are conjugate exponents in H\"older inequality. Moreover $(p-1)\frac{2N}{N+2} \frac{N+2}{4} = p+1$ so
\begin{equation*}
| (f^\prime(V_\e+\bar\phi_1)  - f^\prime(V_\e))\phi_2|_{\frac{2N}{N+2}}^{\frac{2N}{N+2}} \leq c |\bar\phi_1|_{p+1}^{\frac{8N}{(N+2)(N-2)}} |\phi_2|_{p+1}^{\frac{2N}{N+2}},
\end{equation*}
and hence by Sobolev inequality we deduce that
\begin{equation}\label{eq5n2}
| (f^\prime(V_\e+\bar\phi_1)  - f^\prime(V_\e))\phi_2|_{\frac{2N}{N+2}} \leq c \|\bar\phi_1\|^{\frac{4}{N-2}} \|\phi_2\|.
\end{equation}
It remains to estimate the last term. As in the proof of Proposition \ref{errorprop2} we make the decomposition of the domain $\Omega$ as $\Omega=A_0 \sqcup A_1 \sqcup A_2$. Hence we get that:

\begin{eqnarray*}
 |f(V_\e+\bar\phi_1)-f(V_\e)-f(\p\U_{\delta_1}+\bar\phi_1)+f(\p\U_{\delta_1})|_{\frac{2N}{N+2}, A_0}&\leq& |f(V_\e+\bar\phi_1)-f(\p\U_{\delta_1}+\bar\phi_1)|_{\frac{2N}{N+2}, A_0}\\&&
 +|f(V_\e)-f(\p\U_{\delta_1})|_{\frac{2N}{N+2}, A_0}
\end{eqnarray*}

Then, by using the definition of $\delta_1, \delta_2$, the usual elementary inequalities, the computations made in \eqref{elemcompututil} and Sobolev inequality, we get that
\begin{eqnarray*}
|f(V_\e+\bar\phi_1)-f(\p\U_{\delta_1}+\bar\phi_1)|_{\frac{2N}{N+2}, A_0}&\leq &c_1\left( |\p\U_{\delta_2}|_{p+1, A_0}^p+|\p\U_{\delta_1}^{p-1}\p\U_{\delta_2}|_{\frac{2N}{N+2}, A_0}+\Big||\bar\phi_1|^{p-1}\p\U_{\delta_2}\Big|_{\frac{2N}{N+2}, A_0}\right)\\
&\leq & c_2\left(\delta_2^{\frac{N+2}{2}}+\delta_1^2\delta_2^{\frac{N-2}{2}}+\|\bar\phi_1\|^{p-1}\delta_2^{\frac{N-2}{2}}\right)\\
&\leq & c_3 \e^{\frac{\theta_2}{2}+\sigma},
\end{eqnarray*}
for some $\sigma>0$.\\\\
Moreover, as in the previous estimate, we get that
\begin{eqnarray*}
|f(V_\e)-f(\p\U_{\delta_1})|_{\frac{2N}{N+2}, A_0}&\leq &c_1\left( |\p\U_{\delta_1}^{p-1}\p\U_{\delta_2}|_{\frac{2N}{N+2}, A_0}+|\p\U_{\delta_2}|^p_{\frac{2N}{N+2}, A_0}\right)\\
&\leq & c_2\e^{\frac{\theta_2}{2}+\sigma}.
\end{eqnarray*}
%
In $A_1$ we argue as in the previous case. The various terms now can be estimated as done in \eqref{eq2I1A1} and \eqref{eq3I1A1} and hence the same conclusion holds.\\

 For $A_2$, by using the usual elementary inequalities, Lemma \ref{techlembehvphi1} and remembering the choice of $\delta_1$, $\delta_2$, we have:
\begin{eqnarray*}
&& |f(V_\e+\bar\phi_1)-f(V_\e)-f(\p\U_{\delta_1}+\bar\phi_1)+f(\p\U_{\delta_1})|_{\frac{2N}{N+2}, A_2}\\
&\leq& |f(V_\e+\bar\phi_1)-f(V_\e)-f'(V_\e)\bar\phi_1|_{\frac{2N}{N+2}, A_2}+|f(\p\U_{\delta_1}+\bar\phi_1)-f(\p\U_{\d_1})-f'(\p\U_{\d_1})\bar\phi_1|_{\frac{2N}{N+2}, A_2}\\
&&+|[f'(V_\e)-f'(\p\U_{\delta_1})]\bar\phi_1|_{\frac{2N}{N+2}, A_2}\\
&\leq&  c \Big||\bar\phi_1|^p\Big|_{\frac{2N}{N+2}, A_2}+c|\p\U_{\delta_2}^{p-1}\bar\phi_1|_{\frac{2N}{N+2}, A_2}\\
& \leq& c|\bar\phi_1|_\infty^p \left(\int_{A_2}1\, dx \right)^{\frac{N+2}{2N}}+ c |\bar\phi_1|_\infty\left(\int_{A_2}\U_{\delta_2}^{\frac{8N}{N^2-4}}\, dx\right)^{\frac{N+2}{2N}}\\
& \leq&c_1 \delta_1^{-\frac{N-2}{2}p}\left(\int_0^{\sqrt{\delta_1\delta_2}}r^{N-1}\, dr\right)^{\frac{N+2}{2N}} + c_2 |\bar\phi_1|_\infty\left(\int_{A_2} \frac{\delta_2^{\frac{4N}{N+2}}}{(\delta_2^2+|x|^2)^{\frac{4N}{N+2}}}\, dx\right)^{\frac{N+2}{2N}}\\
& \leq&c_3 \delta_1^{-\frac{N+2}{2}} \left(\delta_1\delta_2\right)^{\frac{N+2}{4}} + c_2 |\bar\phi_1|_\infty \delta_2^2\left(\int_{A_2} \frac{1}{|x|^{\frac{8N}{N+2}}}\, dx\right)^{\frac{N+2}{2N}}\\
&\leq&  c_3\left(\frac{\delta_2}{\delta_1}\right)^{\frac{N+2}{4}}+ c_4 \delta_1^{-\frac{N-2}{2}}\delta_2^2\left(\int_0^{\sqrt{\delta_1\delta_2}}r^{\frac{N^2-7N-2}{N+2}}\, dr\right)^{\frac{N+2}{2N}}\\
&\leq&  c_3\left(\frac{\delta_2}{\delta_1}\right)^{\frac{N+2}{4}}+ c_5 \delta_1^{-\frac{N-2}{2}}\delta_2^2\left(\delta_1 \delta_2\right)^{\frac{N-6}{4}}
=  c_6\left(\frac{\delta_2}{\delta_1}\right)^{\frac{N+2}{4}} \leq  c_7 \e^{\frac{\theta_2}{2}+\sigma}.
 \end{eqnarray*}
Hence, from these estimates, we have
\begin{equation}\label{eqjoln2}
|f(V_\e+\bar\phi_1)-f(V_\e)-f(\p\U_{\delta_1}+\bar\phi_1)+
f(\p\U_{\delta_1})|_{\frac{2N}{N+2}}\leq c\e^{\frac{\theta_2}{2}+\sigma}.
\end{equation}

Since $\phi_2 \in B_{2,\epsilon}$ and thanks to (\ref{eq3n2}), (\ref{eq4n2}), (\ref{eq5n2}) and \eqref{eqjoln2} we get that $$\|\mathcal T_2(\phi_2)\|\leq c \e^{\frac{\theta_2}{2}+\sigma},\qquad \sigma>0$$ and hence $\mathcal{T}_2$ maps $B_{2,\epsilon}$ into itself .

It remains to prove that $\mathcal{T}_2:B_{2,\epsilon} \rightarrow B_{2,\epsilon} $ is a contraction. Thanks to (\ref{eq2n2}) it suffices to estimate $\|\mathcal{N}_2(\bar\phi_1, \phi_2)-\mathcal{N}_2(\bar\phi_1, \psi_2)\|$ for any $\psi_2, \phi_2 \in B_{2,\epsilon}$. To this end, thanks to \eqref{stimaistar}, the definition of $\mathcal{N}_2$ and reasoning as in the proof of Proposition \ref{aux1} we have:
\begin{equation*}
\|\mathcal{N}_2(\bar\phi_1,\phi_2)-\mathcal{N}_2(\bar\phi_1,\psi_2)\| \leq \e^\alpha \|\phi_2-\psi_2\|,
\end{equation*}
for some $\alpha>0$.\\\\ At the end we get that there exists $L\in (0,1)$ such that
$$ \|\mathcal{T}_2(\phi_2)- \mathcal{T}_2(\psi_2)\| \leq L \|\phi_2-\psi_2\|.$$
Finally, taking into account that $d_1 \rightarrow \bar\phi_1(d_1)$ is a $C^1$-map, a standard argument shows that also $(d_1,d_2) \rightarrow \bar\phi_2(d_1,d_2)$ is a $C^1$-map.
The proof is complete.
\end{proof}


\section{The reduced functional}\label{biforcazione}
We are left now to solve \eqref{bif}. Let $(\bar\phi_1, \bar\phi_2)\in \mathcal{K}^\bot_1\times\mathcal{K}^\bot$ be the solution found in Proposition \ref{auxsolving}. Hence $V_\e+\bar\phi_1+\bar\phi_2$ is a solution of our original problem \eqref{BN} if we can find $\bar d_\e=(\bar d_{1\e}, \bar d_{2\e})$  which satisfies condition \eqref{limdj} and solves equation \eqref{bif}.

 To this end we consider the reduced functional $\tilde J_\epsilon:\mathbb R_+^2 \rightarrow \mathbb R$ defined by: $$\tilde J_\epsilon(d_1,d_2):=J_\epsilon(V_\e+\bar\phi_1+\bar\phi_2),$$ where $J_\e$ is the functional defined in \eqref{funzionale}.\\
Our main goal is to show first that solving equation (\ref{bif}) is equivalent to finding critical points $(\bar d_{1,\epsilon},\bar d_{2,\epsilon})$ of the reduced functional $\tilde J_\epsilon(d_1,d_2)$ and then that the reduced functional has a critical point. These facts are stated in the following proposition:
\begin{proposition}\label{ridotto}
\begin{description}
The following facts hold:
\item[(i)] If $(\bar d_{1,\epsilon},\bar d_{2,\epsilon})$ is a critical point of $\tilde J_\epsilon$, then the function $V_\e+\bar\phi_1+\bar\phi_2$ is a solution of \eqref{BN}.
\item[(ii)] For any $\eta>0$, there exists $\epsilon_0>0$ such that for all $\epsilon \in (0,\epsilon_0)$ it holds:
\begin{eqnarray}\label{espfunzridotto}
\begin{array}{lll}
\displaystyle \tilde J_\epsilon(d_1,d_2)&=&\displaystyle \frac{2}{N} S^{N/2} + \epsilon^{\theta_1}\left[ a_1 \tau(0) d_1^{N-2} -a_2 d_1^2\right] + O(\epsilon^{\theta_1 + \sigma}),
\end{array}
\end{eqnarray}
with
\begin{equation}\label{espfunzridottopart2}
O(\epsilon^{\theta_1 + \sigma}) = \epsilon^{\theta_1+\sigma} g(d_1)+
\epsilon^{\theta_2}\left[ a_3 \tau(0)
\left(\frac{d_2}{d_1}\right)^{\frac{N-2}{2}} -a_2 d_2^2\right] +
o\left( \epsilon^{\theta_2 }\right),
\end{equation}
for some function $g$ depending only on $d_1$ (and uniformly bounded with respect to $\epsilon$), where $\theta_1, \theta_2$ are defined in \eqref{thetaj}, $\sigma$ is some positive real number (depending only on $N$), $\tau$ is the Robin's function of the domain $\O$ at the origin and $$a_1:=\frac 12 \alpha_N^{p+1}\int_{\R}\frac{1}{(1+|y|^2)^{\frac{N+2}{2}}}\, dy;\quad a_2:=\frac 12 \alpha_N^2\int_{\R}\frac{1}{(1+|y|^2)^{N-2}}\, dy;$$ $$\qquad a_3:=\alpha_N^{p+1}\int_{\R}\frac{1}{|y|^{N-2}(1+|y|^2)^{\frac{N+2}{2}}}\, dy.$$
The expansions  \eqref{espfunzridotto}, \eqref{espfunzridottopart2} are $C^0$-uniform with respect to $(d_1,d_2)$ satisfying condition \eqref{limdj}.
\end{description}
\end{proposition}
\begin{remark}\label{remarkespfunzrid}
We point out that the term $g$ appearing in \eqref{espfunzridottopart2} does not depend on $d_2$ and this will be used in the sequel, in particular in \eqref{eq21propcrit}.
\end{remark}
The aim of this section is to prove Proposition \ref{ridotto}. First we prove two lemmas about the $C^0$-expansion of the reduced functional $\tilde J_\epsilon (d_1,d_2):=J_\epsilon(V_\e+\bar\phi_1+\bar\phi_2)$, where $\bar\phi_1 \in \mathcal{K}^\perp_1\cap B_{1,\epsilon}$ and $\bar\phi_2 \in  \mathcal{K}^\perp \cap B_{2,\epsilon}$ are the functions given by Proposition \ref{auxsolving}.

\begin{lemma}\label{lem1exp1}
For any $\eta>0$ there exists $\epsilon_0>0$ such that for any $\epsilon \in (0,\epsilon_0)$ it holds:
$$J_\epsilon(V_\e+\bar\phi_1)=J_\epsilon(V_\e) +O (\epsilon^{\theta_1+ \sigma}), $$
with
\begin{equation}\label{esplem1exp1}
O(\epsilon^{\theta_1 + \sigma}) = \epsilon^{\theta_1+\sigma}
g_1(d_1)+  O\left( \epsilon^{\theta_2+\sigma }\right),
\end{equation}
for some function $g_1$ depending only on $d_1$ (and uniformly bounded with respect to $\epsilon$), where $\theta_1, \theta_2$ are defined in \eqref{thetaj}, $\sigma$ is some positive real number (depending only on $N$).
These expansion are $C^0$-uniform with respect to $(d_1,d_2)$ satisfying condition \eqref{limdj}.
\end{lemma}
\begin{proof}
Let us fix $\eta>0$. By direct computation we immediately see that
\begin{equation}\label{eq1exp1}
\begin{array}{lll}
J_\epsilon(V_\e+\bar\phi_1) - J_\epsilon(V_\e)&=&\frac{1}{2} \int_\Omega |\nabla \bar\phi_1|^2 \ dx + \int_\Omega \nabla V_\e \cdot \nabla \bar\phi_1 \ dx -\frac{\epsilon}{2} \int_\Omega |\bar\phi_1|^2 \ dx - \epsilon \int_\Omega V_\e \bar\phi_1 \ dx\\[10pt]
&&- \frac{1}{p+1} \int_\Omega (|V_\e+\bar\phi_1|^{p+1} -|V_\e|^{p+1} ) \ dx.
\end{array}
\end{equation}
By definition we have $$\int_\Omega \nabla V_\e \cdot \nabla \bar\phi_1 \ dx=  \int_\Omega \nabla (\p\U_{\delta_1}- \p\U_{\delta_2}) \cdot \nabla \bar\phi_1 \ dx =  \int_\Omega (\U_{\delta_1}^p - \U_{\delta_2}^p) \bar\phi_1 \ dx =  \int_\Omega [f(\U_{\delta_1}) - f(\U_{\delta_2})] \bar\phi_1 \ dx,$$
moreover, since $F(s)=\frac{1}{p+1} |s|^{p+1}$ is a primitive of $f$, we can write (\ref{eq1exp1}) as

\begin{equation}\label{eq2exp1}
\begin{array}{lll}
J_\epsilon(V_\e+\bar\phi_1) - J_\epsilon(V_\e)&=& \frac{1}{2} \|\bar\phi_1\|^2   -\frac{\epsilon}{2} |\bar\phi_1|_2^2  - \epsilon \int_\Omega V_\e \bar\phi_1 \ dx +  \int_\Omega [f(\U_{\delta_1}) - f(\U_{\delta_2})] \bar\phi_1 \ dx\\[10pt]
&&-  \int_\Omega [F(V_\e+\bar\phi_1) -F(V_\e)] \ dx\\[10pt]
&=& \frac{1}{2} \|\bar\phi_1\|^2   -\frac{\epsilon}{2} |\bar\phi_1|_2^2  - \epsilon \int_\Omega V_\e \bar\phi_1 \ dx +  \int_\Omega [f(\U_{\delta_1}) - f(\U_{\delta_2})-f(V_\e)] \bar\phi_1 \ dx\\[10pt]
&&-  \int_\Omega [F(V_\e+\bar\phi_1) -F(V_\e)-f(V_\e)\bar\phi_1] \ dx\\[10pt]
&&A+B+C+D+E.
\end{array}
\end{equation}
\textbf{A,B:} Thanks to Proposition \ref{auxsolving}, for all sufficiently small $\epsilon$, we have $\|\bar\phi_1\|\leq c \epsilon^{\frac{\theta_1}{2}+ \sigma}$, for some $c>0$ and for some $\sigma>0$ depending only on $N$. Hence we deduce that $A=O(\epsilon^{\theta_1+ 2 \sigma})$, $B=O(\epsilon^{\theta_1+ 2 \sigma+1})$. We point out that, since only $\bar\phi_1$ is involved in $A$ and $B$, these terms depend only on $d_1$.

\textbf{C:}  By definition we have
\begin{equation*}
\epsilon \int_\Omega V_\e \bar\phi_1 \ dx = \epsilon \int_\Omega \p\U_{\delta_1} \bar\phi_1 \ dx - \epsilon \int_\Omega \p\U_{\delta_2} \bar\phi_1 \ dx = I_1+I_2.
\end{equation*}
We observe that in the estimate $I_1$ only $\delta_1$ and $\bar\phi_1$ are involved. Hence $I_1$ depends only on $d_1$.
Thanks to H\"older inequality, we have the following:
\begin{equation*}
\left|I_1 \right| \leq \epsilon |\U_{\delta_1}|_{\frac{2N}{N+2}}|\bar\phi_1|_{\frac{2N}{N-2}}
\end{equation*}
Since $N\geq7$ we have $|\U_{\delta_i}|_{\frac{2N}{N+2}}=O(\delta_i^2)$, for $i=1,2$, so from our choice of $\delta_i$ (see (\ref{deltaj})) and since $\|\bar\phi_1\|\leq c \epsilon^{\frac{\theta_1}{2}+ \sigma}$  we deduce that
 \begin{equation}\label{eq3exp1}
\left|I_1 \right| \leq c \epsilon ( \epsilon^{\frac{2}{N-4}}  \epsilon^{\frac{N-2}{2(N-4)}+\sigma})\leq c \e^{\theta_1+\sigma},
\end{equation}
for all sufficiently small $\epsilon$.
For $I_2$, with similar computations,  we get that
\begin{equation*}
\left|I_2 \right| \leq \epsilon |\U_{\delta_2}|_{\frac{2N}{N+2}}|\bar\phi_1|_{\frac{2N}{N-2}} \leq c \epsilon^{1+\frac{2(3N-10)}{(N-4)(N-6)}} \epsilon^{\frac{N-2}{2(N-4)}+\sigma}.
\end{equation*}
Since $N\geq7$ it is elementary to see that $1+\frac{2(3N-10)}{(N-4)(N-6)}+\frac{N-2}{2(N-4)}> \theta_2 $. From this we deduce that
\begin{equation*}
\left|I_2 \right| \leq c \epsilon^{\theta_2+\sigma},
\end{equation*}
 for all sufficiently small $\epsilon$.

\textbf{D:} we have
\begin{equation}\label{eq4exp1}
\begin{array}{lll}
&& \hskip-1.5cm   \displaystyle\int_\Omega [f(\U_{\delta_1}) - f(\U_{\delta_2})-f(V_\e)] \bar\phi_1 \ dx  = \underbrace{\int_\Omega [f(\p\U_{\delta_1}) - f(\p\U_{\delta_2})-f(V_\e)] \bar\phi_1 \ dx}_{I_1} +\\[10pt]
&&+\underbrace{\int_\Omega [f(\U_{\delta_1})-f(\p\U_{\delta_1})]\bar\phi_1\, dx}_{I_2} +\underbrace{\int_\Omega[f(\p\U_{\delta_2})-f(\U_{\delta_2})]\bar\phi_1\, dx}_{I_3}\\[10pt]
\end{array}
\end{equation}
We evaluate separately the three terms.\\ We divide $\Omega$ into the three regions $A_0, A_1, A_2$ (see the proof of Proposition \ref{errorprop2} for their definition). Then
\begin{eqnarray*}
&&\hskip-1.0cm \int_\Omega [f(\p\U_{\delta_1}) - f(\p\U_{\delta_2})-f(V_\e)] \bar\phi_1 \ dx = \underbrace{\sum_{j=0}^1 \int_{A_j}[f(\p\U_{\delta_1})-f(V_\e)]\bar\phi_1\, dx}_{I_1'}\\
&&-\underbrace{\sum_{j=0}^1\int_{A_j}f(\p\U_{\delta_2})\bar\phi_1, dx}_{I_1''}+\underbrace{\int_{A_2}[f(\p\U_{\delta_1})-f(\p\U_{\delta_2})-f(V_\e)]\bar\phi_1\, dx}_{I_1'''}
\end{eqnarray*}
Now, writing $f(\p\U_{\delta_1})-f(V_\e)=f(\p\U_{\delta_1})-f(V_\e)+f^\prime(\p\U_{\delta_1})\p\U_{\delta_2} - f^\prime(\p\U_{\delta_1})\p\U_{\delta_2}$, applying the usual elementary inequalities, H\"older inequality and taking into account the computations made in \eqref{elemcompututil} , \eqref{eq2I1A1}, \eqref{eq3I1A1}, we get that

\begin{eqnarray*}
|I_1'|&\leq & c|\p\U_{\delta_2}^p|_{\frac{2N}{N+2}, A_0}|\bar\phi_1|_{\frac{2N}{N-2}, A_0} +c | \p\U_{\delta_2}^p|_{\frac{2N}{N+2}, A_1}|\bar\phi_1|_{\frac{2N}{N-2}, A_1}\\
&&+c|\p\U_{\delta_1}|_{\frac{2N}{N-2}, A_0}^{p-1}|\p\U_{\delta_2}|_{\frac{2N}{N+2}, A_0}|\bar\phi_1|_{\frac{2N}{N-2}, A_0}+c|\p\U_{\delta_1}^{p-1}\p\U_{\delta_2}|_{\frac{2N}{N+2}, A_1}|\bar\phi_1|_{\frac{2N}{N-2}, A_1}\\
&\leq & c_1\left(\delta_2^{\frac{N+2}{2}}\e^{\frac{\theta_1}{2}+\sigma}
+\left(\frac{\delta_2}{\delta_1}\right)^{\frac{N+2}{4}}\e^{\frac{\theta_1}{2}+\sigma}+\delta_1^2\delta_2^{\frac{N-2}{2}}\e^{\frac{\theta_1}{2}+\sigma}+\left(\frac{\delta_2}{\delta_1}\right)^2\left(\int_{\sqrt{\frac{\delta_1}{\delta_2}}}^{\frac{\rho}{\delta_2}}r^{\frac{-N^2+5N-2}{N+2}}\, dr\right)^{\frac{N+2}{2N}}\e^{\frac{\theta_1}{2}+\sigma}\right)\\
&\leq & c_2\left( \delta_2^{\frac{N+2}{2}}\e^{\frac{\theta_1}{2}+\sigma}
+\left(\frac{\delta_2}{\delta_1}\right)^{\frac{N+2}{4}}\e^{\frac{\theta_1}{2}+\sigma}+\delta_1^2\delta_2^{\frac{N-2}{2}}\e^{\frac{\theta_1}{2}+\sigma}+\left(\frac{\delta_2}{\delta_1}\right)^2\left(\frac{\delta_2}{\delta_1}\right)^{\frac{N-6}{4}}\e^{\frac{\theta_1}{2}+\sigma}\right)\\
&\leq & c_3 \left(\frac{\delta_2}{\delta_1}\right)^{\frac{N+2}{4}}\e^{\frac{\theta_1}{2}+\sigma} \leq c_4\e^{\theta_2+\sigma}.
\end{eqnarray*}
As before we have
$$|I_1''|\leq \sum_{j=0}^1|f(\p\U_{\delta_2})|_{\frac{2N}{N+2}, A_j}|\bar\phi_1|_{\frac{2N}{N-2}, A_j}\leq c\e^{\theta_2+\sigma}.$$

 Now, by H\"older inequality and reasoning as in (\ref{eq9err2}), \eqref{eq11err2}, (\ref{eqerrore2post}), we get that
 \begin{eqnarray*}
 |I_1'''|&\leq & |f(\p\U_{\delta_1}) - f(\p\U_{\delta_2})-f(V_\e)|_{\frac{2N}{N+2}, A_2}|\bar\phi_1|_{\frac{2N}{N-2}, A_2}\\
  &\leq &c_1 \left(|\p\U_{\delta_1}^p|_{\frac{2N}{N+2}, A_2}+ |\p\U_{\delta_2}^{p-1}\p\U_{\delta_1}|_{\frac{2N}{N+2}, A_2}\right)|\bar\phi_1|_{\frac{2N}{N-2}, A_2}\\
  &\leq&c_2\left(\frac{\delta_2}{\delta_1}\right)^{\frac{N+2}{4}}\e^{\frac{\theta_1}{2}+\sigma} \leq c_3\e^{\theta_2+\sigma}.
  \end{eqnarray*}
  At the end we conclude that
  $$|I_1| \leq c \e^{\theta_2+\sigma}.$$

For the remaining two terms of (\ref{eq4exp1}), reasoning as in the proof of Proposition \ref{errorprop1}, we get that

$$|f(\p\U_{\delta_i})-f(\U_{\delta_i}) |_{\frac{2N}{N+2}} \leq c  \delta_i^{\frac{N+2}{2}}.$$

Hence
$$|I_2|\leq |f(\U_{\delta_1}) - f(\p\U_{\delta_1})|_{\frac{2N}{N+2}}|\bar\phi_1|_{\frac{2N}{N-2}} \leq c \epsilon^{\frac{N+2}{2(N-4)}} \epsilon^{\frac{\theta_1}{2} + \sigma}\leq c \e^{\theta_1+\sigma},$$
for all sufficiently small $\epsilon$. We remark that $I_2$ depends only on $d_1$ and hence it is sufficient that it is of order $\theta_1+\sigma$.\\
At the end
$$|I_3|\leq |f(\U_{\delta_2}) - f(\p\U_{\delta_2})|_{\frac{2N}{N+2}}|\bar\phi_1|_{\frac{2N}{N-2}} \leq c \e^{\theta_2+\sigma},$$
for all sufficiently small $\epsilon$.

\textbf{E:} We decompose $\Omega$ in the three regions $A_j$, $j=0, 1, 2$ used before.\\ For $j=0,1$ we have
\begin{equation}\label{stimadelpezzoE}
\begin{array}{lll}
&&\hskip-1.0cm \displaystyle\int_{A_j} \left[|V_\e+\bar\phi_1|^{p+1}-|V_\e|^{p+1}-(p+1)|V_\e|^{p-1}V_\e \bar\phi_1\right]\, dx \\ &&\displaystyle =\underbrace{\int_{A_j}\left[|\p\U_{\delta_1}-\p\U_{\delta_2}+\bar\phi_1|^{p+1}-|\p\U_{\delta_1}+\bar\phi_1|^{p+1}
+(p+1)|\p\U_{\delta_1}+\bar\phi_1|^{p-1}(\p\U_{\delta_1}+\bar\phi_1)\p\U_{\delta_2}\, \right]dx}_{I_1}\\
&&\displaystyle-\underbrace{\int_{A_j}\left[|\p\U_{\delta_1}-\p\U_{\delta_2}|^{p+1}-\p\U_{\delta_1}^{p+1}+(p+1)\p\U_{\delta_1}^p \p\U_{\delta_2}\right]\, dx}_{I_2}\\
&&\displaystyle+\underbrace{\int_{A_j}\left[|\p\U_{\delta_1}+\bar\phi_1|^{p+1}-\p\U_{\delta_1}^{p+1}-(p+1)\p\U_{\delta_1}^p\bar\phi_1\right]\, dx}_{I_3}\\
&&\displaystyle-(p+1)\underbrace{\int_{A_j}\left[|\p\U_{\delta_1}-\p\U_{\delta_2}|^{p-1}(\p\U_{\delta_1}-\p\U_{\delta_2})-\p\U_{\delta_1}^p\right]\bar\phi_1\, dx}_{I_4}\\
&&\displaystyle-(p+1)\underbrace{\int_{A_j}\left[|\p\U_{\delta_1}+\bar\phi_1|^{p-1}(\p\U_{\delta_1}+\bar\phi_1)-\p\U_{\delta_1}^p\right]\p\U_{\delta_2}\, dx}_{I_5}.
\end{array}
\end{equation}
\\
In order to estimate $I_1$, $I_2$, $I_4$ and $I_5$, applying the usual elementary inequalities, we see that
 $$|I_1|\leq c\left(\int_{A_j}\p\U_{\delta_2}^{p+1}\, dx + \int_{A_j}\p\U_{\delta_1}^{p-1}\p\U_{\delta_2}^2\, dx+\int_{A_j}|\bar\phi_1|^{p-1}\p\U_{\delta_2}^{2}\, dx\right)$$ $$|I_2|\leq c\left( \int_{A_j}\p\U_{\delta_2}^{p+1}\, dx + \int_{A_j}\p\U_{\delta_1}^{p-1}\p\U_{\delta_2}^2\, dx\right)$$
$$|I_4|\leq  c\left( \int_{A_j}\p\U_{\delta_2}^{p}|\bar\phi_1|\, dx+ \int_{A_j}\p\U_{\delta_1}^{p-1}\p\U_{\delta_2}|\bar\phi_1|\, dx\right)$$ $$|I_5|\leq  c\left( \int_{A_j}|\bar\phi_1|^p\p\U_{\delta_2}\, dx+ \int_{A_j}\p\U_{\delta_1}^{p-1}\p\U_{\delta_2}|\bar\phi_1|\, dx\right).$$
Now, as seen in the proof of \eqref{elemcompututil} and thanks to \eqref{eq3I1A1}, we have
$$\int_{A_j}\p\U_{\delta_2}^{p+1}\, dx \leq c \left\{\begin{array}{lr} \delta_2^N \qquad \qquad\qquad\qquad\mbox{if}\,\, j=0\\ \left(\frac{\delta_2}{\delta_1}\right)^{\frac{N}{2}}\qquad\qquad \qquad  \mbox{if}\,\,\ j=1,\end{array}\right.$$
 \begin{equation*}
\int_{A_0}\p\U_{\delta_2}^2\p\U_{\delta_1}^{p-1}\, dx \leq c |\p\U_{\delta_1}|_{p+1, A_0}^{p-1}|\p\U_{\delta_2}|^2_{p+1, A_0}\leq c \delta_1^2\delta_2^{N-2}.
\end{equation*}
Moreover, by analogous computations, we get that 
$$
\int_{A_1}\p\U_{\delta_2}^2\p\U_{\delta_1}^{p-1}\, dx \leq  c_1 \left(\frac{\delta_2}{\delta_1}\right)^2\int_{\sqrt{\frac{\delta_1}{\delta_2}}}^{\frac{\rho}{\delta_2}}\frac{1}{r^{N-3}}\, dr\leq c_2\left(\frac{\delta_2}{\delta_1}\right)^{\frac{N}{2}},$$

$$\int_{A_j}\p\U_{\delta_2}^{2}|\bar\phi_1|^{p-1}\, dx \leq c |\p\U_{\delta_2}|_{p+1, A_j}^2\|\bar\phi_1\|^{p-1}\leq c\left\{\begin{array}{lr} \delta_2^{N-2}\e^{(p-1)(\frac{\theta_1}{2}+\sigma)} \qquad \qquad\qquad\qquad\mbox{if}\,\, j=0\\ \left(\frac{\delta_2}{\delta_1}\right)^{\frac{N-2}{2}}\e^{(p-1)(\frac{\theta_1}{2}+\sigma)}\qquad \qquad \qquad  \mbox{if}\,\,\ j=1,\end{array}\right.$$

$$\int_{A_j}\p\U_{\delta_2}^{p}|\bar\phi_1|\, dx \leq c |\p\U_{\delta_2}|_{p+1, A_j}^p\|\bar\phi_1\|\leq c\left\{\begin{array}{lr} \delta_2^{\frac{N+2}{2}}\e^{\frac{\theta_1}{2}+\sigma} \qquad \qquad\qquad\qquad\mbox{if}\,\, j=0\\ \left(\frac{\delta_2}{\delta_1}\right)^{\frac{N+2}{4}}\e^{\frac{\theta_1}{2}+\sigma}\qquad \qquad \qquad \mbox{if}\,\,\ j=1,\end{array}\right.$$

$$\int_{A_0}\p\U_{\delta_2}\p\U_{\delta_1}^{p-1}|\bar\phi_1|\, dx \leq c |\p\U_{\delta_2}|_{p+1, A_0}|\p\U_{\delta_1}|_{p+1, A_0}^{p-1}\|\bar\phi_1\|\leq c \delta_1^2\delta_2^{\frac{N-2}{2}}\e^{\frac{\theta_1}{2}+\sigma}, $$
and, thanks to \eqref{eq2I1A1}, we have 
\begin{eqnarray*}
\int_{A_1}\p\U_{\delta_2}\p\U_{\delta_1}^{p-1}|\bar\phi_1|\, dx&\leq & |\p\U_{\delta_2}\p\U_{\delta_1}^{p-1}|_{\frac{2N}{N+2}, A_1}|\bar\phi_1|_{\frac{2N}{N-2}, A_1}\\
&\leq &c_1 \left[\int_{A_1}\left(\frac{\delta_2^{\frac{N-2}{2}}\delta_1^2}{(\delta_2^2+|x|^2)^{\frac{N-2}{2}}
(\delta_1^2+|x|^2)^2}\right)^{\frac{2N}{N+2}}\, dx\right]^{\frac{N+2}{2N}}\e^{\frac{\theta_1}{2}+\sigma}\\
&\leq & c_2 \left(\frac{\delta_2}{\delta_1}\right)^2\left(\frac{\delta_2}{\delta_1}\right)^{\frac{N-6}{4}} \e^{\frac{\theta_1}{2}+\sigma}=c_2 \left(\frac{\delta_2}{\delta_1}\right)^{\frac{N+2}{4}} \e^{\frac{\theta_1}{2}+\sigma}
\end{eqnarray*}
At the end
$$\int_{A_0}\p\U_{\delta_2}|\bar\phi_1|^{p}\, dx \leq c_1 |\p\U_{\delta_2}|_{p+1, A_0}\|\bar\phi_1\|^{p}\leq c_2 \delta_2^{\frac{N-2}{2}}\e^{p(\frac{\theta_1}{2}+\sigma)}$$
and, by using Lemma \ref{techlembehvphi1}, we get that
\begin{eqnarray*}
\int_{A_1}\p\U_{\delta_2}|\bar\phi_1|^{p}\, dx&\leq & c_1 |\bar\phi_1|_\infty^{p-1}\left[\int_{A_1}\p\U_{\delta_2}^{\frac{2N}{N+2}}\, dx\right]^{\frac{N+2}{2N}}|\bar\phi_1|_{p+1, A_1}\\
&\leq & c_2 \e^{-\frac{2}{N-4}}\delta_2^2\left[\int_{\sqrt{\frac{\delta_1}{\delta_2}}}^{\frac{\rho}{\delta_2}}r^{\frac{-N^2+5N-2}{N+2}}\, dr\right]^{\frac{N+2}{2N}}\e^{\frac{\theta_1}{2}+\sigma}\\
&= & c_3\e^{-\frac{2}{N-4}}\delta_2^2\left[\left(\frac{\delta_2}{\delta_1}\right)^{\frac{N-6}{4}}-\delta_2^{\frac{N-6}{2}}\right] \e^{\frac{\theta_1}{2}+\sigma}\\
&\leq & c_4\left(\frac{\delta_2}{\delta_1}\right)^{\frac{N+2}{4}} \e^{\frac{\theta_1}{2}+\sigma}.
\end{eqnarray*}

In order to estimate $I_3$ we observe that
\begin{equation}\label{enumeratastarstar}
\left|\int_{A_0}\left[|\p\U_{\delta_1}+\bar\phi_1|^{p+1}-\p\U_{\delta_1}^{p+1}-(p+1)\p\U_{\delta_1}^p\bar\phi_1\right]\, dx\right|
\leq c_1\left( \|\bar\phi_1\|^2|\p\U_{\delta_1}|_{p+1}^{p-1}+\|\bar\phi_1\|^{p+1}\right)\leq c_2 \e^{\theta_1+\sigma}, 
\end{equation}
which is sufficient since this term does not depend on $d_2$.\\
Moreover
\begin{eqnarray*}
&&\int_{A_1}\left[|\p\U_{\delta_1}+\bar\phi_1|^{p+1}-\p\U_{\delta_1}^{p+1}-(p+1)\p\U_{\delta_1}^p\bar\phi_1\right]\, dx = \int_{B(0, \rho)}\left[|\p\U_{\delta_1}+\bar\phi_1|^{p+1}-\p\U_{\delta_1}^{p+1}-(p+1)\p\U_{\delta_1}^p\bar\phi_1\right]\, dx \\ &&-\int_{A_2}\left[|\p\U_{\delta_1}+\bar\phi_1|^{p+1}-\p\U_{\delta_1}^{p+1}-(p+1)\p\U_{\delta_1}^p\bar\phi_1\right]\, dx.
\end{eqnarray*}
We observe that the first integral in the right-hand side of the previous equation depends only on $d_1$. Hence, as in \eqref{enumeratastarstar}, we have
$$\left|\int_{B(0, \rho)}\left[|\p\U_{\delta_1}+\bar\phi_1|^{p+1}-\p\U_{\delta_1}^{p+1}-(p+1)\p\U_{\delta_1}^p\bar\phi_1\right]\, dx\right|\leq c \e^{\theta_1+\sigma}.$$
Furthermore, by using Lemma \ref{techlembehvphi1}, we get that

\begin{equation}\label{prevestimatelast}
\begin{array}{lll}
&&\displaystyle \left|\int_{A_2}\left[|\p\U_{\delta_1}+\bar\phi_1|^{p+1}-\p\U_{\delta_1}^{p+1}-(p+1)\p\U_{\delta_1}^p\bar\phi_1\right]\, dx\right|\leq c_1\left(|\p\U_{\delta_1}|^{p-1}_{p+1, A_2}|\bar\phi_1|_{p+1, A_2}^2+ |\bar\phi_1|_{p+1, A_2}^{p+1}\right)\\[10pt]
&&\displaystyle\leq c_2\left( \left[\int_{B(0, \sqrt{\frac{\d_2}{\d_1}})}\frac{1}{(1+|y|^2)^N}\, dy\right]^{\frac{2}{N}}|\bar\phi_1|_\infty^2\left[\int_{A_2}1\, dx\right]^{\frac{2}{p+1}}+ |\bar\phi_1|_\infty^{p+1}\int_{A_2}1\, dx\right)\\[20pt]
&&\displaystyle\leq c_3\left( \left[\int_0^{\sqrt{\frac{\d_2}{\d_1}}}r^{N-1}\, dr\right]^{\frac{2}{N}}\e^{-\frac{N-2}{N-4}}\left[\int_0^{\sqrt{\d_1\d_2}} r^{N-1}\, dr\right]^{\frac{2}{p+1}}+  \e^{-\frac{N}{N-4}}\int_0^{\sqrt{\d_1\d_2}}r^{N-1}\, dr\right)\\[20pt]
&&\displaystyle\leq c_4\left( \frac{\d_2}{\d_1}\e^{-\frac{N-2}{N-4}}(\d_1\d_2)^{\frac{N-2}{2}}+ \e^{-\frac{N}{N-4}}(\d_1\d_2)^{\frac{N}{2}}\right)\\[12pt]
&&\displaystyle\leq c_5 \e^{{\theta_2}+\sigma}.
\end{array}
\end{equation}

Now, it remains only to estimate the left-hand side of \eqref{stimadelpezzoE} for $j=2$. Hence, thanks to the usual elementary inequalities, we get that
\begin{eqnarray*}
&&\hskip-1.0cm \left|\int_{A_2} \left[|V_\e+\bar\phi_1|^{p+1}-|V_\e|^{p+1}-(p+1)|V_\e|^{p-1}V_\e \bar\phi_1\right]\, dx\right| \\
&&\leq c \left(\int_{A_2}|V_\e|^{p-1} \bar\phi_1^2 \ dx+  \int_{A_2} |\bar\phi_1|^{p+1} \ dx\right)\\
&&\leq c \left(\int_{A_2}\p\U_{\delta_1}^{p-1} \bar\phi_1^2 \ dx+\int_{A_2}\p\U_{\delta_2}^{p-1} \bar\phi_1^2 \ dx+  \int_{A_2} |\bar\phi_1|^{p+1} \ dx\right)
\end{eqnarray*}
For the first and third integrals in the last right-hand side we can reason as in \eqref{prevestimatelast}. For the second integral, using Lemma \ref{techlembehvphi1}, we have
\begin{eqnarray*}
\int_{A_2}\p\U_{\d_2}^{p-1}\bar\phi_1^2\, dx&\leq & c_1|\bar\phi_1|_\infty^2\int_{A_2}\frac{\d_2^2}{(\d_2^2+|x|^2)^2}\, dx \\
&\leq & c_2 \d_1^{-(N-2)}\d_2^2\int_{A_2}\frac{1}{|x|^4}\, dx\\
&\leq & c_3 \d_1^{-N+2}\d_2^2\int_0^{\sqrt{\d_1\d_2}} r^{N-5}\, dr\\
&\leq & c_4 \left(\frac{\d_2}{\d_1}\right)^{\frac{N}{2}}
\end{eqnarray*}


Finally, summing up all the estimates, we conclude that $|{\bf E}| =  \e^{\theta_1+\sigma}g(d_1)+O(\e^{\theta_2+\sigma})$.\\\\
\end{proof}

\begin{lemma}\label{lem2exp1}
For any $\eta>0$ there exists $\epsilon_0>0$ such that for any $\epsilon \in (0,\epsilon_0)$ it holds:
$$J_\epsilon(V_\e+\bar\phi_1+\bar\phi_2)=J_\epsilon(V_\e+\bar\phi_1) +O (\epsilon^{\theta_2+ \sigma}), $$
$C^0$-uniformly with respect to $(d_1,d_2)$ satisfying condition \eqref{limdj}, for some positive real number $\sigma$ depending only on $N$.
\end{lemma}
\begin{proof}
As we have seen in the proof of Lemma \ref{lem1exp1}, by direct computation we get that
\begin{equation}\label{eq1exp2}
\begin{array}{lll}
&&\hskip-1.0cm   J_\epsilon(V_\e+\bar\phi_1+\bar\phi_2) - J_\epsilon(V_\e+\bar\phi_1)
=\frac{1}{2} \int_\Omega |\nabla \bar\phi_2|^2 \ dx + \int_\Omega \nabla (V_\e+\bar\phi_1) \cdot \nabla \bar\phi_2 \ dx \\[10pt]
&&-\frac{\epsilon}{2} \int_\Omega |\bar\phi_2|^2 \ dx - \epsilon \int_\Omega (V_\e +\bar\phi_1) \bar\phi_2 \ dx- \frac{1}{p+1} \int_\Omega (|V_\e+\bar\phi_1+\bar\phi_2|^{p+1} -|V_\e+\bar\phi_1|^{p+1} ) \ dx\\[10pt]
&=&- \frac{1}{2} \|\bar\phi_2\|^2   +\frac{\epsilon}{2} |\bar\phi_2|_2^2 + \int_\Omega \nabla (V_\e + \bar\phi_1 + \bar\phi_2) \cdot \nabla \bar\phi_2 \ dx  \\[10pt]
&& - \epsilon \int_\Omega (V_\e+ \bar\phi_1+\bar\phi_2) \bar\phi_2 \ dx - \int_\Omega f(V_\e+\bar\phi_1)\bar\phi_2 \ dx \\[10pt]
&&-  \int_\Omega [F(V_\e+\bar\phi_1+\bar\phi_2) -F(V_\e+\phi_1)-f(V_\e+\bar\phi_1)\bar\phi_2] \ dx\\[10pt]
\end{array}
\end{equation}
Since $\bar\phi_1+\bar\phi_2$ is a solution of (\ref{aux}) we have
$$\Pi^{\perp}\{V_\e+\bar\phi_1+\bar\phi_2 - i^*[\epsilon (V_\e+\bar\phi_1+\bar\phi_2) + f(V_\e+\bar\phi_1+\bar\phi_2)]\}=0, $$
hence, for some $\psi \in \mathcal K$, we get that $V_\e+\bar\phi_1+\bar\phi_2$ weakly solves
\begin{equation}\label{eq2exp2}
- \Delta (V_\e+ \bar\phi_1+ \bar\phi_2) + \Delta \bar\psi -  [\epsilon  (V_\e+\bar\phi_1+\bar\phi_2)+ f(V_\e+\bar\phi_1+\bar\phi_2)]=0.
\end{equation}
Choosing $\bar\phi_2$ as test function, since $\bar\phi_2 \in \mathcal{K}^\perp$,  $\psi \in \mathcal K$  we deduce that
\begin{equation}\label{eq3exp2}
\int_\Omega \nabla (V_\e+\bar\phi_1+\bar\phi_2) \cdot \nabla \bar\phi_2 \ dx - \epsilon \int_\Omega (V_\e+\bar\phi_1+\bar\phi_2) \bar\phi_2 \ dx = \displaystyle \int_\Omega  f(V_\e+\bar\phi_1+\bar\phi_2)  \bar\phi_2 \ dx
\end{equation}
Thanks to (\ref{eq3exp2}) we rewrite (\ref{eq1exp2}) as
\begin{eqnarray}\label{eq4exp2}
\nonumber
J_\epsilon(V_\e+\bar\phi_1+\bar\phi_2) - J_\epsilon(V_\e+\bar\phi_1)
&=&- \frac{1}{2} \|\bar\phi_2\|^2   +\frac{\epsilon}{2} |\bar\phi_2|_2^2  + \int_\Omega [f(V_\e+\bar\phi_1+\bar\phi_2)-f(V_\e+\bar\phi_1)]\bar\phi_2 \ dx \\
\nonumber
&&-  \int_\Omega [F(V_\e+\bar\phi_1+\bar\phi_2) -F(V_\e+\phi_1)-f(V_\e+\bar\phi_1)\bar\phi_2] \ dx\\
&=&A+B+C+D.
\end{eqnarray}

\textbf{A, B:}  Thanks to Proposition \ref{auxsolving}, for all sufficiently small $\epsilon$, we have $\|\bar\phi_2\|\leq c \epsilon^{\frac{\theta_2}{2}+ \sigma}$, for some $c>0$ and for some $\sigma>0$ depending only on $N$.  Hence we deduce that $A=O(\epsilon^{\theta_2+ 2 \sigma})$, $B=O(\epsilon^{\theta_2+ 2 \sigma+1})$.

 \textbf{C:} By Lemma \ref{lem0n2} we get

\begin{eqnarray*}
\left|\int_\Omega [f(V_\e+\bar\phi_1+\bar\phi_2)-f(V_\e+\bar\phi_1)]\bar\phi_2 \ dx\right| &\leq & \int_{\O}| \bar\phi_2|^{p+1}\, dx +\int_{\O} |V_\e+\bar\phi_1|^{p-1}\bar\phi_2^2\, dx\\
&\leq & c \|\bar\phi_2\|^{p+1}+c |V_\e|^{p-1}_{p+1}|\bar\phi_2|_{p+1}^2+c|\bar\phi_1|^{p-1}_{p+1}|\bar\phi_2|_{p+1}^2\\&\leq& c\e^{\theta_2+\sigma}
\end{eqnarray*}
for all sufficiently small $\epsilon$.

\textbf{D:} Applying Lemma \ref{lem0n2} and H\"older inequality we get that
\begin{equation*}
\left|\int_\Omega [F(V_\e+\bar\phi_1+\bar\phi_2) -F(V_\e+\bar\phi_1)-f(V_\e+\bar\phi_1)\bar\phi_2] \ dx \right|
\leq c |V_\e|_{p+1}^{p-1}|\bar\phi_2|_{p+1}^2 + c |\bar\phi_1|_{p+1}^{p-1} |\bar\phi_2|_{p+1}^2  + c |\bar\phi_2|_{p+1}^{p+1}.
\end{equation*}

Since all the terms from $A$ to $D$ are high order terms with respect to $\epsilon^{{\theta_2}}$ the proof is complete.

\end{proof}
In order to prove Proposition \ref{ridotto} some further preliminary lemmas are needed.
\begin{lemma} \label{lemma1}
Let $\delta_j$ as in \eqref{deltaj} for $j=1, 2$ and $N\geq 7$. For any $\eta>0$ there exists $\epsilon_0>0$ such that for any $\epsilon \in (0,\epsilon_0)$,  it holds
$$\frac 12 \int_{\O}|\nabla\p\U_{\delta_j}|^2\, dx-\frac{1}{p+1}\int_{\O}\p\U_{\delta_j}^{p+1}\, dx = \frac{1}{N} S^{N/2} + a_1 \tau(0) \delta_j^{N-2} + O(\delta_j^{N-1}),$$
$C^0$-uniformly with respect to $(d_1,d_2)$ satisfying condition \eqref{limdj},
where  $a_1:=\frac{1}{2} \alpha_{N}^{p+1} \int_{\R} \frac{1}{(1+|y|^2)^{\frac{N+2}{2}}} \ dy$ and $\tau(0)$ is the Robin's function of the domain $\Omega$ at the origin.
\end{lemma}
\begin{proof}
By using \eqref{proiezione}, \eqref{pbproiezione} and \eqref{expvarphi} we have that
\begin{eqnarray*}
&&\hskip-1.0cm \frac 12 \int_{\O}|\nabla \p\U_{\delta_j}|^2\, dx-\frac{1}{p+1}\int_{\O}\p\U_{\delta_j}^{p+1}\, dx= \frac{1}{2}  \int_\Omega \U_{\delta_j}^{p} \p\U_{\delta_j} \ dx - \frac{1}{p+1}\int_\Omega \p\U_{\delta_j}^{p+1} \ dx\\
&=& \frac{1}{2}  \int_\Omega \U_{\delta_j}^{p} (\U_{\delta_j}-\varphi_{\delta_j}) \ dx - \frac{1}{p+1}\int_\Omega (\U_{\delta_j}-\varphi_{\delta_j})^{p+1} \ dx\\
&=& \frac{1}{2}  \int_\Omega \U_{\delta_j}^{p+1} \ dx - \frac{1}{2}  \int_\Omega \U_{\delta_j}^p \varphi_{\delta_j} \ dx - \frac{1}{p+1}\int_\Omega \U_{\delta_j}^{p+1}  \ dx + \int_\Omega \U_{\delta_j}^p \varphi_{\delta_j} \ dx + O\left( \int_\Omega \U_{\delta_j}^{p-1} \varphi_{\delta_j}^2 \ dx\right)\\
&=& \left(\frac{1}{2}- \frac{1}{p+1}\right)  \int_\Omega \U_{\delta_j}^{p+1} \ dx + \frac{1}{2}  \int_\Omega \U_{\delta_j}^p \varphi_{\delta_j} \ dx  + O\left( \int_\Omega \U_{\delta_j}^{p-1} \varphi_{\delta_j}^2 \ dx\right)\\
\end{eqnarray*}
Now it's easy to see that
\begin{eqnarray}\label{primalem1}
\int_\Omega \U_{\delta_j}^{p+1} \ dx=  \int_{\R} \frac{\alpha_N^{p+1}}{(1+|y|^2)^{N}} \ dy + O(\delta_j^N),
\end{eqnarray}
while
\begin{eqnarray}\label{secondalem1}
\nonumber
\int_\Omega \U_{\delta_j}^p \varphi_{\delta_j} \ dx&=& \int_\Omega \U_{\delta_j}^p \left(\alpha_N \delta_j^{\frac{N-2}{2}}H(0,x)+O(\delta_j^{\frac{N+2}{2}})\right)  \ dx\\
\nonumber
&=& \alpha_N \delta_j^{\frac{N-2}{2}} \int_\Omega \U_{\delta_j}^p H(0,x) \ dx + O\left(\delta_j^{\frac{N+2}{2}} \int_\Omega \U_{\delta_j}^p  \ dx\right)\\
&=& \alpha_{N}^{p+1} \tau(0) \delta_j^{N-2} \int_{\R} \frac{1}{(1+|y|^2)^{\frac{N+2}{2}}} \ dy + O(\delta_j^{N-1}).
\end{eqnarray}
Moreover
\begin{equation}\label{terzalem1}
O\left( \int_\Omega \U_{\delta_j}^{p-1} \varphi_{\delta_j}^2 \ dx\right)=O(\delta_j^{N-1}).
\end{equation}
Indeed, we get
\begin{eqnarray*}
\int_{\O}\U_{\delta_j}^{p-1}\varphi_{\delta_j}^2\, dx &=& \int_{B_{\sqrt{\delta_j}}(0)}\U_{\delta_j}^{p-1}\varphi_{\delta_j}^2\, dx +\int_{\O\setminus B_{\sqrt{\delta_j}}(0)}\U_{\delta_j}^{p-1}\varphi_{\delta_j}^2\, dx\\
&\leq & c_1 \delta_j^{N-2} \int_{B_{\sqrt{\delta_j}}(0)}\U_{\delta_j}^{p-1}\, dx+ |\varphi_{\delta_j}|_{p+1}^2\left(\int_{\O\setminus B_{\sqrt{\delta_j}}(0)}\U_{\delta_j}^{p+1}\, dx\right)^{\frac{p-1}{p+1}}\\
&\leq & c_2 \delta_j ^{2N-4}\int_0^{\frac{1}{\sqrt{\delta_j}}}\frac{r^{N-1}}{(1+r^2)^2}\, dr+c_3 \delta_j^{N-2}\left(\int_{\frac{1}{\sqrt{\delta_j}}}^{+\infty}\frac{r^{N-1}}{(1+r^2)^N}\ dr\right)^{\frac{p-1}{p+1}}\\
&\leq & c_4\delta_j ^{\frac{3N-4}{2}}+c_5\delta_j^{N-1}\leq c_6\delta_j^{N-1}.
\end{eqnarray*}
Hence, from (\ref{primalem1}), (\ref{secondalem1}), (\ref{terzalem1}) we get the thesis.
\end{proof}

\begin{lemma} \label{lemma6}
Let $\delta_j$ as in \eqref{deltaj} for $j=1,2$ and $N\geq 7$. For any $\eta>0$ there exists $\epsilon_0>0$ such that for any $\epsilon \in (0,\epsilon_0)$,  it holds
$$\frac{\e}{2}\int_{\O}\p\U_{\delta_j}^2\, dx=a_2 \e \delta_j^2+O(\e\delta_j^{\frac{N}{2}}),$$
$C^0$-uniformly with respect to $(d_1,d_2)$ satisfying condition \eqref{limdj}, where
$a_2:= \frac{1}{2}\alpha_N^2\int_{\R} \frac{1}{(1+|y|^2)^{N-2}} \ dy$.
\end{lemma}
\begin{proof}
From (\ref{proiezione}) we get that
\begin{equation}\label{primalem6}
\frac{\e}{2} \int_\Omega (\p\U_{\delta_j})^2 \ dx= \displaystyle \frac{\e}{2} \int_\Omega (\U_{\delta_j}-\varphi_{\delta_j})^2 \ dx=\frac{\e}{2} \int_\Omega \U_{\delta_j}^2 \ dx -\epsilon \int_\Omega \U_{\delta_j}\varphi_{\delta_j} \ dx + \frac{\e}{2} \int_\Omega \varphi_{\delta_j}^2 \ dx.
\end{equation}
The principal term is the first one, in fact we have:
\begin{eqnarray}\label{secondalem6}
\begin{array}{lll}
\displaystyle \frac{\e}{2} \int_\Omega \U_{\delta_j}^2 \ dx &=&\displaystyle\frac{\e}{2} \ \alpha_N^2 \int_\Omega \frac{\delta_j^{N-2}}{(\delta_1^2+|x|^2)^{N-2}} \ dx=\displaystyle\frac{\e}{2} \ \alpha_N^2 \int_\Omega \frac{\delta_j^{-(N-2)}}{(1+|x/\delta_1|^2)^{N-2}} \ dx\\[12pt]
  &=&\displaystyle\frac{\e}{2} \ \alpha_N^2 \int_{\Omega/\delta_j} \frac{\delta_j^{-(N-2)}}{(1+|y|^2)^{N-2}} \delta_j^N \ dy=\displaystyle\frac{\e}{2} \ \alpha_N^2  \delta_j^2 \int_{\R} \frac{1}{(1+|y|^2)^{N-2}} \ dy +\\[16pt]
   &&+\displaystyle O\left( \epsilon\delta_j^2 \int_{1/\delta_j}^{+\infty} \frac{r^{N-1}}{(1+r^2)^{N-2}} \ dr \right) \\[12pt]
  &=&\displaystyle\frac{\e}{2} \ \alpha_N^2  \delta_j^2 \int_{\R} \frac{1}{(1+|y|^2)^{N-2}} \ dy + O\left( \epsilon\delta_j^{N-2} \right).
   \end{array}
\end{eqnarray}
For the remaining terms, by using also \eqref{stimavarphi}, we deduce that
\begin{equation}\label{ultimalem6}
\e\int_\Omega \U_{\delta_j}\varphi_{\delta_j} \ dx \leq \e |\U_{\delta_j}|_2|\varphi_{\delta_j}|_2\leq c\e\delta_j\delta_j^{\frac{N-2}{2}}\leq c\e\delta_j^{\frac{N}{2}}.
\end{equation}
Moreover by using again \eqref{stimavarphi} $$\frac{\e}{2}\int_\Omega \varphi_{\delta_j}^2\, dx=\frac{\e}{2}|\varphi_{\delta_j}|_2^2\leq C \e\delta_j^{N-2}$$ and the lemma in proved.
\end{proof}

\begin{lemma} \label{lemma7}
Let $\delta_j$ as in \eqref{deltaj} for $j=1,2$ and $N\geq 7$. For any $\eta>0$ there exists $\e_0>0$ such that for any $\e\in (0, \e_0)$ it holds
$$\e\int_\Omega \p\U_{\delta_1}\p\U_{\delta_2}\, dx =O\left(\e\left(\frac{\delta_2}{\delta_1}\right)^{\frac{N-2}{2}}\delta_1^2\right),$$
$C^0$-uniformly with respect to $(d_1,d_2)$ satisfying condition \eqref{limdj}.
\end{lemma}
\begin{proof}
From (\ref{proiezione}) we get that
\begin{eqnarray}\label{primalem7}
\begin{array}{lll}
\displaystyle \epsilon \int_\Omega \p\U_{\delta_1}\ \p\U_{\delta_2} \ dx &=& \displaystyle \epsilon \int_\Omega (\U_{\delta_1}-\varphi_{\delta_1})(\U_{\delta_2}-\varphi_{\delta_2}) \ dx\\[12pt]
&=&\displaystyle \epsilon \int_\Omega \U_{\delta_1}\U_{\delta_2} \ dx - \epsilon \int_\Omega \U_{\delta_1}\varphi_{\delta_2} \ dx - \epsilon \int_\Omega \U_{\delta_2}\varphi_{\delta_1} \ dx + \epsilon \int_\Omega \varphi_{\delta_1}\varphi_{\delta_2} \ dx.\\[12pt]
    \end{array}
\end{eqnarray}
We analyze every term.
\begin{eqnarray}\label{secondalem7}
\begin{array}{lll}
\displaystyle \epsilon \int_\Omega \U_{\delta_1}\U_{\delta_2} \ dx
&=&\displaystyle\epsilon \ \alpha_N^2 \int_\Omega \frac{\delta_1^{-\frac{N-2}{2}}}{(1+|x/\delta_1|^2)^{\frac{N-2}{2}}}  \frac{\delta_2^{\frac{N-2}{2}}}{(\delta_2^2+|x|^2)^{\frac{N-2}{2}}} \ dx\\[12pt]
&=&\displaystyle\epsilon \ \alpha_N^2 \int_{\Omega/\delta_1} \frac{\delta_1^{\frac{N+2}{2}}}{(1+|y|^2)^{\frac{N-2}{2}}}  \frac{\delta_2^{\frac{N-2}{2}}}{(\delta_2^2+\delta_1^2|y|^2)^{\frac{N-2}{2}}} \ dy\\[12pt]
&=&\displaystyle\epsilon \ \alpha_N^2 \int_{\Omega/\delta_1} \frac{\delta_1^{-\frac{N-6}{2}}}{(1+|y|^2)^{\frac{N-2}{2}}}  \frac{\delta_2^{\frac{N-2}{2}}}{\left(\left(\frac{\delta_2}{\delta_1}\right)^2+|y|^2\right)^{\frac{N-2}{2}}} \ dy\\[12pt]
&\leq&\displaystyle\epsilon \ \alpha_N^2 \left(\frac{\delta_2}{\delta_1}\right)^{\frac{N-2}{2}} \delta_1^2 \int_{\Omega/\delta_1} \frac{1}{(1+|y|^2)^{\frac{N-2}{2}}|y|^{N-2}} \ dy \\[12pt]
&=&\displaystyle\epsilon \ \alpha_N^2 \left(\frac{\delta_2}{\delta_1}\right)^{\frac{N-2}{2}} \delta_1^2 \int_{\R} \frac{1}{(1+|y|^2)^{\frac{N-2}{2}}|y|^{N-2}} \ dy\\[12pt]
&+& \displaystyle O\left(\epsilon \left(\frac{\delta_2}{\delta_1}\right)^{\frac{N-2}{2}} \delta_1^2 \int_{1/\delta_1}^{+\infty} \frac{r^{N-1}}{(1+r^2)^{\frac{N-2}{2}}r^{N-2}} \ dr\right)\\[12pt]
&=&\displaystyle\epsilon \ \alpha_N^2 \left(\frac{\delta_2}{\delta_1}\right)^{\frac{N-2}{2}} \delta_1^2 \int_{\R} \frac{1}{(1+|y|^2)^{\frac{N-2}{2}}|y|^{N-2}} \ dy+ \displaystyle O\left(\epsilon \left(\frac{\delta_2}{\delta_1}\right)^{\frac{N-2}{2}} \delta_1^{N-2}\right).
\end{array}
\end{eqnarray}
Hence $\epsilon \int_\Omega \U_{\delta_1}\U_{\delta_2} \ dx =O\left(\epsilon \left(\frac{\delta_2}{\delta_1}\right)^{\frac{N-2}{2}} \delta_1^2\right)$. Thanks to (\ref{ultimalem6}) we deduce that
$\epsilon \int_\Omega \U_{\delta_1}\varphi_{\delta_2} \ dx = O\left(\epsilon\ \delta_1^{\frac{N-2}{2}}  \delta_2^{\frac{N-2}{2}}\right)$, $ \epsilon \int_\Omega \U_{\delta_2}\varphi_{\delta_1} \ dx=O\left(\epsilon \ \delta_1^{\frac{N-2}{2}}  \delta_2^{\frac{N-2}{2}}\right)$. Moreover it's clear that $\epsilon \int_\Omega \varphi_{\delta_1}\varphi_{\delta_2} \ dx = O\left(\epsilon\ \delta_1^{\frac{N-2}{2}}  \delta_2^{\frac{N-2}{2}}\right)$. Since these last three terms are high order terms compared to $\epsilon \left(\frac{\delta_2}{\delta_1}\right)^{\frac{N-2}{2}} \delta_1^2$, we deduce the thesis, and the proof is complete.
\end{proof}

We are ready to prove Proposition \ref{ridotto}.
\\\\
\begin{proof}[Proof of Proposition \ref{ridotto}]
\begin{itemize}
\item[{(i):}] One can reason as Part 1 of Proposition 2.2 of \cite{Musso1}.\\
\item[{(ii):}] Let us fix $\eta>0$. From Lemma \ref{lem1exp1} and Lemma \ref{lem2exp1}, for all sufficiently small $\epsilon$, we get that
$$J_\e(V_\e+\bar\phi_1+\bar\phi_2)=J_\e(V_\e)+\e^{\theta_1+\sigma}g(d_1)+O(\e^{\theta_2+\sigma}),$$
for some $\sigma>0$.
We evaluate $J_\e(V_\e)=J_\e(\p\U_{\delta_1}-\p\U_{\delta_2})$.
\begin{eqnarray*}
J_\epsilon(\p\U_{\delta_1}-\p\U_{\delta_2})&=&\frac{1}{2} \int_\Omega |\nabla(\p\U_{\delta_1}-\p\U_{\delta_2})|^2 \ dx -\frac{1}{p+1} \int_\Omega |\p\U_{\delta_1}-\p\U_{\delta_2}|^{p+1} \ dx\\
&-&\frac{\epsilon}{2} \int_\Omega (\p\U_{\delta_1}-\p\U_{\delta_2})^2 dx \\
&=&\frac{1}{2} \int_\Omega |\nabla \p\U_{\delta_1}|^2 \ dx +\frac{1}{2} \int_\Omega |\nabla \p\U_{\delta_2}|^2 \ dx - \int_\Omega \nabla \p\U_{\delta_1} \cdot  \nabla \p\U_{\delta_2} \ dx \\
&-&\frac{1}{p+1} \int_\Omega |\p\U_{\delta_1}-\p\U_{\delta_2}|^{p+1} \ dx-\frac{\epsilon}{2} \int_\Omega (\p\U_{\delta_1})^2\ dx
-\frac{\epsilon}{2} \int_\Omega (\p\U_{\delta_2})^2\ dx \\
&+&\epsilon \int_\Omega \p\U_{\delta_1}\ \p\U_{\delta_2} \ dx\\
&=&\underbrace{\sum_{j=1}^2\left(\frac{1}{2}\int_\Omega |\nabla \p\U_{\delta_j}|^2\, dx-\frac{1}{p+1}\int_{\O}\p\U_{\delta_j}^{p+1}\, dx\right)}_{(I)}-\underbrace{\sum_{j=1}^2 \frac{\e}{2}\int_{\Omega} \p\U_{\delta_j}^2\, dx}_{(II)}\\
&&+\underbrace{\e\int_{\Omega} \p\U_{\delta_1}\p\U_{\delta_2}\, dx}_{(III)} \underbrace{-\int_\Omega \nabla\p\U_{\delta_1}\nabla\p\U_{\delta_2}\, dx}_{(IV)}\\
&&\underbrace{-\frac{1}{p+1}\int_{\O}\left[|\p\U_{\delta_1}-\p\U_{\delta_2}|^{p+1}-\p\U_{\delta_1}^{p+1}-\p\U_{\delta_2}^{p+1}\right]\, dx}_{(IV)}.
\end{eqnarray*}
By Lemma \ref{lemma1}, Lemma  \ref{lemma6} and Lemma  \ref{lemma7} we get
\begin{eqnarray*}
(I)&=&\frac{2}{N}S^{N/2}+a_1\tau(0)\delta_1^{N-2}+a_1\tau(0)\delta_2^{N-2}+O(\delta_1^{N-1})+O(\delta_2^{N-1}),\\
(II)&=& a_2\e\delta_1^2+a_2\e\delta_2^2+ O(\e\delta_1^{\frac{N}{2}})+O(\e\delta_2^{\frac{N}{2}}),\\
(III)&=&O\left(\e\left(\frac{\delta_2}{\delta_1}\right)^{\frac{N-2}{2}}\delta_1^2\right).
\end{eqnarray*}
Now since $-\Delta\p\U_{\delta_2}=\U_{\delta_2}^p$ then $\int_\Omega\nabla\p\U_{\delta_1}\nabla\p\U_{\delta_2}\, dx=\int_{\O}\U_{\delta_2}^p\p\U_{\delta_1}\, dx$ and hence
\begin{eqnarray*}
(IV)&=&\underbrace{-\frac{1}{p+1}\int_{\O}\left[|\p\U_{\delta_1}-\p\U_{\delta_2}|^{p+1}-\p\U_{\delta_1}^{p+1}-\p\U_{\delta_2}^{p+1}+(p+1)\p\U_{\delta_2}^p\p\U_{\delta_1}\right]\, dx}_{I_1}\\
&&+\underbrace{\int_{\O}\left[\p\U_{\delta_2}^p-\U_{\delta_2}^p\right]\p\U_{\delta_1}\, dx}_{I_2}.
\end{eqnarray*}
By (\ref{proiezione}) and Lemma \ref{lem0n2} we deduce that
\begin{equation*}
|I_2|\leq  C \int_{\O}\U_{\delta_2}^{p-1}\varphi_{\delta_2}\p\U_{\delta_1}\, dx + C\int_{\O}\varphi_{\delta_2}^p \p\U_{\delta_1}\, dx.
\end{equation*}
Now let $\rho>0$ such that $B(0,\rho)\subset \O$.
\begin{eqnarray*}
\int_\Omega \varphi_{\delta_2}^p\p\U_{\delta_1}\, dx &\leq & \int_\Omega \varphi_{\delta_2}^p\U_{\delta_1}\, dx =\int_{\O\setminus B(0,\rho)}\varphi_{\delta_2}^p\U_{\delta_1}\, dx +\int_{B(0,\rho)}\varphi_{\delta_2}^p\U_{\delta_1}\, dx\\
&\leq & |\varphi_{\delta_2}|^p_{p+1}\left(\int_{\O\setminus B(0,\rho)}\U_{\delta_1}^{p+1}\, dx\right)^{\frac{1}{p+1}}+ C\delta_2^{\frac{N+2}{2}}\int_{B(0,\rho)}\frac{1}{\left(1+\left|\frac{x}{\delta_1}\right|^2\right)^{\frac{N-2}{2}}}\, dx\\
&\leq & C_1\delta_2^{\frac{N+2}{2}}\delta_1^{\frac{N-2}{2}}+C_2\delta_2^{\frac{N+2}{2}}\delta_1^N\int_0^{\frac{\rho}{\delta_1}}\frac{r^{N-1}}{(1+r^2)^{\frac{N-2}{2}}}\, dr\\
&\leq & C_3\left[\delta_2^{\frac{N+2}{2}}\delta_1^{\frac{N-2}{2}}+\delta_2^{\frac{N+2}{2}}\delta_1^{N-2}\right]\leq C_3\left(\frac{\delta_2}{\delta_1}\right)^{\frac{N+2}{2}}\delta_1^{\frac{N}{2}}.
\end{eqnarray*}
Moreover, since $\int_{\O}\U_{\delta_j}^p\, dx =O(\delta_j^{\frac{N-2}{2}})$, we get
\begin{eqnarray*}
\int_{\O}\U_{\delta_2}^{p-1}\varphi_{\delta_2}\p\U_{\delta_1}\, dx& \leq &\|\varphi_{\delta_2}\|_{\infty} \int_{\O}\U_{\delta_2}^{p-1}\U_{\delta_1}\, dx \leq C \delta_2^{\frac{N-2}{2}}\int_{\O}\U_{\delta_2}^{p-1}\U_{\delta_1}\, dx\\
&\leq &C \delta_2^{\frac{N-2}{2}}\left(\int_{\O}\U_{\delta_2}^p\, dx\right)^{\frac{p-1}{p}}\left(\int_{\O}\U_{\delta_1}^p\, dx\right)^{\frac{1}{p}}\\
&\leq & C_1 \delta_2^{\frac{N^2+4N-12}{2(N+2)}}\delta_1^{\frac{(N-2)^2}{2(N+2)}}= C_1 \left(\frac{\delta_2}{\delta_1}\right)^{\frac{N-2}{2}} \left(\frac{\delta_2}{\delta_1}\right)^{\frac{2(N-2)}{N+2}} \delta_1^{{N-2}}\\
 &\leq&C_1 \left(\frac{\delta_2}{\delta_1}\right)^{\frac{N}{2}} \delta_1^{{N-2}}.
\end{eqnarray*}
Now let $\rho>0$ and we decompose the domain $\O$ as $\O=A_0\cup A_1\cup A_2$ where $A_0=\O\setminus B(0,\rho)$, $A_1=B(0,\rho)\setminus B(0,\sqrt{\delta_1\delta_2})$, $A_2=B(0,\sqrt{\delta_1\delta_2})$. Then we define $$L_j:=-\frac{1}{p+1}\int_{A_j}\left[|\p\U_{\delta_1}-\p\U_{\delta_2}|^{p+1}-\p\U_{\delta_1}^{p+1}-\p\U_{\delta_2}^{p+1}+(p+1)\p\U_{\delta_2}^p\p\U_{\delta_1}\right]\, dx$$
for $j=0,1,2$.\\
Now, by using Lemma \ref{lem0n2} and H\"older inequality, we see that
\begin{eqnarray*}
|L_0|&\leq &\frac{1}{p+1} \left[\int_{A_0}\left(|\p\U_{\delta_1}-\p\U_{\delta_2}|^{p+1}-\p\U_{\delta_1}^{p+1}\right)\, dx +\int_{A_0}\p\U_{\delta_2}^{p+1}\, dx + \int_{A_0}\p\U_{\delta_2}^p\p\U_{\delta_1}\, dx \right]\\
&\leq & C \left(  \int_{A_0}\p\U_{\delta_1}^p\p\U_{\delta_2}\, dx +\int_{A_0}\p\U_{\delta_2}^{p+1}\, dx + \int_{A_0}\p\U_{\delta_2}^p\p\U_{\delta_1}\, dx\right)\\
&\leq & C \left(\int_{A_0}\U_{\delta_1}^p\U_{\delta_2}\, dx+\int_{A_0}\U_{\delta_2}^{p+1}\, dx + \int_{A_0}\U_{\delta_2}^p\U_{\delta_1}\, dx \right)\\
&\leq & C\left(\int_{A_0}\U_{\delta_1}^{p+1}\, dx \right)^{\frac{p}{p+1}}\left(\int_{A_0}\U_{\delta_2}^{p+1}\, dx \right)^{\frac{1}{p+1}}+ C_1\int_{\frac{\rho}{\delta_2}}^{+\infty}\frac{r^{N-1}}{(1+r^2)^N}\, dr\\
&&+ C\left(\int_{A_0}\U_{\delta_2}^{p+1}\, dx\right)^{\frac{p}{p+1}}\left(\int_{A_0}\U_{\delta_1}^{p+1}\, dx\right)^{\frac{1}{p+1}}\\
&\leq & C_2\left(\int_{\frac{\rho}{\delta_1}}^{+\infty}\frac{r^{N-1}}{(1+r^2)^N}\, dr\right)^{\frac{p}{p+1}}\left(\int_{\frac{\rho}{\delta_2}}^{+\infty}\frac{r^{N-1}}{(1+r^2)^N}\, dr\right)^{\frac{1}{p+1}}+ C_3\delta_2^N\\
&&+C_2\left(\int_{\frac{\rho}{\delta_2}}^{+\infty}\frac{r^{N-1}}{(1+r^2)^N}\, dr\right)^{\frac{p}{p+1}}\left(\int_{\frac{\rho}{\delta_1}}^{+\infty}\frac{r^{N-1}}{(1+r^2)^N}\, dr\right)^{\frac{1}{p+1}}\\
&\leq & C_4\left(\delta_1^{\frac{N+2}{2}}\delta_2^{\frac{N-2}{2}}+\delta_2^N+\delta_2^{\frac{N+2}{2}}\delta_1^{\frac{N-2}{2}}\right)\leq C_5 \left(\frac{\delta_2}{\delta_1}\right)^{\frac{N-2}{2}}\delta_1^N.
\end{eqnarray*}
Now
\begin{eqnarray*}
L_1&=& -\frac{1}{p+1}\int_{A_1}\left[|\p\U_{\delta_1}-\p\U_{\delta_2}|^{p+1}-\p\U_{\delta_1}^{p+1}+(p+1)\p\U_{\delta_1}^p\p\U_{\delta_2}\right]\, dx\\
&&+\int_{A_1}\p\U_{\delta_1}^p\p\U_{\delta_2}\, dx-\int_{A_1}\p\U_{\delta_2}^p\p\U_{\delta_1}\, dx-\frac{1}{p+1}\int_{A_1}\p\U_{\delta_2}^{p+1}\, dx.
\end{eqnarray*}
Applying Lemma \ref{lem0n2} we get
\begin{eqnarray*}
&&\left|\int_{A_1}\left[|\p\U_{\delta_1}-\p\U_{\delta_2}|^{p+1}-\p\U_{\delta_1}^{p+1}+(p+1)\p\U_{\delta_1}^p\p\U_{\delta_2}\right]\, dx\right|\leq C \left(\int_{A_1}\p\U_{\delta_1}^{p-1}\p\U_{\delta_2}^{2}\, dx+\int_{A_1}\p\U_{\delta_2}^{p+1}\, dx\right)\\
&&\leq C_1 \left( \left(\frac{\delta_2}{\delta_1}\right)^2\int_{\sqrt{\frac{\delta_1}{\delta_2}}}^{\frac{\rho}{\delta_2}}\frac{1}{r^{N-3}}\, dr +  \int_{\sqrt{\frac{\delta_2}{\delta_1}}}^{\frac{\rho}{\delta_2}}\frac{r^{N-1}}{(1+r^2)^N} \ dr\right) \leq C_2\left(\frac{\delta_2}{\delta_1}\right)^{\frac{N}{2}}.
\end{eqnarray*}

Thanks to (\ref{proiezione}) and  Lemma \ref{lem0n2} we have

\begin{eqnarray*}
\int_{A_1}\p\U_{\delta_1}^p\p\U_{\delta_2}\, dx&=& \int_{A_1}\U_{\delta_1}^p \p\U_{\delta_2}\, dx + O\left(\int_{A_1}\U_{\delta_1}^{p-1}\varphi_{\delta_1}\p\U_{\delta_2}\, dx\right)+O\left(\int_{A_1}\varphi_{\delta_1}^p\p\U_{\delta_2}\, dx\right)\\
&=&\int_{A_1}\U_{\delta_1}^p \U_{\delta_2}\, dx + O\left(\int_{\O}\U_{\delta_1}^p\varphi_{\delta_2}\, dx\right)+ O\left(\int_{\O}\U_{\delta_1}^{p-1}\varphi_{\delta_1}\p\U_{\delta_2}\, dx\right)+\\
&&+O\left(\int_{\O}\varphi_{\delta_1}^p\p\U_{\delta_2}\, dx\right)
\end{eqnarray*}
By definition we have:

\begin{eqnarray*}
\int_{A_1}\U_{\delta_1}^p\U_{\delta_2}\, dx &=& \alpha_N^{p+1}\int_{A_1}\frac{\delta_1^{-\frac{N+2}{2}}}{\left(1+\left|\frac{x}{\delta_1}\right|^2\right)^{\frac{N+2}{2}}}
\frac{\delta_2^{-\frac{N-2}{2}}}{\left(1+\left|\frac{x}{\delta_2}\right|^2\right)^{\frac{N-2}{2}}}\, dx\\
&=&\alpha_N^{p+1}\delta_1^{-\frac{N+2}{2}+N}\delta_2^{-\frac{N-2}{2}}
\left(\frac{\delta_2}{\delta_1}\right)^{N-2}\int_{\sqrt{\frac{\delta_2}{\delta_1}}\leq |x|\leq \frac{\rho}{\delta_1}}\frac{1}{(1+|y|^2)^{\frac{N+2}{2}}}\frac{1}{|y|^{N-2}}\, dy+o\left(\left(\frac{\delta_2}{\delta_1}\right)^{\frac{N-2}{2}}\right)\\
&=& a_3\left(\frac{\delta_2}{\delta_1}\right)^{\frac{N-2}{2}}+o\left(\left(\frac{\delta_2}{\delta_1}\right)^{\frac{N-2}{2}}\right)
\end{eqnarray*}
Moreover by using \eqref{stimaproiezione} we get
\begin{equation*}
\int_{\O}\U_{\delta_1}^p\varphi_{\delta_2}\, dx\leq C \delta_2^{\frac{N-2}{2}}\int_{\O}\U_{\delta_1}^p\, dx \leq C_1 \delta_1^{\frac{N-2}{2}}\delta_2^{\frac{N-2}{2}}
\end{equation*}
and by using again \eqref{stimaproiezione} we have
\begin{eqnarray*}
\int_{\O}\U_{\delta_1}^{p-1}\varphi_{\delta_1}\p\U_{\delta_2}\, dx&\leq & \int_{\O}\U_{\delta_1}^{p-1}\varphi_{\delta_1}\U_{\delta_2}\, dx\\
&\leq & C \delta_1^{\frac{N-2}{2}}\left(\int_{\O}\U_{\delta_1}^{p+1}\, dx\right)^{\frac{p-1}{p+1}} \left(\int_{\O}\U_{\delta_2}^{\frac{p+1}{2}}\, dx\right)^{\frac{2}{p+1}}\leq C_1 \delta_1^{\frac{N-2}{2}}\delta_2^{\frac{N-2}{2}}.
\end{eqnarray*}
Finally
\begin{eqnarray*}
\int_{\O}\varphi_{\delta_1}^p \p\U_{\delta_2}\, dx &\leq & \int_{\O}\varphi_{\delta_1}^p\U_{\delta_2}\, dx = \int_{B(0,\rho)}\varphi_{\delta_1}^p\U_{\delta_2}\, dx+\int_{\O\setminus B(0,\rho)}\varphi_{\delta_1}^p\U_{\delta_2}\, dx\\
&\leq &C \delta_1^{\frac{N+2}{2}}\delta_2^{-\frac{N-2}{2}+N}\int_0^{\frac{\rho}{\delta_2}}
\frac{r^{N-1}}{(1+r^2)^{\frac{N-2}{2}}}\, dr+|\varphi_{\delta_1}|^p_{p+1}\left(\int_{\O\setminus B(0,\rho)}\U_{\delta_2}^{p+1}\, dx\right)^{\frac{1}{p+1}}\\
&\leq & C\delta_1^{\frac{N+2}{2}}\delta_2^{\frac{N+2}{2}}+C_1\delta_1^{\frac{N+2}{2}}\left(\int_{\frac{\rho}{\delta_2}}^{+\infty}\frac{r^{N-1}}{(1+r^2)^N}\, dr\right)^{\frac{N-2}{2N}}\\
&\leq & C_2 \delta_1^{\frac{N+2}{2}}\delta_2^{\frac{N-2}{2}}.
\end{eqnarray*}
At the end
\begin{eqnarray*}
\int_{A_1}\p\U_{\delta_2}^p\p\U_{\delta_1}\, dx &\leq & \alpha_N^{p+1} \delta_2^{-\frac{N+2}{2}+N}\delta_1^{-\frac{N-2}{2}}\int_{\sqrt{\frac{\delta_1}{\delta_2}}\leq |y|\leq \frac{\rho}{\delta_2}}\frac{1}{\left(1+\left|\frac{\delta_2}{\delta_1}y\right|^2\right)^{\frac{N-2}{2}}}\frac{1}{(1+|y|^2)^{\frac{N+2}{2}}}\, dy\\
&\leq & C_1\left(\frac{\delta_2}{\delta_1}\right)^{\frac{N-2}{2}}\int_{\sqrt{\frac{\delta_1}{\delta_2}}}^{\frac{\rho}{\delta_2}}\frac{r^{N-1}}{(1+r^2)^{\frac{N+2}{2}}}\, dr \leq  C_2\left(\frac{\delta_2}{\delta_1}\right)^{\frac{N}{2}}.
\end{eqnarray*}
Finally, thanks to Lemma \ref{lem0n2} we get that
\begin{eqnarray*}
|L_2|&\leq&\frac{1}{p+1}\left\{\left|\int_{A_2}\left[|\p\U_{\delta_1}-\p\U_{\delta_2}|^{p+1}-\p\U_{\delta_2}^{p+1}+(p+1)\p\U_{\delta_2}^p\p\U_{\delta_1}\right]\, dx\right|+\int_{A_2}\p\U_{\delta_1}^{p+1}\, dx \right\}\\
&\leq & C\left(\int_{A_2}\p\U_{\delta_2}^{p-1}\p\U_{\delta_1}^{2}\, dx+ \int_{A_2}\p\U_{\delta_1}^{p+1}\, dx\right) \leq C\left(\int_{A_2}\U_{\delta_2}^{p-1}\U_{\delta_1}^{2}\, dx+ \int_{A_2}\U_{\delta_1}^{p+1}\, dx\right)\\
&\leq & C_1\left(\left(\frac{\delta_2}{\delta_1}\right)^2\int_0^{\sqrt{\frac{\delta_2}{\delta_1}}}\frac{r^{N-5}}{(1+r^2)^{N-2}}\, dr+ \int_0^{\sqrt{\frac{\delta_2}{\delta_1}}}\frac{r^{N-1}}{(1+r^2)^N}\, dr \right)\\
&\leq& C_2\left(\left(\frac{\delta_2}{\delta_1}\right)^2\int_0^{\sqrt{\frac{\delta_2}{\delta_1}}}r^{N-5}\, dr+ \int_0^{\sqrt{\frac{\delta_2}{\delta_1}}}r^{N-1}\, dr\right) \leq C_2\left(\frac{\delta_2}{\delta_1}\right)^{\frac{N}{2}}.
\end{eqnarray*}
From Lemma \ref{lemma1} to Lemma \ref{lemma7} summing up all the terms we get that
\begin{eqnarray}\label{funzionalenellebubble}
\begin{array}{lll}
\displaystyle J_\epsilon(\p\U_{\delta_1}-\p\U_{\delta_2})&=&\displaystyle \frac{2}{N} S^{N/2} + a_1 \tau(0) \delta_1^{N-2} + a_1 \tau(0) \delta_2^{N-2} + O(\delta_1^{N-1}) + O(\delta_2^{N-1})\\[12pt]
&+&\displaystyle O\left( \left(\frac{\delta_2}{\delta_1}\right)^{\frac{N-2}{2}} \delta_1^{2}\right)+ O\left( \left(\frac{\delta_2}{\delta_1}\right)^{\frac{N}{2}}
\right) + a_3 \left(\frac{\delta_2}{\delta_1}\right)^{\frac{N-2}{2}}\\[12pt]
&-&a_2 \epsilon \delta_1^2 +  O\left( \epsilon \delta_1^{N-2}\right) - a_2 \epsilon \delta_2^2 +  O\left( \epsilon \delta_2^{N-2}\right),
\end{array}
\end{eqnarray}
where $a_1=\frac{1}{2} \alpha_{N}^{p+1} \int_{\R} \frac{1}{(1+|y|^2)^{\frac{N+2}{2}}} \ dy$, $a_2=\frac{1}{2}\alpha_N^2\int_{\R} \frac{1}{(1+|y|^2)^{N-2}} \ dy$, $a_3=\alpha_N^{p+1}  \int_{\R} \frac{1}{(1+|y|^2)^{\frac{N+2}{2}}} \ dy$.
 Recalling the choice of $\delta_j$, $j=1,2$ we get
 \begin{eqnarray}\label{funzionalenellebubble2}
\begin{array}{lll}
 \displaystyle J_\epsilon(\p\U_{\delta_1}-\p\U_{\delta_2})&=&\displaystyle \frac{2}{N} S^{N/2} +  a_1 \tau(0) d_1^{N-2} \epsilon^{\frac{N-2}{N-4}} + a_1 \tau(0) d_2^{N-2} \epsilon^{\frac{(3N-10)(N-2)}{(N-4)(N-6)}}  \\[12pt] &+&\displaystyle O\left(\epsilon^{\frac{N-1}{N-4}}\right)
+ O\left(\epsilon^{\frac{(3N-10)(N-1)}{(N-4)(N-6)}}\right)+O\left(\epsilon^{\frac{N+2}{N-6}}\right) +  O\left(\epsilon^{\frac{(N-2)N}{(N-4)(N-6)}}\right)\\[12pt]
&+& a_3 \left(\frac{d_2}{d_1}\right)^{\frac{N-2}{2}} \epsilon^{\frac{(N-2)^2}{(N-4)(N-6)}}-a_2d_1^2 \epsilon^{\frac{N-2}{N-4}} +O\left(\epsilon^{\frac{2(N-3)}{N-4}}\right)-a_2d_2^2\epsilon^{\frac{(N-2)^2}{(N-4)(N-6)}}\\[12pt]
&+&\displaystyle O\left(\epsilon^{\frac{2(2N^2-13N+22)}{(N-4)(N-6))}}\right)\\[12pt]
&=&\displaystyle \frac{2}{N} S^{N/2} +  [a_1 \tau(0) d_1^{N-2}-a_2d_1^2 ] \epsilon^{\frac{N-2}{N-4}} + O\left(\epsilon^{\frac{N-1}{N-4}}\right)\\[12pt]
&+&\displaystyle \left[a_3 \left(\frac{d_2}{d_1}\right)^{\frac{N-2}{2}} -a_2d_2^2 \right]\epsilon^{\frac{(N-2)^2}{(N-4)(N-6)}} +  O\left(\epsilon^{\frac{N+2}{N-6}}\right).
\end{array}
\end{eqnarray}
We point out that the term $O(\e^\frac{N-1}{N-4})$ depends only on $d_1$.
\end{itemize}

\end{proof}

\section{proof of Theorems \ref{principale} and \ref{principale1}}\label{teorema}

\begin{proof}[Proof of Theorem \ref{principale}]
Let us set $G_1(d_1):= a_1 \tau(0) d_1^{N-2} -a_2 d_1^2$, where $a_1$, $a_2$  are the positive constants appearing in Proposition \ref{ridotto} and $\tau(0)$ is the Robin's function of the domain $\Omega$ at the origin, so by definition it follows that $\tau(0)$ is positive. It's elementary to see that the function $G_1:\mathbb R^+ \rightarrow \mathbb R$ has a strictly local minimum point at $\bar d_1=\left(\frac{2a_2}{(N-2)a_1\tau(0)}\right)^{\frac{1}{N-4}}$.

Since $\bar d_1$ is a strictly local minimum for $G_1$, then, for any sufficiently small $\gamma>0$ there exists an open interval $I_{1,\sigma_1}$ such that $\overline I_{1,\sigma_1} \subset \mathbb R^+$, $I_{1,\sigma_1}$  has diameter $\sigma_1$, $\bar d_1 \in I_{1,\sigma_1}$ and  for all $d_1 \in \partial I_{1,\sigma_1}$
\begin{equation}\label{eq00propcrit}
G_1(d_1) \geq G_1(\bar d_1) + \gamma.
\end{equation}
Clearly as $\gamma \rightarrow 0$   we can choose $\sigma_1$ so that $\sigma_1 \rightarrow 0$.

We set $G_2(d_1,d_2):=a_3 \tau(0) \left(\frac{d_2}{d_1}\right)^{\frac{N-2}{2}} -a_2 d_2^2$, $G_2: \mathbb
R^2_+ \rightarrow \mathbb R$, where $a_3>0$ is the same constant appearing in Proposition \ref{ridotto}. If we fix $d_1=\bar d_1$ then $\hat G_2(d_2):=G( \bar d_1, d_2)$ has a strictly local minimum point at $\bar d_2:=\left(\frac{2a_2  \bar d_1^{\frac{N-2}{2}}}{a_3\tau(0) \frac{N-2}{2}}\right)^{\frac{2}{N-6}}$. As in the previous case
there exists an open interval $I_{2,\sigma_2}$ such that $\overline I_{2,\sigma_2} \subset \mathbb R^+$, $I_{2,\sigma_2}$  has diameter $\sigma_2$, $\bar d_2 \in I_{1,\sigma_1}$ and  for all $d_2 \in \partial I_{2,\sigma_2}$
\begin{equation}\label{eq01propcrit}
\hat G_2(d_2) \geq \hat G_2(\bar d_2) + \gamma.
\end{equation}
As $\gamma \rightarrow 0$   we can choose $\sigma_2$ so that $\sigma_2 \rightarrow 0$.

Let us set $K:= \overline{ I_{1,\sigma_1} \times I_{2,\sigma_2}}$ and let $\eta>0$ be small enough so that $K\subset ]\eta,\frac{1}{\eta}[\times]\eta,\frac{1}{\eta}[$. Thanks to Proposition \ref{auxsolving}, for all sufficiently small $\epsilon$, $\tilde J_\epsilon:\mathbb R_+^2 \rightarrow \mathbb R$ is defined  and it is of class $C^1$. By Weierstrass theorem we know there exists a global minimum point for $\tilde J_\epsilon$ in $K$. Let  $(d_{1,\epsilon},d_{2,\epsilon})$ be that point, we want to show that there exists $\epsilon_1$ such that, for all $\epsilon < \epsilon_1$, $(d_{1,\epsilon},d_{2,\epsilon})$ lies in the interior of $K$.

Assume by contradiction there exists a sequence $\epsilon_n\rightarrow 0$ such that for all $n \in \mathbb N$ $$ (d_{1,\epsilon_n},d_{2,\epsilon_n}) \in \partial K.$$ There are only two possibilities:
\begin{description}
\item[(a)] $d_{1,\epsilon_n} \in \partial I_{1,\sigma_1}$, $d_{2,\epsilon_n} \in I_{2,\sigma_2}$,
\item[(b)] $d_{1,\epsilon_n} \in  I_{1,\sigma_1}$, $d_{2,\epsilon_n} \in \partial I_{2,\sigma_2}$.
\end{description}

Thanks to (ii) of Proposition \ref{ridotto} we have the uniform expansion
\begin{equation}\label{eq11propcrit}
\tilde J_\epsilon (d_1, d_2) - \tilde J_\epsilon ( \bar d_1,d_2)=  \epsilon^{\theta_1}\left[ G_1(d_1)-G_1(\bar d_1)\right]
+  o\left( \epsilon^{\theta_1}\right).
\end{equation}
for all $\epsilon<\epsilon_0$, $(d_1,d_2) \in K$. We point out that we have incorporated the other high order terms in $o\left( \epsilon^{\theta_1}\right)$.
Thanks to (\ref{eq00propcrit}) and (\ref{eq11propcrit}), for all sufficiently small $\epsilon$ we have
\begin{equation} \label{eq20propcrit}
 \tilde J_\epsilon (d_1, d_2) - \tilde J_\epsilon ( \bar d_1,d_2) >0,
\end{equation}
for all $d_1 \in \partial I_{1,\sigma_1}$, for all $d_2 \in I_{2,\sigma_2}$. So for $n$ sufficiently large if (a) holds, since by definition $\tilde J_{\epsilon_n} (d_{1,\epsilon_n},d_{2,\epsilon_n})=\min_{K} \tilde J_{\epsilon_n}$, then
$$  \tilde J_{\epsilon_n} (d_{1,\epsilon_n},d_{2,\epsilon_n}) \leq \tilde J_{\epsilon_n} ( \bar d_1,d_{2,\epsilon_n}),$$
which contradicts (\ref{eq20propcrit}).
Assume (b). Thanks to (ii) of Proposition \ref{ridotto} (see also Remark \ref{remarkespfunzrid}) we have the uniform expansion
\begin{equation}\label{eq21propcrit}
\tilde J_\epsilon (d_1, d_2) - \tilde J_\epsilon ( d_1,\bar d_2)=  \epsilon^{\theta_2}\left[ G_2(d_1,d_2)-G_2( d_1,\bar d_2)\right] + o\left(\epsilon^{\theta_2}\right),
\end{equation}
for all $\epsilon \in (0,\epsilon_0)$, for all $(d_1,d_2) \in K$.

For $n$ sufficiently large so that $\epsilon_n < \epsilon_0$ we have
\begin{equation}\label{eq22propcrit}
\begin{array}{lll}
 \displaystyle \tilde J_{\epsilon_n}(d_{1,\epsilon_n},d_{2,\epsilon_n})-\tilde J_{\epsilon_n} ( d_{1,\epsilon_n},\bar d_2)&=&\displaystyle  \epsilon^{\theta_2}\left[ G_2(d_{1,\epsilon_n},d_{2,\epsilon_n})-G_2( d_{1,\epsilon_n},\bar d_2)\right] + o\left(\epsilon^{\theta_2}\right)\\[12pt]
&=&\displaystyle  \epsilon^{\theta_2}\left[ G_2(d_{1,\epsilon_n},d_{2,\epsilon_n})-G_2(\bar d_1,d_{2,\epsilon_n}) + G_2(\bar d_1,d_{2,\epsilon_n})-G_2(\bar d_1,\bar d_{2})\right. \\[10pt]
&&\displaystyle \ \  \ \ \ \ \left.  G_2(\bar d_1,\bar d_{2}) - G_2( d_{1,\epsilon_n},\bar d_2)\right] + o\left(\epsilon^{\theta_2}\right)\\[12pt]
&=&\displaystyle  \epsilon^{\theta_2}\left[ a_3 \tau(0) d_{2,\epsilon_n}^{\frac{N-2}{2}} \left(\frac{1}{d_{1,\epsilon_n}^{\frac{N-2}{2}}}-\frac{1}{\bar d_{1}^{\frac{N-2}{2}}}\right) + G_2(\bar d_1,d_{2,\epsilon_n})-G_2(\bar d_1,\bar d_{2})\right. \\[18pt]
&&\displaystyle  \ \ \ \ \left. + \ a_3 \tau(0) \bar d_{2}^{\frac{N-2}{2}} \left(\frac{1}{\bar d_{1}^{\frac{N-2}{2}}}-\frac{1}{d_{1,\epsilon_n}^{\frac{N-2}{2}}}\right)\right] + o\left(\epsilon_n^{\theta_2}\right)
\end{array}
\end{equation}

We observe now that, up to a subsequence, $d_{1,\epsilon_n} \rightarrow \bar d_1$ as $n\rightarrow + \infty$. This is a consequence of the uniform expansion given by  (ii) of Proposition \ref{ridotto}, in fact
\begin{equation}\label{eq31propcrit}
\tilde J_{\epsilon_n} (d_{1,\epsilon_n}, d_{2,\epsilon_n}) - \tilde J_{\epsilon_n} ( \bar d_1,\bar d_2)=  \epsilon_n^{\theta_1}\left[ G_1(d_{1,\epsilon_n})-G_1(\bar d_1)\right]
+  o\left( \epsilon_n^{\theta_1}\right).
\end{equation}
 Since $(d_{1,\epsilon_n}, d_{2,\epsilon_n}) $ is the minimum point we have $\tilde J_\epsilon (d_{1,\epsilon_n}, d_{2,\epsilon_n}) - \tilde J_\epsilon ( \bar d_1,\bar d_2) \leq 0$, hence, dividing (\ref{eq31propcrit}) by $\epsilon_n^{\theta_1}$, for all sufficiently large $n$ we get that
 $ G_1(d_{1,\epsilon_n})-G_1(\bar d_1) \leq - \frac{o\left( \epsilon_n^{\theta_1}\right)}{\epsilon_n^{\theta_1}}$. On the other side, since $\bar d_1$ is the minimum of $G_1$, we get that
$ G_1(d_{1,\epsilon_n})-G_1(\bar d_1) \geq 0$. So we have proved that
$$0\leq G_1(d_{1,\epsilon_n})-G_1(\bar d_1) \leq -\frac{o\left( \epsilon_n^{\theta_1}\right)}{\epsilon_n^{\theta_1}}, $$
and passing to the limit we deduce that $\lim_{n \rightarrow + \infty}  G_1(d_{1,\epsilon_n})=G_1(\bar d_1)$. Hence, up to a subsequence, since $\bar d_1$ is a strict local minimum, the only possibility is  $d_{1,\epsilon_n} \rightarrow \bar d_1$.

Since we are assuming (b), from (\ref{eq01propcrit}) we get that
$$ G_2(\bar d_1,d_{2,\epsilon_n})-G_2(\bar d_1,\bar d_{2}) \geq \gamma.$$
From this last inequality, (\ref{eq22propcrit}) and since $(d_{2,\epsilon_n})_n$ is bounded, then, choosing $\bar n$ sufficiently large so that $a_3 \tau(0) d_{2,\epsilon_n}^{\frac{N-2}{2}} \left|\frac{1}{\bar d_{1}^{\frac{N-2}{2}}}-\frac{1}{d_{1,\epsilon_n}^{\frac{N-2}{2}}}\right|$ and $a_3 \tau(0) \bar d_{2}^{\frac{N-2}{2}} \left|\frac{1}{\bar d_{1}^{\frac{N-2}{2}}}-\frac{1}{d_{1,\epsilon_n}^{\frac{N-2}{2}}}\right|$ are small enough, we deduce that
$$ \tilde J_{\epsilon_n}(d_{1,\epsilon_n},d_{2,\epsilon_n})-\tilde J_{\epsilon_n} ( d_{1,\epsilon_n},\bar d_2) > 0,$$
for all  $n> \bar n$. Since $(d_{1,\epsilon_n},d_{2,\epsilon_n})$ is the minimum point it also holds $$\tilde J_{\epsilon_n}(d_{1,\epsilon_n},d_{2,\epsilon_n})-\tilde J_{\epsilon_n} ( d_{1,\epsilon_n},\bar d_2) \leq 0,$$
and we get a contradiction.

To complete the proof we point out that, as observed before, up to a subsequence $d_{1,\epsilon} \rightarrow \bar d_1$ as $\epsilon \rightarrow 0$. With a similar argument we prove that $d_{2,\epsilon} \rightarrow \bar d_2$. In fact, from the same argument of (\ref{eq22propcrit}), since $d_{1,\epsilon} \rightarrow \bar d_1$ and $(d_{2,\epsilon})_\epsilon$ is bounded,  we have
\begin{equation} \label{eq41propcrit}
\begin{array}{lll}
 \displaystyle 0\geq \frac{\tilde J_{\epsilon}(d_{1,\epsilon},d_{2,\epsilon})-\tilde J_{\epsilon} ( d_{1,\epsilon},\bar d_2)}{ \epsilon^{\theta_2}}&=&\displaystyle  G_2(d_{1,\epsilon},d_{2,\epsilon})-G_2( d_{1,\epsilon},\bar d_2) + \frac{o\left(\epsilon^{\theta_2}\right)}{ \epsilon^{\theta_2}}\\[12pt]
&=&\displaystyle   a_3 \tau(0) d_{2,\epsilon}^{\frac{N-2}{2}} \left(\frac{1}{d_{1,\epsilon}^{\frac{N-2}{2}}}-\frac{1}{\bar d_{1}^{\frac{N-2}{2}}}\right) + G_2(\bar d_1,d_{2,\epsilon})-G_2(\bar d_1,\bar d_{2}) \\[18pt]
&&\displaystyle + \ a_3 \tau(0) \bar d_{2}^{\frac{N-2}{2}} \left(\frac{1}{\bar d_{1}^{\frac{N-2}{2}}}-\frac{1}{d_{1,\epsilon}^{\frac{N-2}{2}}}\right) + \frac{o\left(\epsilon^{\theta_2}\right)}{ \epsilon^{\theta_2}}\\[18pt]
&=& \displaystyle o(1) + G_2(\bar d_1,d_{2,\epsilon})-G_2(\bar d_1,\bar d_{2}).
\end{array}
\end{equation}
Since $\bar d_2$ is a local maximum point for $d_2 \rightarrow \hat G_2(d_2)$ we have $G_2(\bar d_1,d_{2,\epsilon})-G_2(\bar d_1,\bar d_{2})\geq 0$ and so from (\ref{eq41propcrit}) we get that
$$0 \leq G_2(\bar d_1,d_{2,\epsilon})-G_2(\bar d_1,\bar d_{2}) \leq - o(1). $$
Passing to the limit as $\epsilon \rightarrow 0$ we deduce that $ \hat G_2(d_{2,\epsilon}) \rightarrow \hat G_2(\bar d_{2}) $.  Hence, up to a subsequence, since $\bar d_2$ is a strict local minimum, the only possibility is  $d_{2,\epsilon} \rightarrow \bar d_2$. \\ Hence by (i) of Proposition \ref{ridotto} we have that $V_{\bar{d}_\e}+\bar\phi_1+\bar\phi_2$ is a solution  of \eqref{BN} and the proof is complete.
\end{proof}
We are ready also to prove Theorem \ref{principale1}. We reason as in \cite{Pistoia}.
\begin{proof}[Proof of Theorem \ref{principale1}]
Let $u_\e$ be a solution of \eqref{BN} as in Theorem \ref{principale} and assume that $\Phi_\e\to 0$ uniformly in compact subsets of $\Omega$. We set
\begin{eqnarray*}
\tilde u_\e(x)&:=&\left(\frac{d_{1\e} \e^{\frac{1}{N-4}}}{d_{1\e}^2\e^{\frac{2}{N-4}}+|x|^2}\right)^{\frac{N-2}{2}}-
\left(\frac{d_{1\e} \e^{\frac{3N-10}{(N-4)(N-6)}}}{d_{1\e}^2\e^{2\frac{3N-10}{(N-4)(N-6)}}+|x|^2}\right)^{\frac{N-2}{2}}\\
&=& \left(\frac{1}{d_{1\e}\e^{\frac{1}{N-4}}+d_{1\e}^{-1}\e^{-\frac{1}{N-4}}|x|^2}\right)^{\frac{N-2}{2}}-\left(\frac{1}{d_{2\e}
\e^{\frac{3N-10}{(N-4)(N-6)}}+d_{2\e}^{-1}\e^{-\frac{3N-10}{(N-4)(N-6)}}|x|^2}\right)^{\frac{N-2}{2}}
\end{eqnarray*}
Then, by Theorem \ref{principale} and by using the assumption on the remainder term $\Phi_\e$ we get
\begin{equation}\label{pallino}
u_\e(x)=\alpha_N \tilde u_\e(x)(1+o(1)), \qquad x\in\O,
\end{equation}
where $o(1)\rightarrow 0$ uniformly on compact subsets of $\O$.\\ We consider the spheres
$$\mathcal S_\e^1:=\{x\in \R\,\,:\,\, |x|=\e^{\frac{1}{N-4}}\}$$ and $$\mathcal S_\e^2:=\{x\in\R\,\,:\,\, |x|=\e^{\frac{3N-10}{(N-4)(N-6)}}\}.$$ We may fix a compact subset $K\subset \O$ such that $\mathcal S_\e^j\subset K $, $j=1,2$ and $\e>0$ sufficiently small.\\ For $x\in \mathcal S_\e^1$ we get
\begin{eqnarray*}
\tilde u_\e(x)&=& \left(\frac{1}{d_{1\e}\e^{\frac{1}{N-4}}+d_{1\e}^{-1}\e^{\frac{1}{N-4}}}\right)^{\frac{N-2}{2}}-
\left(\frac{1}{d_{2\e}\e^{\frac{3N-10}{(N-4)(N-6)}}+d_{2\e}^{-1}\e^{-\frac{N+2}{(N-4)(N-6)}}}\right)^{\frac{N-2}{2}}\\
&=&\e^{-\frac{N-2}{2(N-4)}}\left[\left(\frac{1}{d_{1\e}+d_{1\e}^{-1}}\right)^{\frac{N-2}{2}}-
\left(\frac{1}{d_{2\e}\e^{\frac{2(N-2)}{(N-4)(N-6)}}+d_{2\e}^{-1}\e^{-\frac{8}{(N-4)(N-6)}}}\right)^{\frac{N-2}{2}}\right]\\
&=&\e^{-\frac{N-2}{2(N-4)}}\left[\left(\frac{1}{d_{1\e}+d_{1\e}^{-1}}\right)^{\frac{N-2}{2}}+o(1)\right]
\end{eqnarray*}
as $\e\rightarrow 0$. Hence $\tilde u_\e>0$ on $\mathcal S_\e^1$ for $\e$ small.\\ Analogously if $x\in\mathcal S_\e^2$ then
$$\tilde u_\e(x)=-\e^{-\frac{(3N-10)(N-2)}{2(N-4)(N-6)}}\left[\left(\frac{1}{d_{2\e}+d_{2\e}^{-1}}\right)^{\frac{N-2}{2}}+o(1)\right]$$
as $\e\rightarrow 0$ and hence $\tilde u_\e<0$ on $\mathcal S_\e^2$ for $\e$ small.\\ Since \eqref{pallino} holds, this implies that $u_\e>0$ on $\mathcal S_\e^1$ and $u_\e<0$ on $\mathcal S_\e^2$ for $\e$ small.\\ Then $u_\e$ has at least two nodal domains $\O_1, \O_2$ such that $\O_j$ contains the sphere $\mathcal S_\e^j$, $j=1,2$.\\ Next we show that $u_\e$ has not more than two nodal domains for $\e$ small.\\ We remark that by (ii) of Proposition \ref{ridotto} and by Lemmas \ref{lem1exp1}, \ref{lem2exp1} it follows that
\begin{equation}\label{croce}
J_\e(u_\e)\rightarrow \frac{2}{N}S^{\frac{N}{2}},\qquad \mbox{as}\,\, \e\rightarrow 0
\end{equation}
where $J_\e$ is defined in \eqref{funzionale} and $S$ is the best Sobolev constant for the embedding of $H^1_0(\O)$ into $L^{p+1}(\O)$, namely $$S:=\inf_{u\in H^1_0(\O)\setminus\{0\}}\frac{\int_\O |\nabla u|^2\, dx}{\left(\int_\O |u|^{p+1}\, dx\right)^{\frac{2}{p+1}}}.$$ We set $c_\e:=\inf_{\mathcal N_\e} J_\e$, where $\mathcal N_\e$ is the Nehari manifold, which is defined by
$$\mathcal N_\e:=\left\{ u\in H^1_0(\O)\,\,:\,\, \int_\O |\nabla u|^2\, dx =\int_\O |u|^{p+1}\, dx +\e\int_\O u^2\, dx\right\}.$$ It is easy to see that $c_\e \rightarrow c_0=\frac{1}{N}S^{\frac{N}{2}}$ as $\e\rightarrow 0$ and therefore, by \eqref{croce}, we get that
\begin{equation}\label{striscia}
J_\e(u_\e)<3 c_\e
\end{equation}
for $\e$ small enough.\\ We now suppose by contradiction that $u_\e$ has at least $3$ pairwise different nodal domains $\O_1, \O_2, \O_3$.\\ Let $\chi_i$ be the characteristic function corresponding to the sets $\Omega_i$.\\ Then $u_\e \chi_i\in H^1_0(\O)$ (see \cite{Muller}). Moreover
\begin{eqnarray*}
\int_\O |\nabla(u_\e \chi_i)|^2\, dx &=& \int_\O \nabla u_\e \nabla (u_\e \chi_i)=-\int_\O \Delta u_\e (u_\e\chi_i)\, dx\\
&=& \int_\O |u_\e|^p (u_\e\chi_i)\, dx +\e \int_\O u_\e \cdot u_\e\chi_i\, dx \\
&=& \int_\O |u_\e\chi_i|^{p+1}\, dx +\e \int_\O (u_\e \chi_i)^2\, dx
\end{eqnarray*}
so that $u_\e\chi_i \in \mathcal N_\e$. Since also $u_\e\in\mathcal N_\e$ we obtain
\begin{eqnarray*}
J_\e(u_\e)&=& \left(\frac 12 -\frac{1}{p+1}\right)\int_\O |u_\e|^{p+1}\, dx \\
&\geq & \left(\frac 12 -\frac{1}{p+1}\right)\sum_{i=1}^3 \int_\O |u_\e \chi_i|^{p+1}\, dx\\
&=& \sum_{i=1}^3 J_\e(\chi_i u_\e)\geq 3 c_\e
\end{eqnarray*}
contrary to \eqref{striscia}. The contradiction shows that $u_\e$ has at most two nodal domains for $\e$ small.\\ This completes the proof.
\end{proof}

\end{document}